\documentclass[11pt]{article}
\usepackage[utf8]{inputenc}
\usepackage[T1]{fontenc}

\usepackage{enumerate}
\usepackage{url}
\usepackage{epsf}
\usepackage{epsfig}
\usepackage{graphicx}

\usepackage{amsmath,amsfonts,amssymb,amsthm}
\usepackage{wasysym}
\usepackage[retainorgcmds]{IEEEtrantools}

\usepackage{stmaryrd}

\usepackage{ gensymb }
\usepackage{ esint }

\usepackage{hyperref}
\def\pathPic{Pics}                         

\def\cE{{\cal E}}
\def\cF{{\cal F}}

\def\cL{{\cal L}}

\def\cO{{\cal O}}

\def\cT{{\cal T}}
\def\cU{{\cal U}}
\def\cV{{\cal V}}

\def\cZ{{\cal Z}}

\def\mP{{\mathbb P}}

\def\mR{{\mathbb R}}

\def\mZ{{\mathbb Z}}

\def\R{\mR}

\def\Z{\mZ}

\def\ve{\varepsilon}
\def\vp{\varphi}

\DeclareMathOperator\diam{diam}

\def\<{\langle}
\def\>{\rangle}

\DeclareMathOperator\interp{{\rm I}}

\def\sm{\setminus}

\newtheorem{Theorem}{Theorem}[section]
\newtheorem{Remark}[Theorem]{Remark}
\newtheorem{Lemma}[Theorem]{Lemma}
\newtheorem{Proposition}[Theorem]{Proposition}
\newtheorem{Corollary}[Theorem]{Corollary}
\newtheorem{Definition}[Theorem]{Definition}

\newtheorem*{Theorem*}{Theorem}
\newtheorem*{Remark*}{Remark}
\newtheorem*{Lemma*}{Lemma}
\newtheorem*{Proposition*}{Proposition}
\newtheorem*{Corollary*}{Corollary}
\newtheorem*{Definition*}{Definition}
\newtheorem*{Notation*}{Notation}
\newtheorem*{Assumption*}{Assumption}

\DeclareMathOperator\Cone{Cone}
\DeclareMathOperator\Conv{Conv}
\DeclareMathOperator\DConv{DConv}
\DeclareMathOperator\supp{supp}

\DeclareMathOperator\Hull{Hull}
\DeclareMathOperator\Anc{Anc}
\DeclareMathOperator\GradConv{GradConv}

\renewcommand\subset{\subseteq}
\renewcommand\supset{\supseteq}

\makeatletter
\newcommand{\oset}[2]{%
  {\mathop{#2}\limits^{\vbox to -.5\ex@{\kern-\tw@\ex@
   \hbox{\scriptsize #1}\vss}}}}
\makeatother

\DeclareMathOperator\Cost{Cost}

\DeclareMathOperator\rank{rank}
\DeclareMathOperator\Hessian{Hessian}

\usepackage{algpseudocode}
\usepackage{algorithm}

\def\OmegaC{\Omega}
\def\OmegaD{X}
\def\uC{U}
\def\Candidates{\hat \cV}

\def\FixedGrid{\Z^2}

\def\interior{\oset {$\circ$}}

\author{Jean-Marie Mirebeau%
\footnote{CNRS, University Paris Dauphine, UMR 7534, Laboratory CEREMADE, Paris, France.}
}

\begin{document}
\title{
Adaptive, Anisotropic and Hierarchical\\
cones of Discrete Convex functions%
\thanks{
This work was partly supported by ANR grant NS-LBR ANR-13-JS01-0003-01. 
}
}
\maketitle
\date{}

\begin{abstract}
We introduce a new class of adaptive methods for optimization problems posed on the cone of convex functions. 
Among the various mathematical problems which posses such a formulation, the Monopolist problem \cite{Rochet:1998uj,Ekeland:2010tl} arising in economics is our main motivation.

Consider a two dimensional domain $\OmegaC$, sampled on a grid $\OmegaD$ of $N$ points. 
We show that the cone $\Conv(\OmegaD)$ of restrictions to $\OmegaD$ of convex functions on $\OmegaC$ is typically characterized by $\approx N^2$ linear inequalities; a direct computational use of this description therefore has a prohibitive complexity.
We thus introduce a hierarchy of sub-cones $\Conv(\cV)$ of $\Conv(\OmegaD)$, associated to  stencils $\cV$ which can be adaptively, locally, and anisotropically refined. 
We show, using the arithmetic structure of the grid, that the trace $U_{|X}$ of any convex function $U$ on $\Omega$ is contained in a cone $\Conv(\cV)$ defined by only $\cO( N \ln^2 N)$ linear constraints, in average over grid orientations. 

Numerical experiments for the Monopolist problem, based on adaptive stencil refinement strategies, show that the proposed method offers an unrivaled accuracy/complexity trade-off in comparison with existing methods.
We also obtain, as a side product of our theory, a new average complexity result on edge flipping based mesh generation. \end{abstract}

\begin{figure}
\centering
\includegraphics[width=3.8cm]{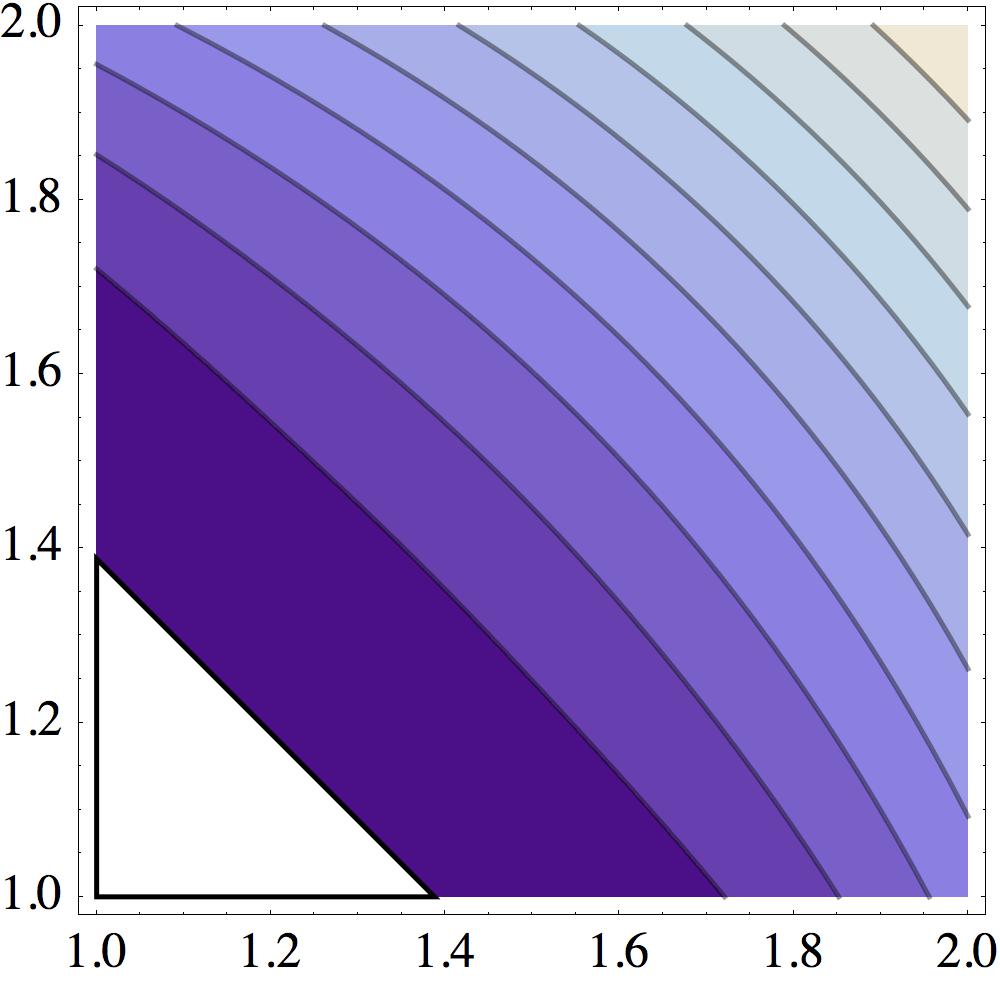}
\includegraphics[width=3.8cm]{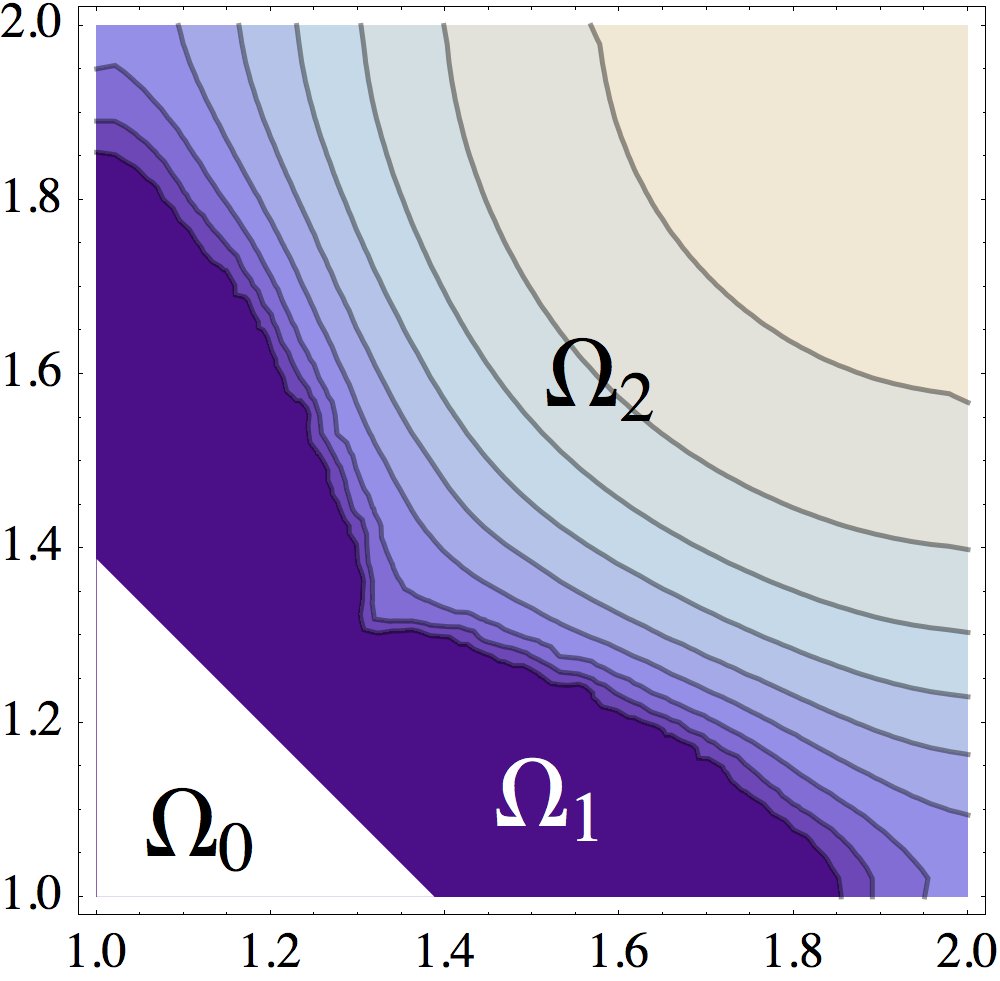}
\includegraphics[width=3.8cm]{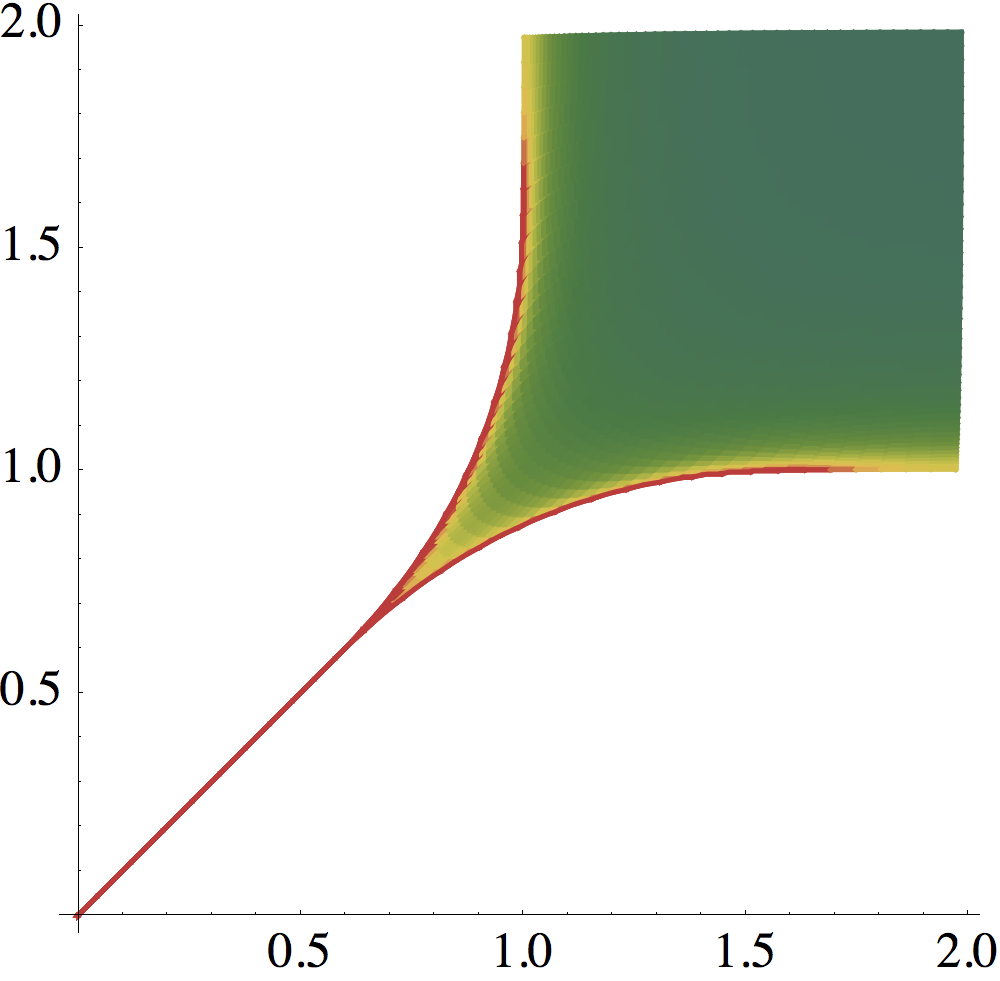}
\includegraphics[width=3.8cm]{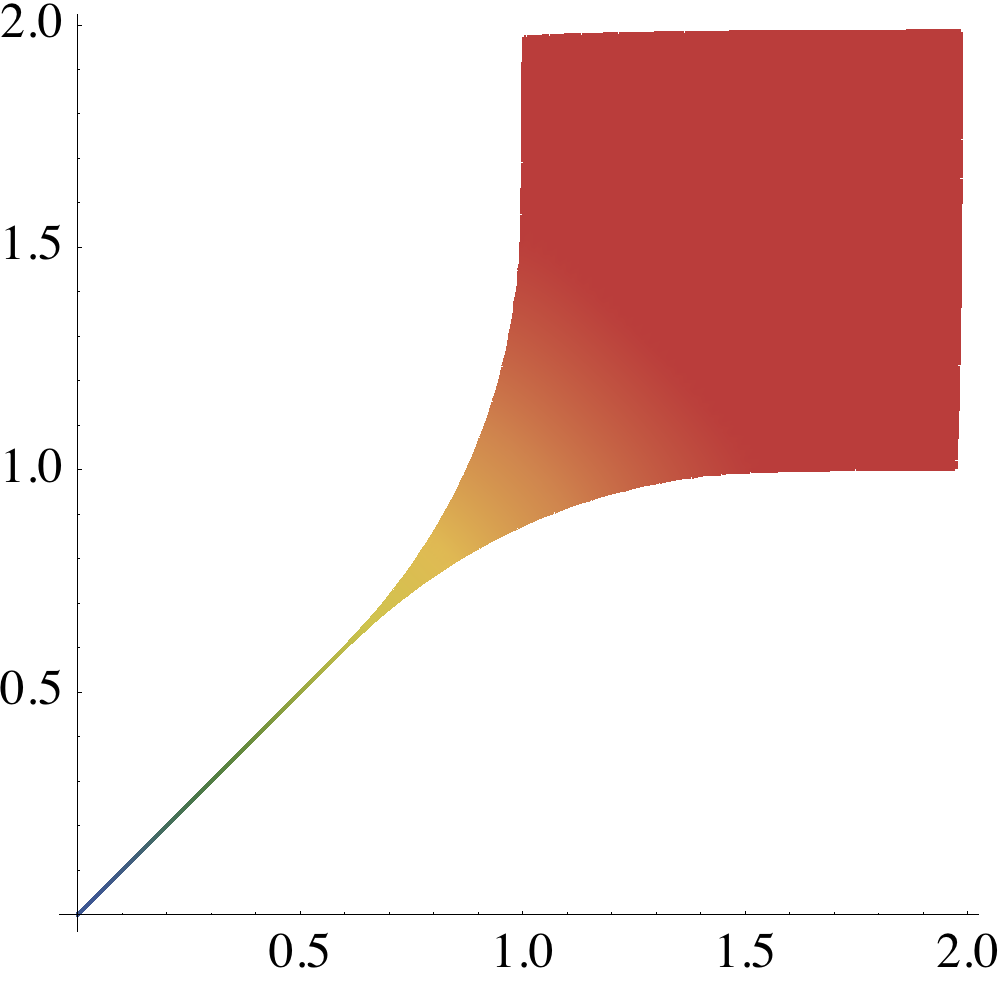}
\caption{
Numerical approximation $U$ of the solution of the classical Monopolist's problem \eqref{eq:PrincipalAgent}, computed on a $50 \times 50$ grid. Left: level sets of $U$, with $U=0$ in white. Center left: level sets of $\det(\Hessian U)$ (with again $U=0$ in white); note the degenerate region $\Omega_1$ where $\det(\Hessian U)=0$. Center right: distribution of products sold by the monopolist. Right: profit margin of the monopolist for each type of product (margins are low on the one dimensional part of the product line, at the bottom left). Color scales on Figure \ref{fig:PARotated}.
}
\label{fig:Monopolist}
\end{figure}

A number of mathematical problems can be formulated as the optimization of a convex functional over the \emph{cone of convex functions} on a domain $\OmegaC$ (here compact and two dimensional):
\begin{equation*}
\Conv(\OmegaC):=\{ \uC : \OmegaC \to \R; \, \uC \text{ is convex}\}.
\end{equation*}
This includes optimal transport, 
as well as various geometrical conjectures such as Newton's problem  \cite{LachandRobert:2005bi,MERIGOT:tr}. 
We choose for concreteness to emphasize an economic application:
the Monopolist (or Principal Agent) problem \cite{Rochet:1998uj}, in which the objective is to design an optimal product line, and an optimal pricing catalog, so as to maximize profit in a captive market. The following minimal instance is numerically studied in \cite{Aguilera:2008uq,Ekeland:2010tl,Oberman:2011wy} and on Figure \ref{fig:Monopolist}. With $\OmegaC = [1,2]^2$
\begin{equation}
\label{eq:PrincipalAgent}
\min \left\{ \int_{\OmegaC} \left( \frac 1 2 \|\nabla \uC(z)\|^2 -\<\nabla \uC(z),z\>+ \uC(z)\, \right) dz; \,  \uC \in \Conv(\OmegaC), \, \uC \geq 0\right\}. 
\end{equation}
We refer to the numerical section \S \ref{sec:Numerics}, and to \cite{Rochet:1998uj} for the economic model details; let us only say here that the Monopolist's optimal product line is $\{\nabla U(z); \, z \in \Omega\}$, and that the optimal prices are given by the Legendre-Fenchel dual of $U$. 
Consider the following three regions, defined for $k\in \{0,1,2\}$ (implicitly excluding points $z \in \Omega$ close to which $U$ is not smooth)
\begin{equation}
\label{def:OmegaK}
\Omega_k := \{z \in \Omega; \, \rank(\Hessian U(z) ) = k\}.
\end{equation}
Strong empirical evidence suggests that these three regions have a non-empty interior, although no qualitative mathematical theory has yet been developed for these problems. 
The optimal product line observed numerically, Figure \ref{fig:Monopolist}, confirms a qualitative (and conjectural) prediction of the economic model \cite{Rochet:1998uj} called ``bunching'': low-end products are less diverse than high-end ones, down to the topological sense. 
(The monopolist willingly limits the variety of cheap products, because they may compete with the more expensive ones, on which he has a higher margin.)

We aim to address numerically optimization problems posed on the cone of convex functions, through numerical schemes which preserve the rich qualitative properties of their solutions, and have a moderate computational cost. In order to put in light the specificity of our approach, we review 
the existing numerical methods for these problems, which fall in the following categories. 
We denote 
by $X$ a grid sampling of the domain $\Omega$, and by $\Conv(X)$ the cone of discrete (restrictions of) convex functions
\begin{equation}
\label{def:ConvX}
\Conv(X) := \{ U_{|X}; \, U \in \Conv(\Omega)\}.
\end{equation}

\begin{itemize}
\item (Interior finite element methods)
For any triangulation $\cT$ of $X$, consider the cone 
\begin{equation*}
\Conv(\cT) := \{ u : X \to \R; \, \interp_\cT u \in \Conv(\Omega)\}.
\end{equation*}
A natural but \emph{invalid} numerical method for \eqref{eq:PrincipalAgent} is to fix a-priori a family $(\cT_h)_{h>0}$ of regular triangulations of $\Omega$, where $h>0$ denotes mesh scale, and to optimize the functional of interest over the associated cones.
Indeed, the union of the cones $\Conv(\cT_h)$ is \emph{not} dense in $\Conv(\Omega)$, see \cite{Chone:2001fa}. Let us also mention that for a given generic $u \in \Conv(X)$, there exists \emph{only one} triangulation $\cT$ of $X$ such that $u \in \Conv(\cT)$, see \S \ref{sec:DelaunayIntro}.
\item (Global constraints methods) 
The functional of interest, suitably discretized, is minimized over the cone $\Conv(X)$ of discrete convex functions \cite{Carlier:2001tq},
or alternatively \cite{Ekeland:2010tl} on the augmented cone 
\begin{equation}
\label{def:GradConvX}
\GradConv(X) := \{ (U_{|X}, \nabla U_{|X}); \, U \in \Conv(\Omega)\},
\end{equation}
in which we refer by $\nabla U$ to arbitrary elements of the subgradient of the convex map $U$.

Both $\Conv(X)$ and $\GradConv(X)$ are characterized by a family of long range linear inequalities, with domain wide supports, and of cardinality growing \emph{quadratically} with $N := \#(X)$, see \S \ref{sec:Characterization} and \cite{Ekeland:2010tl}. Despite rather general convergence results, these two methods are impractical due to their expensive numerical cost, in terms of both computation time and memory.

\item (Local constraints methods) Another cone $\Conv'(X)$ is introduced, usually satisfying neither $\Conv(X) \subset \Conv'(X)$ nor $\Conv'(X) \subset \Conv(X)$, but typically characterized by relatively few  constraints, with short range supports. 
Obermann et al.\ \cite{Oberman:2011wi,Oberman:2011wy} use $\cO(N)$ linear constraints, with $N := \#(X)$. Merigot et al.\ \cite{MERIGOT:tr} use slightly more linear constraints, but provide an efficient optimization algorithm based on proximal operators. Aguilera et al.\ \cite{Aguilera:2008uq} consider $\cO(N)$ constraints of semi-definite type.

Some of these methods benefit from convergence guarantees \cite{MERIGOT:tr,Aguilera:2008uq} as $N \to \infty$. 
Our numerical experiments with \cite{Aguilera:2008uq,Oberman:2011wi,Oberman:2011wy} show however that they suffer from accuracy issues,  see \S \ref{sec:Comparison} and Figure \ref{fig:BadHessian}, which limits the usability of their results.

\item (Geometric methods)
A polygonal convex set can be described as the convex hull of a finite set of points, or as an intersection of half-spaces. Geometric methods approximate a convex function $U$ by representing its epigraph $\{ (z,t); \, z \in \Omega, \, t \geq U(z)\}$ under one of these forms. 
Energy minimization is done by adjusting the points position, or the coefficients of the affine forms defining the half-spaces, see \cite{Wachsmuth:2013ta,LachandRobert:2005bi}. 

The main drawback of these methods lies in the optimization procedure, which is quite non-standard. Indeed
the discretized functional is generally non-convex, and the polygonal structure of the represented convex set changes topology during the optimization.
\end{itemize}
 
We propose an implementation of the constraint of convexity via a limited (typically quasi-linear) number of linear inequalities, featuring both short range and domain wide supports, which are selected locally and anisotropically in an adaptation loop using a-posteriori analysis of solutions to intermediate problems. 
Our approach combines the accuracy of global constraint methods, with the limited cost of local constraint ones, see \S \ref{sec:Comparison}.
It is based on a family of sub-cones 
\begin{equation*}
\Conv(\cV) \subset \Conv(X),
\end{equation*}
each defined by some linear inequalities associated to a family $\cV$ of \emph{stencils}, see Definition \ref{def:Cones}. These stencils are the data  $\cV = (\cV(x))_{x \in X}$ of a collection of offsets $e \in \cV(x)$ pointing to selected neighbors $x+e$ of any point $x\in X$, and satisfying minor structure requirements, see Definition \ref{def:Stencil}. 
The cones satisfy the hierarchy property $\Conv(\cV \cap \cV') = \Conv(\cV) \cap \Conv(\cV')$, see Theorem \ref{th:Hierarchy}. 
Most elements of $\Conv(X)$ belong to a cone $\Conv(\cV)$ defined by only $\cO(N \ln^2 N)$ linear inequalities, in a sense made precise by Theorem \ref{th:Avg}. Regarding both stencils and triangulations as directed graphs on $X$, we show in Theorem \ref{th:Decomp} (under a minor technical condition) that the cone $\Conv(\cV)$ is the union of the cones $\Conv(\cT)$ associated to triangulations $\cT$ included in $\cV$. Our hierarchy of cones has similarities, but also striking differences as discussed in conclusion, with the other multiscale constructions (wavelets, adaptive finite elements) used in numerical analysis. 

The minimizer $u \in \Conv(X)$ of a given convex energy $\cE$ can be obtained without ever listing the inequalities defining $\Conv(X)$ (which would often not fit into computer memory for the problem sizes of interest), but only solving a small sequence of optimization problems over sub-cones $\Conv(\cV_i)$ associated to stencils $\cV_i$, designed through adaptive refinement strategies. 
Our numerical experiments give, we believe, unprecedented numerical insight on the qualitative behavior of the monopolist problem and its variants. Thanks to the adaptivity of our scheme, this accuracy is not at the expense of computation time or memory usage. See \S \ref{sec:Numerics}.

\section{Main results} 

The constructions and results developed in this paper apply to an arbitrary convex and compact domain $\Omega\subset \R^2$, discretized on an orthogonal grid of the form:
\begin{equation}
\label{def:ParametrizedDiscreteDomain}
\OmegaC \cap h R_\theta (\xi+ \Z^2),
\end{equation}
where $h>0$ is a scale parameter, $R_\theta$ is the rotation of angle $\theta\in \R$, and $\xi\in \R^2$ is an offset. The latter two parameters are used in our main approximation result Theorem \ref{th:Avg}, heuristically 
to eliminate by averaging the influence of rare unfavorable cases in which the approximated convex function hessian is degenerate in a direction close to the grid axes.
For simplicity, and up to a linear change of coordinates, we assume unless otherwise mentionned that these parameters take their canonical values: 
\begin{equation*}
\OmegaD := \OmegaC \cap \Z^2.
\end{equation*}

The choice of a grid discretization provides arithmetic tools that would not be available for an unstructured point set. 

\begin{Definition}
\begin{enumerate}
\item An element $e=(\alpha,\beta) \in \Z^2$ is called irreducible iff $\gcd(\alpha,\beta)=1$. 
\item A basis of $\Z^2$ is a pair $(f,g) \in (\Z^2)^2$ such that $|\det(f,g)| = 1$. A basis $(f,g)$ of $\Z^2$ is direct iff $\det(f,g) = 1$, and acute iff $\<f,g\> \geq 0$.
\end{enumerate}
\end{Definition} 
Considering special (non-canonical) bases of $\Z^d$ is relevant when discretizing anisotropic partial differential equations on grids, such as anisotropic diffusion \cite{Fehrenbach:2013ut}, or anisotropic eikonal equations \cite{Mirebeau2012anisotropic}.
In this paper, and in particular in the next proposition, we rely on a specific two dimensional structure called the Stern-Brocot tree \cite{Graham:1994uv}, also used in numerical analysis for anisotropic diffusion \cite{Bonnans:2004ud}, and eikonal equations of Finsler type \cite{Mirebeau:2012wm}.
\begin{Proposition}
\label{prop:Parents}
The application 
$
(f,g) \mapsto e := f+g
$
defines a bijection between direct acute bases $(f,g)$ of $\Z^2$, and irreducible elements $e\in \Z^2$ such that $\|e\| > 1$. The elements $f,g$ are called the parents of $e$. (Unit vectors have no parents.)
\end{Proposition}

\begin{proof}
Existence, for a given irreducible $e$ with $\|e\|>1$, of the direct acute basis $(f,g)$ such that $e=f+g$.
We assume without loss of generality that $e=(\alpha,\beta)$ has non-negative coordinates. Since $\gcd(\alpha,\beta)=1$ and $\|e\|>1$ we obtain that $\alpha \geq 1$ and $\beta \geq 1$. Classical results on the Stern-Brocot tree \cite{Graham:1994uv} state that the irreducible positive fraction $\alpha/\beta$ can be written as the \emph{mediant} $(\alpha'+\alpha'')/(\beta'+\beta'')$ of two irreducible fractions $\alpha'/\beta'$, $\alpha''/\beta''$ (possibly equal to $0$ or $+\infty$), with $\alpha',\beta',\alpha'',\beta'' \in \Z_+$ and $\alpha' \beta''-\beta'\alpha''=1$. Setting $f = (\alpha',\beta')$ and $g=(\alpha'',\beta'')$ concludes the proof.

Uniqueness. Assume that $e=f+g=f'+g'$, where $(f,g)$, $(f',g')$ are direct acute bases of $\Z^2$. One has $\det(f,e)=\det(f,f+g) = 1$, and likewise $\det(f',e)=1$. Hence $\det(f-f',e)=0$, and therefore $f' = f+k e$ for some scalar $k$, which is an integer since $e$ is irreducible.  Subtracting we obtain $g' = e-f' = g-k e$, and therefore
\begin{align*}
\<f',g'\> = \<f+k e, g-k e\> &= \<(k+1) f+k g, -k f +(1-k) g\>\\
 &= -k (k+1) \|f\|^2 - k (k-1) \|g\|^2 + \<f,g\> (1-2 k^2). 
\end{align*}
This expression is negative unless the integer $k$ is zero, hence $f=f'$, and $g=g'$.
\end{proof}

\subsection{Characterization of discrete convexity by linear inequalities}
\label{sec:Characterization}
\begin{figure}
\centering
\includegraphics[width=3.5cm]{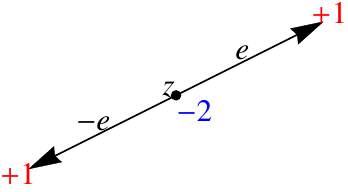}
\includegraphics[width=4cm]{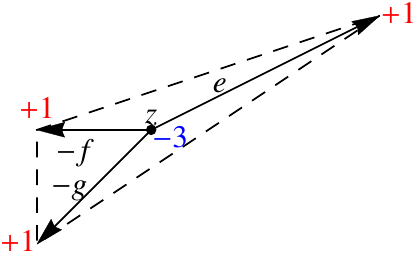}
\includegraphics[width=4cm]{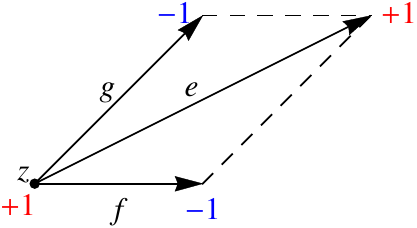}
\caption{The supports, and weights, of the different linear forms $S_x^e$, $T_x^e$, $P_x^e$.}
\label{fig:ConstraintTypes}
\end{figure}

We introduce some linear forms on the vector space $\cF(X) := \{u : X \to \R\}$, which non-negativity characterizes restrictions of convex maps. The convex hulls of their respective supports have respectively the shape of a segment, a triangle, and a parallelogram, see Figure \ref{fig:ConstraintTypes}. 
\begin{Definition}
\label{def:Constraints}
For each $x \in \FixedGrid$, consider the following linear forms of $u \in \cF(\Z^2)$. 
\begin{enumerate}
\item (Segments)
For any irreducible $e \in \FixedGrid$: 
\begin{equation*}
S_x^e (u) := u(x+e)-2 u(x)+u(x-e).
\end{equation*}

\item (Triangles)
For any irreducible $e \in \FixedGrid$, with $\|e\|>1$, of parents $f,g$: 
\begin{equation*}
T_x^e (u) := u(x+e)+u(x-f)+u(x-g) - 3 u(x).
\end{equation*}
\item (Parallelograms)
For any irreducible $e \in \FixedGrid$, with $\|e\|>1$, of parents $f,g$: 
\begin{equation*}
P_x^e(u) := u(x+e) - u(x+f) - u(x+g) + u(x).
\end{equation*}
\end{enumerate}
A linear form $L$ among the above can be regarded as a finite weighted sum of Dirac masses. In this sense we define the support $\supp(L) \subset \FixedGrid$, $\#\supp(L) \in \{3,4\}$. The linear form $L$ is also defined on $\cF(X)$ whenever $\supp(L) \subset X$. 
\end{Definition}

If $u \in \Conv(\OmegaD)$, then by an immediate convexity argument one obtains $S_x^e(u) \geq 0$ and $T_x^e(u) \geq 0$, whenever these linear forms are supported on $\OmegaD$. 
As shown in the next result, this provides a minimal characterization of $\Conv(\OmegaD)$ by means of linear inequalities. The linear forms $P_x^e$ will on the other hand be used to define strict sub-cones of $\Conv(\OmegaD)$.
The following result corrects%
\footnote{%
Precisely, the constraints $T_x^e$ were omitted in \cite{Carlier:2001tq} for $\|e\|>\sqrt 2$.
} 
Corollary 4 in \cite{Carlier:2001tq}.
\begin{Theorem} 
\label{th:AllConstraints}
\begin{itemize}
\item
The cone $\Conv(\OmegaD)$ is characterized by the non-negativity of the linear forms $S_x^e$ and $T_x^e$, introduced in Definition \ref{def:Constraints}, which are supported in $\OmegaD$. 
\item
If one keeps only one representative among the identical linear forms $S_x^e$ and $S_x^{-e}$, then the above characterization of $\Conv(\OmegaD)$ by linear inequalities is minimal.
\end{itemize}
\end{Theorem}

For any given $x \in \OmegaD$, the number of linear inequalities $S_x^e$ (resp.\ $T_x^e$) appearing in the characterization of $\Conv(\OmegaD)$ is bounded by the number of irreducible elements $e \in \Z^2$ such that $x+e \in X$. Hence the $N$-dimensional cone $\Conv(\OmegaD)$ is characterized by at most $2 N^2$ linear inequalities, where $N := \#(\OmegaD)$.

If all the elements of $\OmegaD$ are \emph{aligned}, this turns out to be an over estimate: one easily checks that exactly $N-1$ inequalities of type $S_x^e$ remain, and no inequalities of type $T_x^e$. 
This favorable situation does not extend to the two dimensional case however, because irreducible elements arise frequently in $\FixedGrid$, with positive density \cite{Hardy:1979vq}:
\begin{equation*}
\frac 6 {\pi^2} = \lim_{n \to \infty} n^{-2} \# \left\{(i,j) \in \{1,\cdots,n\}^2; \, \gcd(i,j)=1\right\}.
\end{equation*}
If the domain $\OmegaC$ has a non-empty interior, then one easily checks from this point that the minimal description of $\Conv(\OmegaD)$ given in Theorem \ref{th:AllConstraints} involves no less than $c N^2$ linear constraints%
\footnote{%
This number of constraints is empirically (and slightly erroneously) estimated to $\cO(N^{1.8})$ in \cite{Carlier:2001tq}.
}%
, where the constant $c > 0$ depends on the domain shape but not on its scale (or equivalently, not on the grid scale $h$ in \eqref{def:ParametrizedDiscreteDomain}). This quadratic number of constraints, announced in the description of global constraint methods in the introduction, is a strong drawback for practical applications, which motivates the construction of adaptive sub-cones of $\Conv(X)$ in the next subsection.

\begin{Remark}[Directional convexity]
Several works addressing optimization problems posed on the cone of convex functions \cite{Carlier:2001tq,Oberman:2011wi}, have in the past omitted all or part of the linear constraints $T_x^e$, $x\in X$, $e\in \Z^2$ irreducible with $\|e\|>1$. We consider in Appendix \ref{sec:Directional} this weaker notion of discrete convexity, introducing the cone $\DConv(X)$ of directionally convex functions, defined by the non-negativity of only $S_x^e$, $x \in X$, $e \in \Z^2$ irreducible. 

We show that elements of $\DConv(X)$ cannot in general be extended into globally convex functions, but that one can extend their restriction to a grid coarsened by a factor $2$. We also introduce a hierarchy of sub-cones of $\DConv(X)$, similar to the one presented in the next subsection. 
\end{Remark}

\subsection{Hierarchical cones of discrete convex functions}

We introduce in this section the notion of stencils $\cV = (\cV(x))_{x \in X}$ on $X$, and discuss the properties (hierarchy, complexity) of cones $\Conv(\cV)$ attached to them.
The following family $\cV_{\max}$ of sets is referred to as the ``maximal stencils'': for all $x \in X$
\begin{equation}
\label{def:MaximalStencils}
\cV_{\max}(x) := \{ e \in \Z^2 \text{ irreducible}; \, x+e \in X\}.
\end{equation}
The convex cone generated by a subset $A$ of a vector space is denoted by $\Cone(A)$, with the convention $\Cone(\emptyset) = \{0\}$. 

\begin{Definition}
\label{def:Stencil}
A family $\cV$ of stencils on $\OmegaD$ (or just: ``Stencils on $X$'') is the data, for each $x \in \OmegaD$ of a collection  $\cV(x) \subset \cV_{\max}(x)$ (the stencil at $x$) of irreducible elements of $\Z^2$, satisfying the following properties:
\begin{itemize}
\item (Stability) Any parent $f \in \cV_{\max}(x)$, of any $e \in \cV(x)$, satisfies $f \in \cV(x)$. 
\item (Visibility) One has $\Cone(\cV(x)) = \Cone(\cV_{\max}(x))$. 
\end{itemize}
The set of candidates for refinement $\Candidates(x)$ consists of all elements $e\in \cV_{\max}(x) \sm \cV(x)$ which two parents $f,g$ belong to $\cV(x)$.
\end{Definition}

\begin{figure}
\centering
\includegraphics[width=4cm]{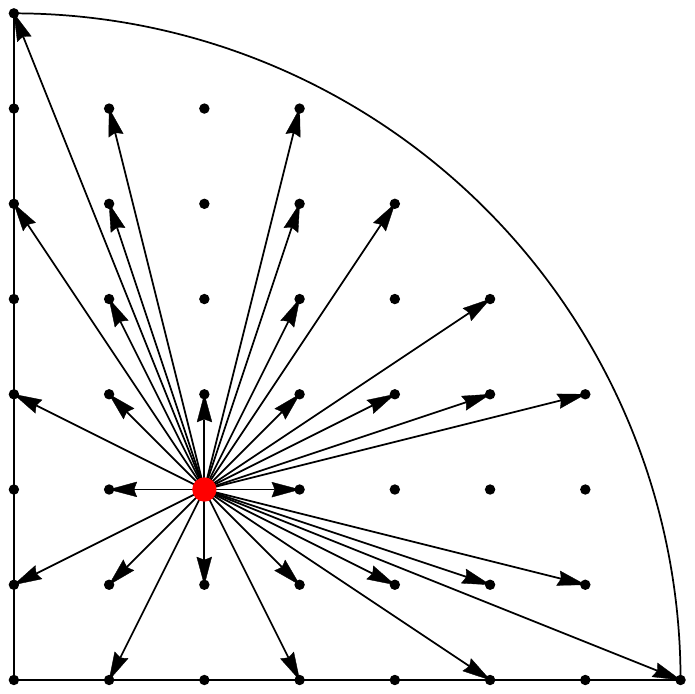}
\includegraphics[width=4cm]{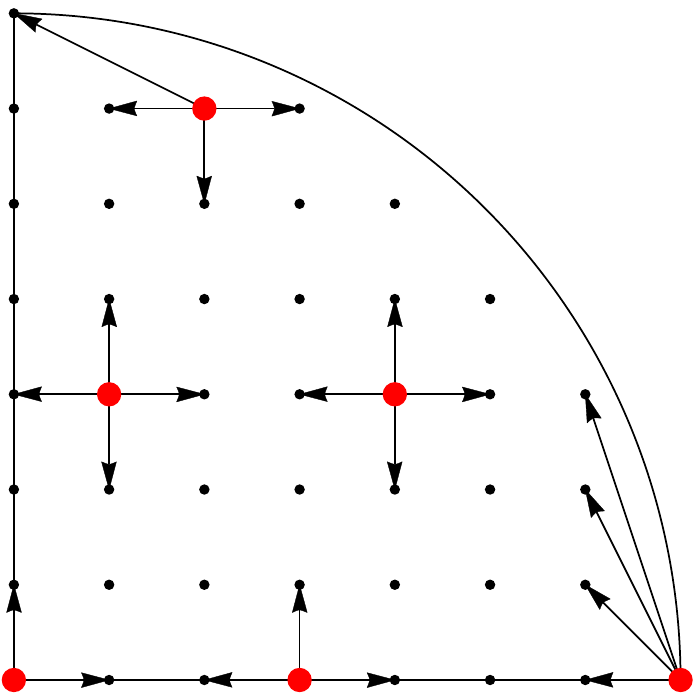}
\includegraphics[width=4cm]{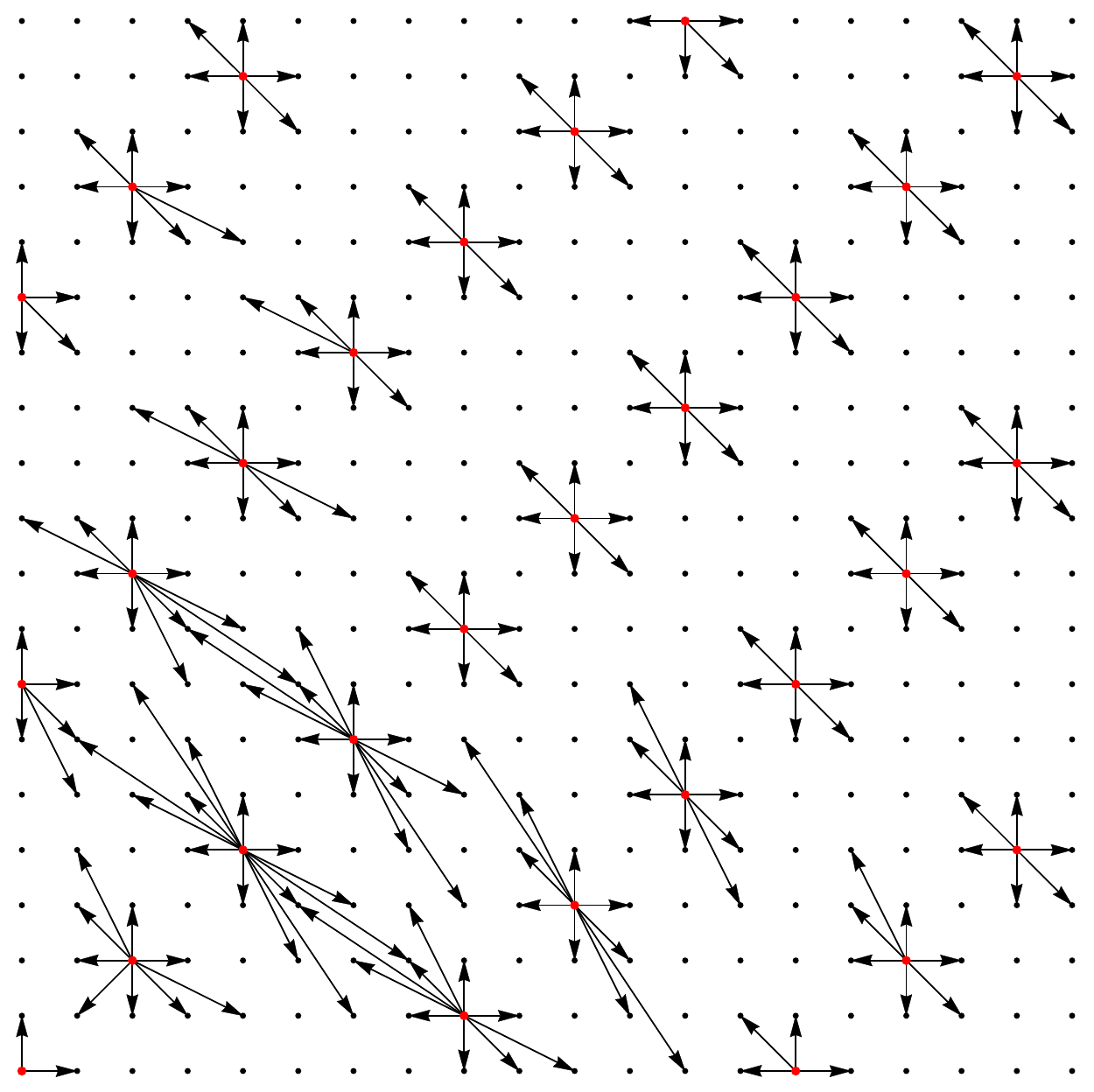}
\caption{Left: a maximal stencil at a point of a domain. Center: some minimal stencils.  Right: some adaptively generated stencils used in the numerical resolution of \eqref{eq:PrincipalAgent}.}
\label{fig:Stencils}
\end{figure}

In other words, a stencil $\cV(x)$ at a point $x\in \Omega$ contains the parents of its members whenever possible (Stability), and covers all possible directions (Visibility). By construction, these properties are still satisfied by the \emph{refined stencil} $\cV(x) \cup \{e\}$, for any candidate for refinement $e \in \Candidates(x)$. The collection $\Candidates(x)$ is easily recovered from $\cV(x)$, see Proposition \ref{prop:Trigo}. 

\begin{Definition}
\label{def:Cones}
We attach to a family $\cV$ of stencils on $\OmegaD$ the cone $\Conv (\cV) \subset \cF(\OmegaD)$, characterized by the non-negativity of the following linear forms: for all $x\in X$
\begin{enumerate}
\item 
$S_x^e$, for all $e \in \cV(x)$ such that $\supp(S_x^e) \subset \OmegaD$.
\item 
$T_x^e$ for all $e \in \cV(x)$, with $\|e\|>1$, such that $\supp(T_x^e) \subset \OmegaD$.
\item 
$P_x^e$ for all $e \in \Candidates(x)$ (by construction $\supp(P_x^e) \subset X$). 
\end{enumerate}
\end{Definition}
When discussing unions, intersections, and cardinalities, we (abusively) identify a family $\cV$ of stencils on $\OmegaD$ with a subset of $X \times \Z^2$:
\begin{equation}
\label{StencilIdentification}
\cV \approx \{(x,e); \, x\in \OmegaD, \, e \in \cV(x)\}.
\end{equation}
Note that the cone $\Conv(\cV)$ is defined by at most $3\#(\cV)$ linear inequalities. 
The sets $\cV_{\max}$ are clearly stencils on $X$, which are maximal for inclusion, and by Theorem \ref{th:AllConstraints} we have $\Conv(\cV_{\max}) = \Conv(X)$.
The cone $\Conv(\cV)$ always contains the quadratic function $q(x) := \frac 1 2 \|x\|^2$, for any family $\cV$ of stencils.  
Indeed, the inequalities $S_x^e(q) \geq 0$, $x \in X$, $e \in \cV(x)$, and $T_x^e(q) \geq 0$, $\|e\|>1$, hold by convexity of $q$. In addition for all $e \in \Candidates(x)$, of parents $f,g$, one has
\begin{equation*}
P_x^e(q)  =  \frac 1 2 \left( \|x+f+g\|^2 - \|x+f\|^2 -\|x+g\|^2+\|x\|^2 \right) =  \<f,g\> \geq 0,
\end{equation*}
since the basis $(f,g)$ of $\FixedGrid$ is acute by definition, see Proposition \ref{prop:Parents}.

\begin{Theorem}[Hierarchy]
\label{th:Hierarchy}
The union $\cV \cup \cV'$, and the intersection $\cV \cap \cV'$ of two families $\cV, \cV'$ of stencils  are also families of stencils on $X$. In addition
\begin{align}
\label{StencilIntersection}
\Conv(\cV) \cap \Conv( \cV') &= \Conv(\cV \cap \cV'), \\
\label{StencilUnion}
\Conv(\cV) \cup \Conv( \cV')&\subset \Conv(\cV \cup \cV') .
\end{align}
\end{Theorem}
As a result, if two families of stencils $\cV, \cV'$ satisfy $\cV \subset \cV'$, then
\begin{equation*}
\Conv(\cV) \subset \Conv(\cV') \subset \Conv(X).
\end{equation*}
The left inclusion follows from \eqref{StencilIntersection}, and the right inclusion from \eqref{StencilUnion} applied to $\cV'$ and $\cV_{\max}$. The intersection rule \eqref{StencilIntersection} also implies the existence of stencils $\cV_{\min}$ minimal for inclusion, which are illustrated on Figure \ref{fig:Stencils} and characterized in Proposition \ref{prop:MinimalStencilsCharacterization}.
 

\begin{Remark}[Optimization strategy]
\label{rem:OptStrategy}
For any $u \in \Conv(x)$, there exists by \eqref{StencilIntersection} a unique smallest (for inclusion) family of stencils $\cV$ such that $u \in \Conv(\cV)$.
If $u$ is the minimizer of an energy $\cE$ on $\Conv(X)$, then it can be recovered by minimizing $\cE$ on the smaller cone $\Conv(\cV)$, defined by $\cO(\#(\cV))$ linear constraints. Algorithm \ref{algo:SubCones} in \S \ref{sec:Numerics}, attempts to find these smallest stencils $\cV$ (or slightly larger ones), starting from $\cV_{\min}$ and performing successive adaptive refinements.
\end{Remark}

In the rest of this subsection, we fix a grid scale $h>0$ and consider for all
$\theta \in \R$, and all $\xi \in \R^2$, the grid 
\begin{equation}
\label{def:XThetaXi}
X_\theta^\xi := \Omega \cap h R_\theta (\xi+\Z^2).
\end{equation}
The notions of stencils and related cones trivially extend to this setting, see \S \ref{sec:Avg} for details.
We denote by $|\Omega|$ the domain area, and by $\diam(\Omega) := \max \{\|y-x\|; \, x,y \in \Omega\}$ its diameter. We also introduce rescaled variants, defined for $h>0$ by 
\begin{equation*}
|\Omega|_h := h^{-2} |\Omega|, \qquad \diam_h(\Omega) := h^{-1} \diam(\Omega).
\end{equation*}
For any parameters $\theta, \xi$, one has denoting $N := \#(X_\theta^\xi)$ (with underlying constants depending only on the shape of $\Omega$)
\begin{equation}
\label{eq:ApproxN}
|\Omega|_h \approx N, \qquad \diam_h(\Omega) \approx \sqrt{N}.
\end{equation}


\begin{Proposition}
\label{prop:WorstCase}
Let $X := X_\theta^\xi$, for some grid position parameters $\theta\in \R$, $\xi\in \R^2$, and let $N := \#(X)$.
Let $u \in \Conv(X)$, and let $\cV$ be the minimal stencils on $X$ such that $u \in \Conv(\cV)$. Then $\#(\cV) \leq C N \diam_h(\Omega)$, for some universal constant $C$ (i.e.\ independent of $\Omega,h,\theta, \xi,u$).
\end{Proposition}

Combining this result with \eqref{eq:ApproxN} we see that an optimization strategy as described in Remark \ref{rem:OptStrategy} should heuristically not require solving optimization problems subject to more than $N \diam_h(\Omega) \approx N^{\frac 3 2}$ linear constraints. This is already a significant improvement over the $\approx N^2$ linear constraints defining $\Conv(X)$. The typical situation is however even more favorable: in average over randomized grid orientations $\theta$ and offsets $\xi$, 
the restriction to $X_\theta^\xi$ of a convex map $U : \Omega \to \R$ (e.g.\ the global continuous solution of the problem \eqref{eq:PrincipalAgent} of interest) belongs to a cone $\Conv(\cV_\theta^\xi)$ defined by a quasi-linear number $\cO(N \ln^2 N)$ of linear inequalities.

\begin{Theorem}
\label{th:Avg}
Let $U \in \Conv(\Omega)$, and let $\cV_\theta^\xi$ be the minimal stencils on $X_\theta^\xi$ such that $U_{|X_\theta^\xi} \in \Conv(\cV_\theta^\xi)$, for all $\theta \in \R$, $\xi \in \R^2$. Assuming $\diam_h(\Omega) \geq 2$, one has for some universal constant $C$ (i.e.\ independent of $h,\Omega,U$): 
\begin{equation}
\label{eq:Avg}
\int_{[0,1]^2} \int_0^{\pi/2} \#(\cV_\theta^\xi) \, d \theta \, d \xi \leq C \, |\Omega|_h  \, (\ln\diam_h(\Omega))^2.
\end{equation}
\end{Theorem}


\subsection{Stencils and triangulations}
\label{sec:DelaunayIntro}

We discuss the connections between stencils and triangulations, which provides in Theorem \ref{th:Decomp} a new insight on the hierarchy of cones $\Conv(\cV)$, and yields in Theorem \ref{th:Delaunay} a new result of algorithmic geometry as a side product of our theory.
We assume in this subsection and \S \ref{sec:Delaunay} that the discrete domain convex hull, denoted by $\Hull(X)$, has a non-empty interior. All triangulations considered in this paper are implicitly assumed to cover $\Hull(X)$ and to have $X$ as collection of vertices.
\begin{Definition}
\label{def:TinV}
Let $\cT$ be a triangulation, and let $\cV$ be a family of stencils on $X$. We write $\cT \prec \cV$ iff the directed graph associated to $\cT$ is included in the one associated to $\cV$. In other words iff for any edge $[x,x+e]$ of $\cT$, one has $e \in \cV(x)$.
\end{Definition}

The next result provides a new interpretation to our approach to optimization problems posed on the cone of convex functions, as a relaxation of the naïve (and flawed without this modification) method via interior finite elements.
\begin{Theorem}
\label{th:Decomp}
Let $\cV$ be a family of stencils on $X$. If $\Conv(\cV)$ has a non-empty interior, then  
\begin{equation}
\label{eq:partitionV}
\Conv(\cV) = \bigcup_{T \prec \cV} \Conv(\cT).
\end{equation}
\end{Theorem}

\begin{figure}
\includegraphics[width=3.9cm]{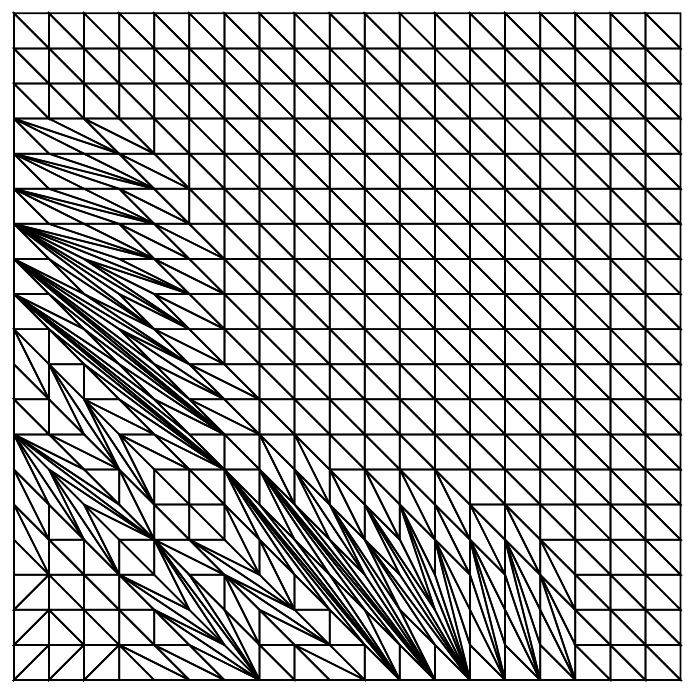}
\includegraphics[width=3.9cm]{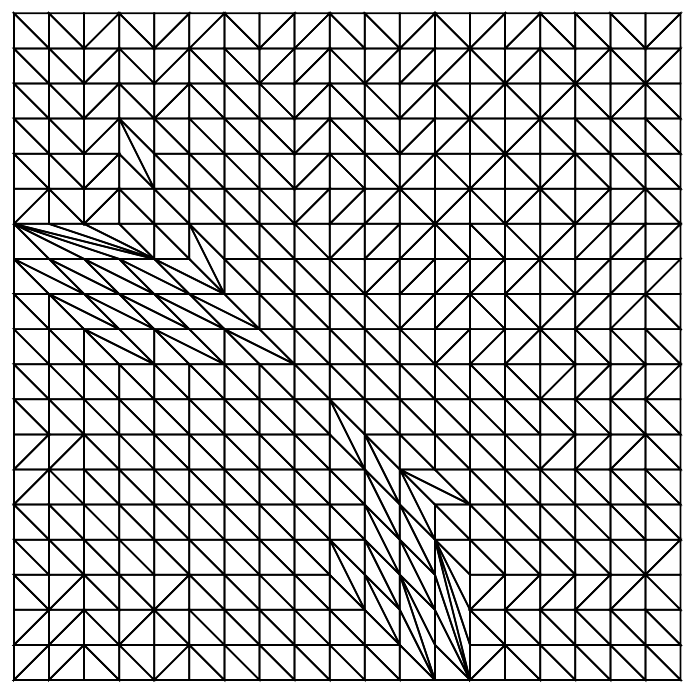}
\includegraphics[width=3.9cm]{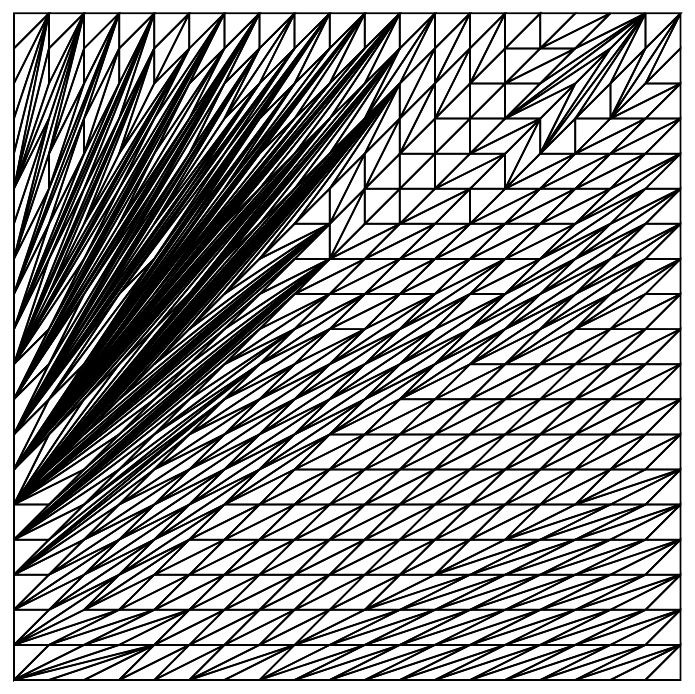}
{\raise 1.5cm \hbox{
\begin{tabular}{cc}
\multicolumn{2}{c}{
\includegraphics[width=3cm]{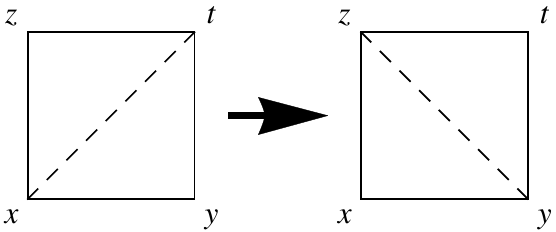}
}
\\
\includegraphics[width=1.2cm]{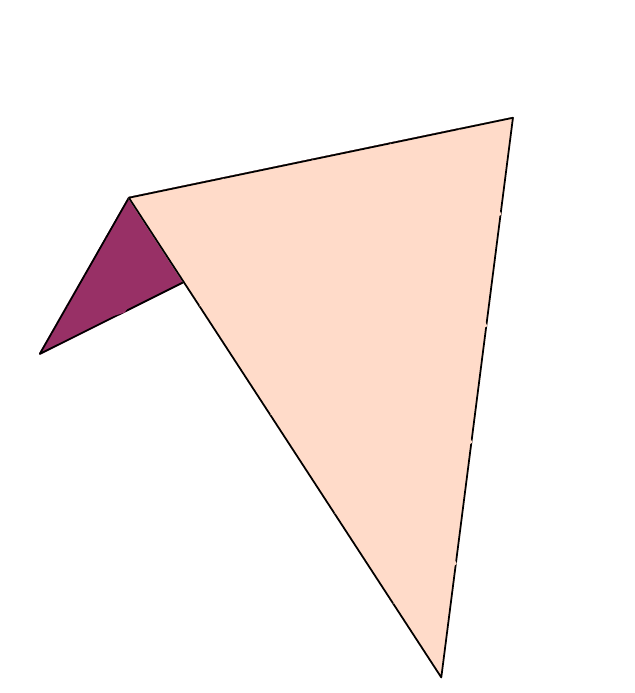} &
\includegraphics[width=1.2cm]{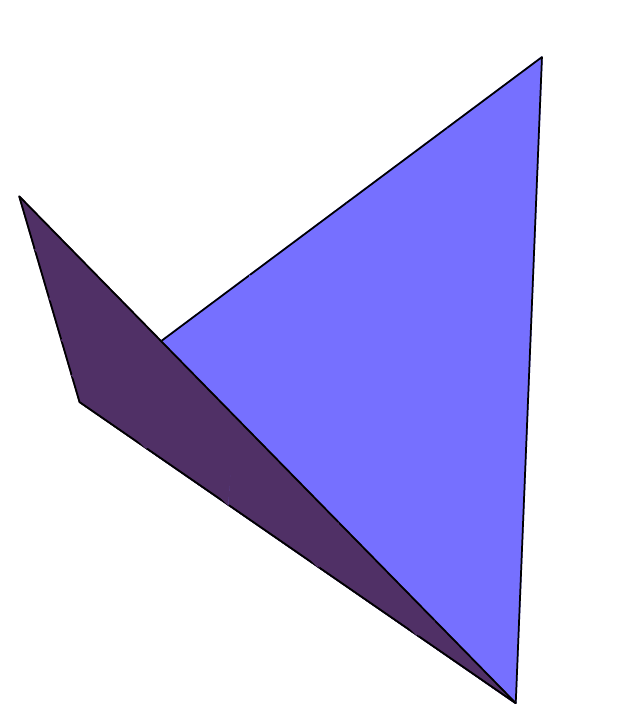}
\end{tabular}
}}
\caption{Delaunay triangulations associated to numerical solutions of (variants of) the monopolist problem, see \ref{sec:Monopolist}. Corresponding convex functions shown on Figure \ref{fig:Monopolist} (left), Figure \ref{fig:Bundles} (left), and Figure \ref{fig:Assets} (left) respectively.
Right: illustration of edge flipping. The diagonal $[x,t]$ shared by two triangles $T:=[x,y,t]$ and $T':=[x,z,t]$ can be flipped into $[y,z]$ if $T \cup T'$ is convex. Right, below: the piecewise linear interpolation of a discrete map $u$ is made convex by this flip.} 
\label{fig:Flipping}
\end{figure}

Delaunay triangulations are a fundamental concept in discrete geometry \cite{Edelsbrunner:1986wt}. We consider in this paper a slight generalization in which the lifting map needs not be the usual paraboloid, but can be an arbitrary convex function, see Definition \ref{def:Delaunay}. Within this paper Delaunay triangulations are simultaneously (i) a theoretical tool for proving results, notably Proposition \ref{prop:WorstCase} and Theorem \ref{th:Decomp}, (ii) an object of study, since in Theorem \ref{th:Delaunay} we derive new results on the cost of their construction, and (iii) a numerical post-processing tool providing global convex extensions of elements of $\Conv(X)$, see Figure \ref{fig:Flipping} and Remark \ref{rem:Subgradients}. 

%
%

\begin{Definition}
\label{def:Delaunay}
We say that $\cT$ is an $u$-Delaunay triangulation iff $u \in \Conv(\cT)$; equivalently the piecewise linear interpolation $\interp_\cT u : \Hull(X) \to \R$ is convex. We refer to $u$ as the \emph{lifting map}. 

A $q$-Delaunay triangulation, with $q(x) := \frac 1 2 \|x\|^2$, is simply called a Delaunay triangulation.
\end{Definition}

Two dimensional Delaunay triangulations, and three dimensional convex hulls, have well known links \cite{Edelsbrunner:1986wt}.
Indeed $\cT$ is an $u$-Delaunay triangulation iff 
the map $x \in \Hull(X) \mapsto (x, \interp_\cT u(x)) \in \R^3$ spans the bottom part of the convex envelope $K$ of lifted set $\{(x,u(x)); \, x \in X\}$. 
As a result of this interpretation, we find that (i) any element $u \in \Conv(X)$ admits an $u$-Delaunay triangulation, and (ii) generic elements of $\Conv(X)$ admit exactly one $u$-Delaunay triangulation (whenever all the faces of $K$ are triangular). In particular, the union \eqref{eq:partitionV} is disjoint up to a set of Hausdorff dimension $N-1$.

Any two triangulations $\cT, \cT'$ of $X$ can be transformed in one another through a sequence of elementary modifications called \emph{edge flips}, see Figure \ref{fig:Flipping}. The minimal number of such operations is called the edge flipping distance between $\cT$ and $\cT'$.
Edge flipping is a simple, robust and flexible procedure, which is used in numerous applications ranging from fluid dynamics simulation \cite{Dobrzynski:2008ix} to GPU accelerated image vectorization \cite{Qi:2012fd}. 
Sustained research has been devoted to estimating edge flipping distances within families of triangulations of interest \cite{Hurtado:1999jj}, although flipping distance bounds are usually quadratic in the number of vertices. 
\begin{Theorem}
\label{th:Delaunay}
Let $\cV$ be a family of stencils on $X$, and let $u \in \Conv(\cV)$. Then any standard Delaunay triangulation of $X$ can be transformed into an $u$-Delaunay triangulation via a sequence of $C\#(\cV)$ edge flips. The constant $C$ is universal (in particular it is independent of $u, \cV, \Omega$).
\end{Theorem}
Combining this result with Theorem \ref{th:Avg}, we obtain that only $\cO(N \ln^2 N)$ edge flips are required to construct an $U$-Delaunay triangulation of $X$, with $N := \#(X)$, for any convex function $U \in \Conv(\Omega)$, in an average sense over grid orientations. 
Note that (more complex and specialized) convex hull algorithms \cite{PrincetonUniversityDeptofComputerScience:1991uv}  could also be used to produce an $U$-Delaunay triangulation, at the slightly lower cost $\cO(N \ln N)$. Theorem \ref{th:Delaunay}  should be understood as a first step in understanding the typical behavior of edge flipping.

\subsection{Outline}

We prove in \S \ref{sec:ConstraintConvexity} the characterization of discrete convexity by linear constraints of Theorem \ref{th:AllConstraints}. The hierarchy properties of Theorem \ref{th:Hierarchy} are established in \S \ref{sec:Hierarchy}. Triangulation related arguments are used in \S \ref{sec:Delaunay} to show Proposition \ref{prop:WorstCase} and Theorems \ref{th:Decomp} and \ref{th:Delaunay}. The average cardinality estimate of Theorem \ref{th:Avg} is proved in \S \ref{sec:Avg}. Numerical experiments, and algorithmic details, are presented in \S \ref{sec:Numerics}. Finally, the weaker notion of directional convexity is discussed in Appendix \ref{sec:Directional}.



\section{Characterization of convexity via linear constraints}
\label{sec:ConstraintConvexity}

This section is devoted to the proof of Theorem \ref{th:AllConstraints}, which characterizes discrete convex functions $u \in \Conv(X)$ via linear inequalities. The key ingredient is its generalization, in \cite{Carlier:2001tq}, to arbitrary unstructured finite point sets $X'$ (in contrast with the grid structure of $X$).

\begin{Theorem}[Carlier, Lachand-Robert, Maury]
\label{th:ConvexityUnstructured}
Let $\OmegaC$ be a convex domain, and let $\OmegaD' \subset \OmegaC$ be an arbitrary finite set. The cone $\Conv(\OmegaD')$, of all restrictions to $\OmegaD'$ of convex functions on $\OmegaC$, is characterized by the following inequalities, none of which can be removed:
\begin{itemize}
\item
For all $x,y,z \in \OmegaD'$, all $\alpha \in ]0,1[$, $\beta := 1-\alpha$,
such that $z = \alpha x+ \beta y$ and $[x,y] \cap \OmegaD' = \{x,y,z\}$:
\begin{equation}
\label{SegmentUnstructured}
\alpha u(x) + \beta u(y) \geq u(z).
\end{equation}
\item  
For all $p,q,r,z \in \OmegaD'$, all
$\alpha, \beta, \gamma\in \R_+^*$ with $\alpha+\beta+\gamma=1$,
such that $z = \alpha p+ \beta q+\gamma r$, the points $p,q,r$ are not aligned and $[p,q,r]\cap \OmegaD' = \{p,q,r,z\}$:
\begin{equation}
\label{TriangleUnstructured}
\alpha u(p) + \beta u(q) + \gamma u(r) \geq u(z).
\end{equation}
\end{itemize}
\end{Theorem}

In the following, we establish Theorem \ref{th:AllConstraints} by specializing Theorem \ref{th:ConvexityUnstructured} to the grid $\OmegaD := \OmegaC \cap \Z^2$, and identifying  \eqref{SegmentUnstructured} and \eqref{TriangleUnstructured} with the constraints $S_z^e(u) \geq 0$ and $T_z^e(u) \geq 0$ respectively, in Propositions \ref{prop:Seg} and \ref{prop:Tri} respectively. Note that the cases \eqref{SegmentUnstructured} and \eqref{TriangleUnstructured} are unified in \cite{Carlier:2001tq}, although we separated them in the above formulation for clarity. For $x_1, \cdots, x_n \in \R^2$, we denote 
\begin{equation*}
[x_1, \cdots, x_n] := \Hull(\{x_1, \cdots, x_n\}).
\end{equation*}

\begin{Proposition}
\label{prop:Seg}
Let $x,y,z \in \OmegaD$, 
$\alpha \in ]0,1[$, $\beta:=1-\alpha$, be
such that $z = \alpha x+ \beta y$ and $[x,y] \cap \OmegaD = \{x,y,z\}$. Then $\alpha = \beta = 1/2$ and there exists an irreducible $e\in \FixedGrid$ such that $x=z+e$, $y=z-e$. (Thus $\alpha u(x)+ \beta u(y)-u(z) =  S_z^e(u)/2$.)
\end{Proposition}

\begin{proof}
We define $e := x-z \in \FixedGrid$, and assume for contradiction that $e$ is not irreducible: $e=k e'$, for some integer $k \geq 2$ and some $e' \in \FixedGrid$. 
Then $z+e' \in [z,x] \cap \Z^2 \subset [x,y] \cap X$, which is a contradiction.
Thus $e$ is irreducible, and likewise $f := y-z \in \FixedGrid$ is irreducible. 
Observing that $f$ is negatively proportional to $e$, namely $f = -(\alpha/ \beta) e$,  we obtain that $f=-e$, which concludes the proof.
\end{proof}

The next lemma, used in Proposition \ref{prop:Tri} to identify the constraints \eqref{TriangleUnstructured}, provides an alternative characterization of the parents of an irreducible vector, see Proposition \ref{prop:Parents}. 

\begin{Lemma}
\label{lem:SumParents}
Let $e,f,g\in \Z^2$ be such that $e=f+g$, $\|e\| \geq \max \{\|f\|,\|g\|\}$, and $\det(f,g)=1$. Then $f,g$ are the parents of $e$.
\end{Lemma}

\begin{proof}
Since $\det(f,e) = \det(f,f+g)=1$, the vector $e$ is irreducible. In addition, $\|e\|>1$, since otherwise $e,f,g$ would be three pairwise linearly independent unit vectors in $\Z^2$.
Let $f_0,g_0$ be the parents of $e$.  Observing that $\det(f,e) = \det(f_0,e)=1$, we obtain that $f = f_0+k e$ for some scalar $k$, which must be an integer since $e$ is irreducible. If $k > 0$ then $\|f\|^2 = k^2 \|e\|^2 + 2 k \<f_0,e\> +\|f_0\|^2 > \|e\|^2$ which is a contradiction (recall that $\<f_0,e\> = \<f_0,f_0+g_0\> \geq 0$). If $k<0$, then observing that $g = e-f = g_0-k e$ we obtain likewise a contradiction. Hence $k=0$ and $f,g$ are the parents of $e$, which concludes the proof.
\end{proof}

\begin{Corollary}
\label{corol:SumParent}
Let $(f,g)$ be a basis of $\Z^2$ which is \emph{not} acute. If $\|f\| \geq \|g\|$, then $f+g$ is a parent of $f$, and otherwise it is a parent ot $g$.
\end{Corollary}

\begin{proof}
Up to exchanging the roles of $f,g$, we may assume that $\det(f,g)=1$. Denoting $m := \max\{\|f\|,\|g\|,\|f+g\|\}$, we have by Lemma \ref{lem:SumParents} three possibilities: (i) $\|f+g\|=m$, and $f,g$ are the parents of $f+g$, (ii) $\|f\|=m$, and $-g,f+g$ are the parents of $f$, (iii) $\|g\|=m$, and $f+g,-f$ are the parents of $g$. Excluding (i), since $(f,g)$ is not an acute basis, we conclude the proof.
\end{proof}

\begin{Proposition}
\label{prop:Tri}
Let $p,q,r,z \in \OmegaD$, $\alpha, \beta, \gamma \in \R_+^*$, with $\alpha+\beta+\gamma=1$, be such that $z = \alpha p+ \beta q + \gamma r$, the points $p,q,r$ are not aligned, and $[p,q,r] \cap \OmegaD = \{p,q,r,z\}$.
Then $\alpha = \beta = \gamma = 1/3$ and, up to permuting $p,q,r$, there exists an irreducible $e\in \FixedGrid$ with $\|e\| > 1$, of parents $f,g$, such that  $p =z+e$, $q=z-f$, and $r=z-g$. (Thus $\alpha u(p) + \beta u(q) + \gamma u(r) - u(z) =  T_z^e(u)/3$.)

\end{Proposition}

\begin{proof}
Let $e:=p-z$, $f:=z-q$, $g:=z-r$. Up to permuting $p,q,r$, we may assume that $\|e\| \geq \max \{\|f\|, \|g\|\}$ and $\det(f,g) \geq 0$.
Note that $f$ and $g$ are not collinear since $z$ lies in the interior of $[p,q,r]$. We claim that $(f,g)$ is a basis of $\FixedGrid$. Indeed, otherwise, the triangle $[0,f,g]$ would contain an element of $e' \in \FixedGrid$ distinct from its vertices. Since $\OmegaC$ is convex, this implies $[p,q,r] \cap \OmegaD \supset \{p,q,r,z,z+e'\}$, which is a contradiction. 

Thus $(f,g)$, and likewise $(e,f)$ and $(e,g)$, are bases of $\FixedGrid$, and therefore
\begin{equation}
\label{eq:Det}
|\det(e,f)| = |\det(f,g)| = |\det(g,e)|=1.
\end{equation}
Injecting in the above equation the identity $e = (\beta/\alpha) f + (\gamma/ \alpha) g$, we obtain that $|\beta/\alpha |=1$ and $|\gamma/\alpha| = 1$. Thus $\alpha = \beta = \gamma = 1/3$ since these coefficients are positive and sum to one.
Finally, we have $e = f+g$, $\|e\| \geq \max \{\|f\|, \|g\|\}$, and $f,g$ is a direct basis. Using Lemma \ref{lem:SumParents} we conclude as announced that $f,g$ are the parents of $e$. \end{proof}


\section{Hierarchy of the cones $\Conv(\cV)$} 
\label{sec:Hierarchy}

This section is devoted to the proof of Theorem \ref{th:Hierarchy}, which is split into two parts: the proof that an intersection (or union) of stencils is still a stencil, and the hierarchy properties \eqref{StencilIntersection} and \eqref{StencilUnion}.

\subsection{An intersection of stencils is still a stencil}

Let $\cV, \cV'$ be families of stencils on $X$. Property (Stability) of stencils immediately holds for the intersection $\cV \cap \cV'$ and union, $\cV \cup \cV'$, while property (Visibility) is also clear for the union $\cV \cup \cV'$. In order to establish property (Visibility) for the intersection $\cV \cap \cV'$, we identify in Proposition \ref{prop:MinimalStencils} a family $\cV_{\min}$ of stencils which included in any other. From $\cV_{\min} \subset \cV$ and $\cV_{\min} \subset \cV'$ we obtain $\cV_{\min} \subset \cV \cap \cV'$, so that (Visibility) for $\cV_{\min}$ implies the same property for $\cV \cap \cV'$.



\begin{Definition}
\label{def:CyclicOrder}
The cyclic strict trigonometric order on $\R^2 \sm \{0\}$ is denoted by $\prec$. 
\end{Definition}
In other words $e_1 \prec e_2 \prec e_3$ iff there exists $\theta_1, \theta_2, \theta_3 >0$, such that $\theta_1+\theta_2+\theta_3 = 2 \pi$ and $e_{i+1}/\|e_{i+1}\| = R_{\theta_i} e_i / \|e_i\|$ for all $1 \leq i \leq 3$, with $e_4 : =e_1$. 
The following lemma, and Corollary \ref{corol:ChildrenInCone}, discuss the combination of the cyclic ordering with the notions of parents (and children) of an irreducible vector.

\begin{Definition}[Collection of ancestors of a vector]
\label{def:Ancestors}
For any irreducible $e \in \Z^2$, let $\Anc(e)$ be the smallest set containing $e$ and the parents of any element $e' \in \Anc(e)$ such that $\|e'\|>1$.
\end{Definition}

\begin{Lemma}
\label{lem:ConsecutiveParents}
Let $e \in \Z^2\sm \{0\}$, let $(f,g)$ be a direct basis of $\Z^2$ such that $f \prec e \prec g$, and let us consider the triangle $T := [e,f,g]$. Then (i) $f+g \in T$. If in addition $e$ is irreducible, $\|e\|>1$ and (ii.a) $\<f,g\> \geq 0$ or (ii.b) $e \notin \Anc(f) \cup \Anc(g)$, then the parents of $e$ also belong to $T$.
\end{Lemma}

\begin{proof}
Point (i). By construction, we have $e=\alpha x+\beta y$ for some positive integers $\alpha, \beta$. One easily checks that $e + (\beta-1) f + (\alpha-1) g = (\alpha+\beta-1) (f+g)$. This expression of $f+g$ as a weighted barycenter of the points $e,f,g,$ establishes (i).

Points (ii.a) and (ii.b).
We fix $e$ and show these points by \emph{decreasing} induction on the integer $k = \<f,g\>$.
Initialization: Assume that $k = \<f,g\> \geq \frac 1 2 \|e\|^2$. Then $\|e\|^2 = \|\alpha f + \beta g\|^2 \geq 2 \alpha \beta \<f,g\> \geq 2 \<f,g\>$, which is a contradiction. No basis $(f,g)$ satisfies simultaneously $f\prec e \prec g$ and $\<f,g\> = k$. The statement is vacuous, hence true.

Case $k=\<f,g\> \geq 0$. If $e=f+g$, then $f,g$ are the parents of $e$, and the result follows. Otherwise, we have either $f \prec e \prec (f+g)$ or $(f+g)\prec e \prec g$. Since $\<f,f+g\> > \<f,g\>$ and $\<f+g,g\> > \<f,g\>$, we may apply our induction hypothesis to the bases $(f,f+g)$ and $(f+g,g)$ which satisfy (ii.a). Thus the parents of $e$ belong to $T_1 := [e,f,f+g]$ or $T_2 := [e,g,f+g]$. Finally, Point (i) implies that $f+g \in T$, thus $T_1 \cup T_2 \subset T$ which concludes the proof of this case.

Case $k = \<f,g\> < 0$. Assumption (ii.b) must hold, since (ii.a) contradicts this case. By corollary \ref{corol:SumParent}, $f+g$ is a parent of $f$ or of $g$. Hence $e \neq f+g$ and $\Anc(f+g) \subset \Anc(f) \cup \Anc(g)$. We apply our induction hypothesis to the bases $(f,f+g)$ and $(f+g,g)$ which satisfy (ii.b), and conclude the proof similarly to the case $k \geq 0$. 
\end{proof}

\begin{Lemma}
\label{lem:AllChildren}
Consider an irreducible $e \in \Z^2$, $\|e\|>1$, and let $f,g$ be its parents. The children of $e$ (i.e. the vectors $e' \in \Z^2$ of which $e$ is a parent) have the form $f+k e$ and $k e + g$, $k \geq 1$.
\end{Lemma}

\begin{proof}
Let $e'$ be a children of $e$, and let $f'$, $g'$ be its parents. Without loss of generality, we assume that $g'=e$. Then $\det(f',e) = 1 = \det(f,e)$, thus $f' = f+k e$ for some $k \in \R$. Since $e$ is irreducible, one has $k \in \mZ$. Since $0 \leq \<f',e\> = \<f,e\> + k \|e\|^2 \leq (k-1)\|e\|^2$, one has $k \geq 1$. The result follows.
\end{proof}

\begin{Corollary}
\label{corol:ChildrenInCone}
Let $e \in \Z^2\sm \{0\}$, let $(f,g)$ be a direct acute basis of $\Z^2$ such that $f \prec e \prec g$. Then any child $e'$ of $e$ (i.e.\ $e$ is a parent of $e'$) satisfies $f \prec e' \prec g$. 
\end{Corollary}

\begin{proof}
Let $K := \Cone(\{f,g\})$, and let $\interior K$ be its interior. Let also $f',g'$ be the parents of $e$. By assumption $e \in \interior K$, and by Lemma \ref{lem:ConsecutiveParents} (ii.a) one has $f',g' \in [e,f,g] \subset K$. 
Thus $f'+ k e, k e+g' \in \interior K$ for any integer $k \geq 1$, which by Lemma \ref{lem:AllChildren} concludes the proof.
\end{proof}


\begin{Lemma}[Consecutive elements of a stencil]
\label{lem:ConsecutiveBasis}
Let $\cV$ be a family of stencils on $\OmegaD$, let $x \in \OmegaD$, and let $f,g$ be two trigonometrically consecutive elements of $\cV(x)$. Then either (i) $(f,g)$ form a direct acute basis, or (ii) 
no element $e \in \cV_{\max}(x)$ satisfies  $f \prec e \prec g$.
\end{Lemma}

\begin{proof}
We distinguish three cases, depending on the value of $\det(f,g)$.
In case $\det(f,g) \leq 0$, property Visibility of stencils implies (ii). 

Case $\det(f,g)=1$. 
Assuming that (i) does not hold, Corollary \ref{corol:SumParent} implies that $f+g$ is a parent of $f$ or of $g$. Assuming that (ii) does not hold, we have $f+g \in [e,f,g]$ for some $e\in \cV_{\max}(x)$ by Lemma \ref{lem:ConsecutiveParents} (i), thus $f+g \in \cV_{\max}(x)$ by convexity of $\Omega$, thus $f+g \in \cV(x)$ by (Stability), which contradicts our assumption $f,g$ are trigonometrically consecutive in $\cV(x)$.

Case $k:=\det(f,g) > 1$.
We assume without loss of generality that $\|f\| \geq \|g\|$, hence $\|f\|^2 \geq \det(f,g) > 1$. Let $f',g'$ be the parents of $f$, so that $g = \alpha f' + \beta g'$ for some $\alpha, \beta \in \Z$.
We obtain $k = \det(f,\alpha f' + \beta g') = \beta-\alpha$. If $\alpha=0$ or $\beta=0$, then $g$ is not irreducible, which is a contradiction. If $\alpha$ and $\beta$ have the same sign, then $\|g\|>\|f\|$, which again is a contradiction. Hence $\alpha, \beta$ have opposite signs, and since $\beta-\alpha=k>1$ we obtain $\beta>0>\alpha$. 
Finally we have $g=\alpha (f-g')+\beta g'$, thus $g - \alpha f = (\beta - \alpha) g'$, and therefore $g' \in [0,f,g]$. By convexity $g' \in \cV_{\max}(x)$, by (Stability) $g' \in \cV(x)$, which contradicts our assumption that $f,g$ are trigonometrically consecutive in $\cV(x)$. This concludes the proof.
\end{proof}

\begin{Proposition}[Characterization of the smallest stencils]
\label{prop:MinimalStencils}
For all $x \in X$, define $\cV_{\min}(x)$ as the collection of all $e \in \cV_{\max} (x)$ which have none or just one parent in $\cV_{\max}(x)$ (this includes all unit vectors in $\cV_{\max}(x)$). 
Then $\cV_{\min} := (\cV_{\min}(x))_{x \in X}$ is a family of stencils, which is contained in any other family of stencils. 
\end{Proposition}

\begin{proof}
Property (Stability) of stencils. Consider $x\in X$ and $e \in \cV_{\min}(x)$. Assume for contradiction that $e$ has one parent $e' \in \cV_{\max}(x)$ which is not an element of $\cV_{\min}(x)$. Hence $e'$ has two parents $f,g \in \cV_{\max}(x)$. 
By Corollary \ref{corol:ChildrenInCone} we have $f \prec e \prec g$, thus by Lemma \ref{lem:ConsecutiveParents} (ii.a) the \emph{two} parents of $e$ belong to the triangle $[e,f,g]$, hence also to $\cV_{\max}(x)$ by convexity of $\Omega$. This contradicts our assumption that $e \in \cV_{\min}(x)$.

Property (Visibility). 
We consider $e\in \cV_{\max}(x)$, and prove by induction on the norm $\|e\|$ that $e\in K := \Cone(\cV_{\min}(x))$. If $\|e\|=1$ or if none of just one parent of $e$ belongs to $\cV_{\max}(x)$, then $e \in \cV_{\min}(x) \subset K$. If both parents $f,g$ of $e$ belong to $\cV_{\max}(x)$, then by induction $f,g \in K$, and by additivity $f+g \in K$, which concludes the proof. 

Minimality for inclusion of $\cV_{\min}$. Let $\cV$ be a family of stencils, let $x \in X$, $e \in \cV_{\min}(x)$, and let us assume for contradiction that $e \notin \cV(x)$.
By property (Visibility) of stencils, the vector $e$ belongs to the cone generated by two elements $f,g \in \cV(x)$, which can be chosen trigonometrically consecutive in $\cV(x)$. By lemma \ref{lem:ConsecutiveBasis}, $(f,g)$ is a direct acute basis of $\Z^2$. By Lemma \ref{lem:ConsecutiveParents} (ii.a) the parents of $e$ belong to the triangle $[e,f,g]$, hence to $\cV_{\max}(x)$ by convexity of $\Omega$, which contradicts the definition of $\cV_{\min}(x)$.
\end{proof}

\begin{Proposition}[Structure of candidates for refinement]
\label{prop:Trigo}
Let $\cV$ be a family of stencils on $\OmegaD$, and let $x \in \OmegaD$. Then the parents $f,g$, of any candidate for refinement $e \in \Candidates(x)$, are consecutive elements of $\cV(x)$ in trigonometric order.
\end{Proposition}

\begin{proof}
Since $e \notin \cV(x)$, there exists by (Visibility) two elements $f,g \in \cV(x)$ such that $f\prec e \prec g$, and which we can choose trigonometrically consecutive in $\cV(x)$. By Lemma \ref{lem:ConsecutiveBasis}, $(f,g)$ is a direct acute basis. By Lemma \ref{lem:ConsecutiveParents} (ii.a) the parents $f',g'$ of $e$ between satisfy $f \preceq f' \prec g' \preceq g$. Recalling that $f',g' \in \cV(x)$, by definition of $\hat \cV(x)$, we obtain $f=f'$ and $g=g'$ which concludes the proof.
\end{proof}

\subsection{Combining and intersecting constraints}
\label{sec:CombiningIntersecting}

The  following characterization of the cones $\Conv(\cV)$ implies the announced hierarchy properties. 

\begin{Proposition}
\label{prop:ConvLComplement}
For any family $\cV$ of stencils on $X$ one has
\begin{equation}
\Conv(\cV) = \{u \in \Conv(X); \, P_x^e(u) \geq 0 \text{ for all } x \in X, \, e \in \cV(x)\}. 
\end{equation}
\end{Proposition}


Before turning to the proof of this proposition, we use it to conclude the proof of Theorem \ref{th:Hierarchy}. The sub-cone $\Conv(\cV)$, of $\Conv(\Omega)$, is characterized by the non-negativity of a family of linear forms indexed by $\cV_{\max} \sm \cV$, with the convention \eqref{StencilIdentification}. 
Observing that 
\begin{equation*}
\cV_{\max} \sm (\cV \cup \cV') = (\cV_{\max} \sm \cV) \cap (\cV_{\max} \sm \cV'), \qquad \cV_{\max} \sm (\cV \cap \cV') = (\cV_{\max} \sm \cV) \cup (\cV_{\max} \sm \cV'),
\end{equation*}
we find that $\Conv(\cV \cup \cV')$ is characterized, as a subset of $\Conv(\OmegaD)$, by the intersection of the families of constraints defining $\Conv(\cV)$ and $\Conv( \cV')$, while $\Conv(\cV\cap \cV')$ is defined by their union. Hence we conclude as announced
\begin{equation*}
\Conv(\cV \cup \cV') \supseteq \Conv(\cV) \cup \Conv( \cV'), \qquad \Conv(\cV \cap \cV') = \Conv(\cV) \cap \Conv( \cV').
\end{equation*}


\begin{proof}[Proof of Proposition \ref{prop:ConvLComplement}]
We proceed by decreasing induction on the cardinality $\#(\cV)$.

Initialization. If $\#(\cV) = \#(\cV_{\max})$, then $\cV = \cV_{\max}$, and therefore $\Conv(\cV) = \Conv(\cV_{\max}) = \Conv(\OmegaD)$ and $\cV_{\max} \sm \cV = \emptyset$. The result follows.

Induction. Assume that $\#(\cV) < \#(\cV_{\max})$, thus $\cV \subsetneq \cV_{\max}$. Let $x\in \OmegaD$ and $e \in \cV_{\max}(x) \sm \cV(x)$ be such that $\|e\|$ is minimal. 
Since $e \notin \cV_{\min}(x) \subset \cV(x)$, the two parents $f,g$ of $e$ belong to $\cV_{\max}(x)$. Since $\|e\| > \max \{\|f\|,\|g\|\}$, and by minimality of the norm of $e$, we have $f,g \in \cV(x)$. Hence $e$ is a candidate for refinement: $e\in \hat \cV(x)$.

Consider the extended stencils $\cV'$ defined by $\cV'(x) := \cV(x) \cup \{e\}$, and $\cV'(y) := \cV(y)$ for all $y \in \OmegaD\sm \{x\}$. 
Let $\cL$ and $\cL'$ be the collections of linear forms enumerated in Definition \ref{def:Cones}, which non-negativity respectively defines the cones $\Conv(\cV)$ and $\Conv(\cV')$ as subsets of $\cF(X)$. Let also $\cL_0 := \cL \cap \cL'$.
Since $e\in \hat \cV(x)$ we have  $\cL = \cL_0 \cup \{P_x^e\}$. Using Proposition \ref{prop:Trigo} we obtain $\hat \cV'(x) \sm \hat \cV(x) \subset \{e+f,e+g\}$, hence $\cL'$ is the union of $\cL_0$ and of those of the following constraints which are supported on $\OmegaD$:
\begin{equation}
\label{AdditionalForms}
S_x^e, \ T_x^e, \ P_x^{e+f}, \ P_x^{e+g}.
\end{equation}
We next show that $\Cone(\cL) = \Cone(\cL' \cup \{P_x^e\})$, by expressing the linear forms \eqref{AdditionalForms} in terms of the elements of $\cL$.
\begin{itemize}
\item If $S_x^e$ is supported on $X$, then $-e \in \cV_{\max}(x)$. Assuming that $-e \in \cV(x)$, we obtain $S_x^e = S_x^{-e} \in \cL$. On the other hand, assuming that $-e \notin \cV(x)$, we obtain $-f,-g \in \cV_{\max}(x)$ by Proposition \ref{prop:MinimalStencils}, since otherwise $-e \in \cV_{\min}(x) \subset \cV(x)$. Therefore $S_x^f, S_x^g$ are supported on $X$, hence they belong to $\cL$. By minimality of the norm of $e$, we have $-f,-g \in \cV(x)$, hence $-e \in \hat \cV(x)$ and therefore $P_x^{-e} \in \cL$. As a result $S_x^e = P_x^e+P_x^{-e}+S_x^f+S_x^g \in \Cone(\cL)$.
\item If $T_x^e$ is supported on $X$, then $-f, -g \in \cV_{\max}(x)$. Therefore $S_x^f, S_x^g$ are supported on $X$, hence they belong to $\cL$. As a result $T_x^e = P_x^e+S_x^f+S_x^g \in \Cone(\cL)$.
\item If $P_x^{e+f}$ is supported on $X$, then $x+e+f \in \OmegaD$, thus $f\in \cV_{\max}(x+e)$ and therefore $f\in \cV(x+e)$ by minimality of $\|e\|$. The linear form $S_{x+e}^f$ belongs to $\cL$, since it has support $\{x+g,x+e,x+e+f\} \subset \OmegaD$.
Observing that the parents of $e+f$ are $e$ and $f$, we find that $P_x^{e+f} = P_x^e + S_{x+e}^f \in \Cone(\cL)$. The case of $P_x^{e+g}$ is similar.
\end{itemize}
Denoting by $K^*$ the dual cone of a cone $K$, we obtain
\begin{equation*}
\Conv(\cV) = 
\Cone(\cL)^* = 
\Cone(\cL' \cup \{P_x^e\})^* = 
\{u \in \Conv(\cV'); \, P_x^e(u) \geq 0\}.
\end{equation*}
Applying the induction hypothesis to $\cV'$, we conclude the proof.
\end{proof}

\section{Stencils and triangulations}
\label{sec:Delaunay}
Using the interplay between stencils $\cV$ and triangulations $\cT$, we prove Proposition \ref{prop:WorstCase} and Theorems \ref{th:Decomp}, \ref{th:Delaunay}. By convention, all stencils $\cV$ are on $X$, and all triangulations $\cT$ have $X$ as vertices and cover $\Hull(X)$.

\subsection{Minimal stencils containing a triangulation}
We characterize in Proposition \ref{prop:VofT} the minimal stencils $\cV$ containing a triangulation, in the sense of Definition \ref{def:TinV}, and we estimate their cardinality, proving Proposition \ref{prop:WorstCase}. In the way, we establish in Proposition \ref{prop:UinConvV} ``half'' (one inclusion) of the decomposition of $\Conv(\cV)$ announced in Theorem \ref{th:Decomp}.

\begin{Lemma}
\label{lem:EdgeIneq}
Let $\cT$ be a triangulation, and let $u \in \Conv(\cT)$.
Let $p,q,r \in X$. Assume that $[p,q]$ is an edge of $\cT$, and that $s := p+q-r \in X$. Then 
$
u(r)+u(s) \geq u(p)+u(q).
$
\end{Lemma}

\begin{proof}
The interpolating function $U := \interp_\cT u$ is convex on $\Hull(X)$, and linear on the edge $[p, q]$. Introducing the edge midpoint $m := (p+q)/2 = (r+s)/2$ we obtain $u(p)+u(q) = 2 U(m) \leq u(r)+u(s)$, as announced.
\end{proof}

The inequalities $u(r)+u(s) \geq u(p)+u(q)$ identified in the previous lemma are closely tied with the linear constraints $P_x^e$, since $[p,q,r,s]$ is a parallelogram, and as shown in the next lemma. The set $\Anc(e)$ of ancestors of an irreducible vector $e\in \Z^2$ was introduced Definition \ref{def:Ancestors}.

\begin{Lemma}
\label{lem:PHierarchy}
Let $e \in \Z^2$ be irreducible, with $\|e\|> 1$, and let $(f,g)$ be a direct basis such that $f \prec e \prec g$ and $e \notin \Anc(f) \cup \Anc(g)$. Let $x \in X$ be such that $f,g,f+g,e \in \cV_{\max}(x)$. If $u \in \Conv(X)$ satisfies $u(x)+u(x+f+g) \geq (x+f)+u(x+g)$, then $P_x^e(u) \geq 0$.
\end{Lemma}

\begin{proof}
Without loss of generality, up to adding a global affine map to $u$, we may assume that $u(x+e) = u(x+f)=u(x+g)=0$. 
Denoting by $f',g'$ the parents of $e$, we have by Lemma \ref{lem:ConsecutiveParents} (ii.b) $f',g', f+g \in [e,f,g]$, hence by convexity $u(x+f'), u(x+g'), u(x+f+g) \leq 0$.
Our hypothesis implies $u(x) \geq -u(x+f+g) \geq 0$, therefore $P_x^e(u) = u(x)-u(x+f')-u(x+g')+u(x+e) \geq 0$.
\end{proof}

\begin{Proposition}
\label{prop:UinConvV}
If a triangulation $\cT$, and stencils $\cV$, satisfy $\cT \prec \cV$, then $\Conv(\cT) \subset \Conv(\cV)$.
\end{Proposition}

\begin{proof}
The inequalities $S_x^e(u) \geq 0$, and $T_x^e(u) \geq 0$, for $x \in X$, $e \in \cV(x)$, hold by convexity of $u$. We thus consider an arbitrary refinement candidate $e \in \hat \cV(x)$, $x\in X$, and establish below that $P_x^e(u) \geq 0$. 

Since the triangulation $\cT$ covers $\Hull(X)$, there exists a triangle $T \in \cT$, containing $x$, and such that $e \in \Cone(T-x)$. 
Since $\cT \prec \cV$ and $e \notin \cV(x)$, the segment $[x,x+e]$ is not an edge of $\cT$. Denoting the vertices of $T$ by $[x,x+f,x+g]$ we have $f \prec e \prec g$. Since $e \in \hat \cV(x)$, one has $e \notin (\Anc(f) \cup \Anc(g)) \subset \cV(x)$. Applying Lemma \ref{lem:EdgeIneq} to the edge $[x+f,x+g]$ we obtain $u(x)+u(x+f+g) \geq u(x+f)+u(x+g)$. Finally, Lemma \ref{lem:PHierarchy} implies $P_x^e(u) \geq 0$ as announced. 
\end{proof}

\begin{Proposition}
\label{prop:VofT}
Let $\cT$ be a triangulation, and for all $x \in X$ let $V_x$ be the collection of all $e\in \Z^2$ such that $[x,x+e]$ is an edge of $T$. The minimal family of stencils satisfying $\cT \prec \cV$ is given by 
\begin{equation*}
\cV(x) := \cV_{\max}(x) \cap \bigcup_{e \in V_x} \Anc(e).
\end{equation*}
\end{Proposition}
\begin{proof}
The family of sets $\cV$ satisfies the (Stability) property by construction. Since the triangulation $\cT$ covers $\Hull(X)$, the sets $(V_x)_{x\in X}$ satisfy the (Visibility) property, hence also the larger sets $\cV(x) \supset V_x$. 

Minimality. Consider \emph{arbitrary} stencils $\cV$ such that $\cT \prec \cV$.
Let also $x \in X$, $e \in V_x$, $e' \in \cV_{\max}(x) \cap \Anc(e)$, and let us assume for contradiction that $e' \notin \cV(x)$. By property (Visibility) there exists $f,g$, trigonometrically consecutive elements of $\cV(x)$, such that $f\prec e \prec g$ (where $\prec$ refers to the cyclic trigonometric order, see Definition \ref{def:CyclicOrder}). 
By Lemma \ref{lem:ConsecutiveBasis}, $(f,g)$ is a direct acute basis of $\Z^2$. By Corollary \ref{corol:ChildrenInCone}, and an immediate induction argument, we have $f \prec e \prec g$, hence $e \notin \cV(x)$, which contradicts our assumption that $\cT \prec \cV$.
\end{proof}

Given a triangulation $\cT$, our next objective is to estimate the cardinality of the minimal stencils $\cV$ such that $\cT \prec \cV$. We begin by counting the ancestors of an irreducible vector.

\begin{Lemma}
\label{lem:AncestorsSimple}
\begin{enumerate}
\item
Let $(f,g)$ be an acute basis of $\Z^2$. Then either (i) $f$ is a parent of $g$, (ii) $g$ is a parent of $f$, or (iii) $\|f\|=\|g\|=1$.
\item 
For any irreducible $e \in \Z^2$ one has $\#(\Anc(e)) \leq \|e\|_\infty+2$, where $\|(\alpha,\beta)\|_\infty := \max\{|\alpha|, |\beta|\}$.
\end{enumerate}
\end{Lemma}

\begin{proof}
Point 1. 
If $\<f,g\> > 0$, then applying Corollary \ref{corol:SumParent} to the non-acute basis $(f,-g)$ we find that either (i) $(f,g-f)$ are the parents of $g$, or (ii) $(f-g,g)$ are the parents of $f$. On the other hand if $\<f,g\>=0$, then $1 = |\det(f,g)| = \|f\| \|g\|$, hence (iii) $\|f\|=\|g\|=1$. 

Before proving Point 2, we introduce the cone $K$ generated by $(1,0),(1,1)$, so that $\|(\alpha, \beta)\|_\infty = \alpha$ for any $(\alpha,\beta) \in K$. If $e\in \Z^2$ irreducible belongs to the interior of $K$, then its parents $f,g \in K$, and we have $\|e\|_\infty = \|f\|_\infty+\|g\|_\infty$.

Point 2 is proved by induction on $\|e\|_\infty$. It is immediate if $\|e\|_\infty = 1$, hence we may assume that $\|e\|_\infty \geq 2$, and denote its parents by $f,g$. 
We have $\|e\|_\infty = \|f\|_\infty+ \|g\|_\infty$, since without loss of generality we may assume that $e\in K$.
Applying Point 1 we find that either (i) $\Anc(e) = \Anc(g) \cup \{e\}$, (ii) $\Anc(e) = \Anc(f) \cup \{e\}$, or (iii) $\|f\|=\|g\|=1$, so that $\|e\|_\infty=1$, a case which we have excluded. Thus $\#\Anc(e) \leq \max\{\# \Anc(f), \#\Anc(g)\} +1$, which implies the announced result by induction.
\end{proof}

\begin{Proposition}
\label{prop:VTCard}
Let $\cT$ be a triangulation, and let $\cV$ be the minimal family of stencils such that $\cT \prec \cV$. Then $\#(\cV) \leq 6 (N-2) (\diam(\OmegaC) + 2)$, with $N:=\#(X)$. A sharper estimate holds for (standard) Delaunay triangulations: $\#(\cV) \leq 6 (N-2)$, 
\end{Proposition}
\begin{proof}
Let $E,F$ be respectively the number of edges and faces of $\cT$, where faces refer to both triangles and the infinite exterior face. By Euler's theorem,  $N-E+F=2$. Since each edge is shared by two faces, and each face has at least three edges, one gets $2 E \geq 3 F$, hence $E \leq 3 (N-2)$, and therefore, with the notation $V_x$ of Proposition \ref{prop:VofT},
\begin{equation}
\label{eq:Euler}
\sum_{x\in X} \#(V_x) = 2 E \leq 6(N-2).
\end{equation}
Combining lemma \ref{lem:AncestorsSimple}, Proposition \ref{prop:VofT}, and observing that any edge $[x,x+e]$ of $\cT$ satisfies $\|e\|_\infty \leq \diam(\Omega)$, we obtain $\#(\cV(x)) \leq \#(V_x) (\diam(\Omega)+2)$, which in combination with \eqref{eq:Euler} implies the first estimate on $\#(\cV)$.

In the case of a Delaunay triangulation, we claim that $\cV(x) = V_x$. Indeed, consider an edge $[x,x+e]$ of $\cT$, and a parent $f\in \cV_{\max}(x)$ of $e$. Since $\cT$ covers $\Hull(X)$, it contains a triangle $[x,x+e,x+f']$ with $f$ and $f'$ on the same side of the edge $[x,x+e]$. Thus the determinants $\det(e,f)$ and $\det(e,f')$ have the same sign, and therefore the same value since their magnitude is $1$. As a result $f'=k e + f$, for some integer $k$. Since $\cT$ is Delaunay, the point $x+f$ is outside of the circumcircle of $[x,x+e,x+f']$. This property is equivalent to the non-positivity of the following determinant, called the in-circle predicate: assuming without loss of generality that $\det(e,f)=1$ so that the vertices $(x,x+e,x+f')$ are in trigonometric order
\begin{align}
\nonumber
\det
\left(
\begin{array}{ccc}
e_1 & e_2 & \|e\|^2\\
f_1 & f_2 & \|f\|^2 \\
k e_1+f_1 & k e_2+ f_1 & \| k e + f\|^2
\end{array}
\right)
&=
\det
\left(
\begin{array}{ccc}
e_1 & e_2 & \|e\|^2\\
f_1 & f_2 & \|f\|^2\\
0 & 0 & \| k e + f\|^2 - k \|e\|^2 - \|f\|^2
\end{array}
\right),\\
\nonumber
& = \| k e + f\|^2 - k \|e\|^2 - \|f\|^2,\\
\label{inCircle}
& = k(k-1)\|e\|^2 + 2 k \<e,f\>,
\end{align}
where we denoted $e=(e_1,e_2)$, $f=(f_1,f_2)$.
Observing that $0 < \<e,f\> \leq \|e\| \|f\| < \|e\|^2$, we find that \eqref{inCircle} is non-positive only for $k=0$. Thus $f=f'$, hence $f \in V_x$, and therefore $\cV(x)=V_x$ as announced. Finally, the announced estimate of $\#(\cV)$ immediately follows from \eqref{eq:Euler}.
\end{proof}

Let us conclude the proof of Proposition \ref{prop:WorstCase}.
Let $u \in \Conv(X)$, and let $\cV_u$ be the minimal stencils such that $u \in \Conv(\cV_u)$. Let $\cT$ be an $u$-Delaunay triangulation, and let $\cV_\cT$ be the minimal stencils such that $\cT \prec \cV$. By Proposition \ref{prop:UinConvV} we have $u \in \Conv(\cT) \subset \Conv(\cV_\cT)$, hence $\cV_u \subset \cV_\cT$. Estimating $\#(\cV_T)$ with Proposition \ref{prop:VTCard}, we obtain as announced $\#(\cV_u) \leq 6 (N-2) (\diam(\Omega)+2)$.

\subsection{Decomposition of the cone $\Conv(\cV)$, and edge-flipping distances}

We conclude in this section the proof of Theorem \ref{th:Decomp}, 
and establish the complexity result Theorem \ref{th:Delaunay} on the edge-flipping generation of $u$-Delaunay triangulations.

\begin{Definition}
We say (abusively) that a discrete map $u : X \to \R$ is generic iff, for all $x\in X$ and all $e \in \Z^2$ such that the linear form $P_x^e$ is supported on $X$, one has $P_x^e(u) \neq 0$.
\end{Definition}

Generic elements are dense in $\Conv(X)$, since this set is convex, has non-empty interior, and since non-generic elements lie on a union of hyperplanes. The quadratic function $q(x) := \frac 1 2 \|x\|^2$ is not generic however, since choosing $e=(1,1)$ one gets $P_x^e(q) = 0$. 

\begin{Lemma}
\label{lem:genericTV}
Consider stencils $\cV$, a generic $u \in \Conv(\cV)$, and an $u$-Delaunay triangulation $\cT$. Then $\cT \prec \cV$.
\end{Lemma}

\begin{proof}

Consider an edge $[x,x+e]$ of $\cT$. If the linear form $P_x^e$ is not supported on $\cT$, then $e \in \cV_{\min}(x) \subset \cV(x)$ by Proposition \ref{prop:MinimalStencils}. On the other hand if $P_x^e$ is supported on $X$, then $P_x^e(u) \leq 0$ by Lemma \ref{lem:EdgeIneq}. By genericity of $u$, we have $P_x^e(u) < 0$, hence $e \in \cV(x)$ by Proposition \ref{prop:ConvLComplement}. This concludes the proof.
\end{proof}

We established in Proposition \ref{prop:UinConvV} that 
$\Conv(\cV) \supset \cup_{\cT \prec \cV} \Conv(\cT)$.
The next corollary, stating the reverse inclusion, concludes the proof of Theorem \ref{th:Decomp}.

\begin{Corollary}
\label{corol:DecompSub}
If $\Conv(\cV)$ has a non-empty interior, then 
$\Conv(\cV) \subset \cup_{\cT \prec \cV} \Conv(\cT)$.
\end{Corollary}

\begin{proof}
The set $K := \cup_{\cT \prec \cV} \Conv(\cT)$ contains all generic elements of $\Conv(\cV)$, by Lemma \ref{lem:genericTV}. Observing that $K$ is closed, and recalling that generic elements are dense in $\Conv(\cV)$, we obtain the announced inclusion.
\end{proof}

The next lemma characterizes the obstructions to the convexity of the piecewise linear interpolant $\interp_\cT u$ of a convex function $u \in \Conv(X)$ on a triangulation $\cT$. See also Figure \ref{fig:Flipping} (right).
\begin{Lemma}
\label{lem:BadEdge}
Consider $u \in \Conv(X)$, and a triangulation $\cT$ which is \emph{not} $u$-Delaunay. Then there exists $x\in X$, and a direct basis $(f,g)$ of $\Z^2$, such that the triangles $[x,x+f,x+g]$ and $[x+f,x+g,x+f+g]$ belong to $\cT$, and satisfy $u(x)+u(x+f+g) < u(x+f)+u(x+g)$. 
\end{Lemma}

\begin{proof}
Since convexity is a local property, there exists two triangles $T,T' \in \cT$, sharing an edge, such that the interpolant $\interp_\cT u$ is not convex on $T \cup T'$ (i.e.\ convexity fails on the edge by $T$ and $T'$). Up to a translation of the domain, we may assume that $T=[0,f,g]$ and $T'=[f,g,e]$, for some $e,f,g \in \Z^2$. The pair $(f,g)$ is a basis of $\Z^2$ because the triangle $T$ contains no point of $\Z^2$ except its vertices; up to exchanging $f$ and $g$ we may assume that it is a direct basis. 
Up to adding an affine function 
to $u$, we may assume that $u$  vanishes at the vertices $0,f,g$ of $T$.

If $f$ lies in the triangle $[0,e,g]$, then since $u\in \Conv(X)$, and recalling that 
$u(0)=u(f)=u(g)=0$, we obtain $u(e) \geq 0$. 
This implies that $\interp_\cT u$ is convex on $T \cup T'$, which contradicts our assumption.
Likewise $g \notin [0,f,e]$, thus $f \prec e \prec g$ and therefore $\det(f,e)> 0$ and $\det(e,g)>0$. We next observe that 
\begin{equation*}
\det(f,g) + \det(g-e,f-e) = \det(f,e)+\det(e,g).
\end{equation*}
The four members of this equation are integers, the two left being equal to $2|T| = 2|T'| = 1$, and the two right being positive. Hence $\det(f,e)=\det(e,g)=1$, and therefore $e=f+g$ as announced. From this point, the inequality $u(x)+u(x+f+g) < u(x+f)+u(x+g)$ is easily checked to be equivalent to the non-convexity of $\interp_\cT u$ on $T \cup T'$.
\end{proof}

\begin{Proposition}
\label{prop:Flipping}
Consider stencils $\cV$, a triangulation $\cT \prec \cV$, and $u \in \Conv(\cV)$. Define a sequence of triangulations $\cT_0 := \cT,\ \cT_1,\ \cT_2 \cdots$ as follows: if $\cT_i$ is $u$-Delaunay, then the sequence ends, otherwise $\cT_{i+1}$ is obtained by flipping an arbitrary edge of $\cT_i$ satisfying Lemma \ref{lem:BadEdge}. Then the sequence is finite, contains at most $\#(\cV)$ elements, and $\cT_i \prec \cV$ for all $0 \leq i \leq n$.
\end{Proposition}

\begin{proof}
Proof that $\cT_i \prec \cV$, by induction on $i \geq 0$. Initialization: $\cT_0 := \cT \prec \cV$ by assumption. Induction: adopting the notations of Lemma \ref{lem:BadEdge}, the ``flipped'' edge $[x+f,x+g]$ of $\cT_i$ is replaced with $[x,x+e]$ in $\cT_{i+1}$, with $e:=f+g$. We only need to check that $e \in \cV(x)$, and for that purpose we distinguish two cases. 
If the basis $(f,g)$ is acute, then $f,g$ are the parents of $e$, and we have $P_x^e(u) < 0$ by Lemma \ref{lem:BadEdge}. This implies $e \in \cV(x)$ by Proposition \ref{prop:ConvLComplement}.
On the other hand, if the basis $(f,g)$ is not acute, then by Corollary \ref{corol:SumParent} the vector $e$ is a parent of either $f$ or $g$, thus $e \in \cV(x)$ by property (Stability) of stencils.

Bound on the number $n$ of edge flips. For all $0 \leq i < n$ one has $\interp_{\cT_{i+1}} u \leq \interp_{\cT_i} u$ on $\Hull(X)$, and this inequality is strict at the common midpoint of the flipped edges $[x_i+f_i, x_i+g_i]$ and $[x_i,x_i+e_i]$, with the above conventions. 
Hence the edge $[x_i,x_i+e_i]$ appears in the triangulation $\cT_{i+1}$ but not in any of the $\cT_j$, for all $0 \leq j \leq i$. It follows that $i \mapsto (x_i,e_i)$ is injective, and since $e_i \in \cV(x_i)$ this implies $n \leq \#(\cV)$.
\end{proof}

We finally prove Theorem \ref{th:Delaunay}. Consider a Delaunay triangulation $\cT$, and the minimal stencils $\cV_{\cT}$ such that $\cT \prec \cV$. Let also $u \in \Conv(X)$, and let $\cV_u$ be the minimal stencils such that $u \in \Conv (\cV_u)$. Then, by Proposition \ref{prop:Flipping}, $\cT$ can be transformed into an $u$-Delaunay triangulation via $\#(\cV_\cT \cup \cV_u)$ edge flips. Furthermore $\#(\cV_\cT) = \cO(\#(X))$ by Proposition \ref{prop:VTCard} and $\#(\cV_u) \geq \#(X)$, as follows e.g.\ from property (Visibility) of stencils. Thus $\#(\cV_\cT \cup \cV_u) = \cO(\#(\cV_u))$, and the result follows.


\section{Average case estimate of the cardinality of minimal stencils}
\label{sec:Avg}

The minimal stencils $\cV$, such that the cone $\Conv(\cV)$ contains a given discrete convex map, admit a simple characterization described in the following proposition. 
\begin{Proposition}
\label{prop:MinimalStencilsCharacterization}
Let $u\in \Conv(X)$, and let $\cV$ be the minimal stencils on $X$ such that $u \in \Conv(\cV)$. 
For any $x \in X$, and any irreducible $e \in \Z^2$ with $\|e\|>1$, one has: 
\begin{center}
$e \in \cV(x)$ $\Leftrightarrow$ ($P_x^e$ is not supported on $X$, or $P_x^e(u) < 0$).
\end{center}
\end{Proposition}

\begin{proof}
Proof of implication $\Leftarrow$. If $P_x^e$ is not supported on $X$, then 
$e \in \cV_{\min}(x) \subset \cV(x)$ by Proposition \ref{prop:MinimalStencils}.  On the other hand if $P_x^e(u) < 0$, then $e \in \cV(x)$ by Proposition \ref{prop:ConvLComplement}. 

Proof of implication $\Rightarrow$. Consider $x\in X$, $e \in \cV(x)$, with $\|e\| > 1$, and such that $P_x^e$ is supported on X. Assume for contradiction that $P_x^e(u) \geq 0$, and denote by $f,g$ the parents of $e$. Let $E := \{e' \in \cV(x); \, f \prec e' \prec g\}$. The parents of any $e' \in E$ belong to $E \cup \{f,g\}$ by Lemma \ref{lem:ConsecutiveParents} (ii), and one has $P_z^{e'} (u) \geq 0$ by Lemma \ref{lem:PHierarchy}. Defining new stencils by $\cV'(x) := \cV(x) \sm E$, and $\cV'(y) := \cV(y)$ for $y \neq x$, we contradict the minimality of $\cV$.
\end{proof}

The rest of this section is devoted to the proof of Theorem \ref{th:Avg}, and for that purpose we consider the rotated and translated grids $X_\theta^\xi$, defined in \eqref{def:XThetaXi}. 
For simplicity, but without loss of generality, we assume a unit grid scale $h:=1$.
For each rotation angle $\theta \in \R$, and each offset $\xi \in \R^2$, we introduce an affine transform $A_\theta^\xi$: for all $x\in \R^2$
\begin{equation*}
A_\theta^\xi (x) := R_\theta (\xi+x).
\end{equation*}
For any set $E \subset \R^2$, and any affine transform $A$, we denote $A(E) := \{A(e); \, e\in E\}$. For instance, the displaced grids \eqref{def:XThetaXi} are given by $X_\theta^\xi := \Omega \cap A_\theta^\xi(\Z^2)$.

The maximal stencils on the grid $X_\theta^\xi$ are defined by: for all $x\in X_\theta^\xi$
\begin{equation*}
\cV_{\max}^{\theta,\xi}(x) := \{e \in \Z^2 \text{ irreducible}; \ x+ R_\theta e\in X_\theta^\xi\}.
\end{equation*}
A family $\cV_\theta^\xi$ of stencils on $X_\theta^\xi$ is a collection of sets $\cV_\theta^\xi(x) \subset \cV_{\max}^{\theta,\xi}(x)$, $x \in X_\theta^\xi$ which satisfies the usual (Stability) and (Visibility) properties of Definition \ref{def:Stencil} (replacing, obviously, instances of $\cV_{\max}$ with $\cV_{\max}^{\theta, \xi}$).
For $x \in X_\theta^\xi$, and $e \in \cV_{\max}^{\theta,\xi}$ we consider the linear forms $S_{x, \theta}^e(u) := u(x+R_\theta e)-2 u(x) +u(x-R_\theta e)$, and likewise $T_{x, \theta}^e$, $P_{x, \theta}^e$, which are used to define cones $\Conv(\cV_\theta^\xi) \subset \Conv(X_\theta^\xi)$.
In a nutshell, when embedding a stencil element $e \in \cV_\theta^\xi (x) \subset \Z^2$, where $x\in X_\theta^\xi$, into the physical domain $\Omega$ (e.g.\ considering $x+R_\theta e \in X_\theta^\xi$), one should never forget to apply the rotation $R_\theta$.

Consistently with the notations of Theorem \ref{th:Avg}, we consider a fixed convex map $U\in \Conv(\Omega)$, and study the smallest stencils $\cV_\theta^\xi \subset \Z^2$ on $X_\theta^\xi$ such that the restriction of $U$ to $X_\theta^\xi$ belongs to $\Conv(\cV_\theta^\xi)$.
The \emph{midpoints} $m = x+R_\theta e/2$ of ``stencil edges'' $[x,x+R_\theta e]$, $x\in X_\theta^\xi$, $e\in \cV_\theta^\xi(x)$, play a central role in our proof. 
\begin{Definition}
For any $m \in \Omega$, and any irreducible $e \in \Z^2$, let 
\begin{equation*}
\Lambda_m^e := \{(\theta, \xi)\in [0, \pi/2[ \times [0,1[^2; \, m = x+R_\theta e/2, \text{ for some } x \in X_\theta^\xi \text{ such that } e \in \cV_\theta^\xi(x)\}.
\end{equation*}
\end{Definition}
We introduce offsetted grids, of points with half-integer coordinates
\begin{equation*}
\cZ : =(\textstyle {\frac 1 2},0) +\Z^2, \qquad \cZ' : =(0,\textstyle {\frac 1 2}) +\Z^2, \qquad \cZ'' : =(\textstyle {\frac 1 2},\textstyle {\frac 1 2}) +\Z^2.
\end{equation*}
For any $x,y \in \Z^2$ with $x-y$ irreducible, the midpoint $(x+y)/2$ of the segment $[x,y]$ belongs to the disjoint union $\cZ \sqcup \cZ' \sqcup \cZ''$. 
\begin{Lemma}
\label{lem:HalfGrid}
For any $m \in \Omega$ and any $(\theta,\xi) \in \Lambda_m^e$, one has $m \in A_\theta^\xi (\cZ \sqcup \cZ' \sqcup \cZ'')$. 
\end{Lemma}
\begin{proof}
Let $x \in X_\theta^\xi$, and $e \in \cV_\theta^\xi(x)$, be such that $m = x+R_\theta e/2$. 
Observing that the coordinates of $e$ are not both even, since $e$ is irreducible, we obtain $e/2 \in \cZ \sqcup \cZ' \sqcup \cZ''$. Adding $R_\theta (e/2)$ to $x \in A_\theta^\xi (\Z^2)$ yields as announced a point $m\in A_\theta^\xi (\cZ \sqcup \cZ' \sqcup \cZ'')$.
\end{proof}
For any point $m \in \R^2$, and any angle $\theta$, there exists exactly one offset $\xi\in [0,1[^2$ such that $m \in A_\theta^\xi(\cZ)$; and likewise for $\cZ'$, $\cZ''$.
Hence the set $\Lambda_m^e$ contains redundant information, which 
motivates the following definition: for any $m \in \Omega$, and any irreducible $e \in \Z^2$
\begin{equation}
\label{def:ThetaME}
\Theta_m^e := \{ \theta \in [0, \pi/2[; \, \exists \xi \in [0,1[^2, \, (\theta, \xi) \in \Lambda_m^e \text{ and } m \in A_\theta^\xi (\cZ) \}.
\end{equation}
and similarly we define $\Theta_m'^e$, $\Theta_m''^e$, by replacing $\cZ$ with $\cZ'$, $\cZ''$ respectively in \eqref{def:ThetaME}. By convention, $\Theta_m^e = \Theta_m'^e=\Theta_m''^e=\emptyset$ for \emph{non} irreducible vectors $e \in \Z^2$. The following lemma accounts in analytical terms for a simple combinatorial identity: one can count stencil edges by looking at their endpoints or their midpoints.

\begin{Lemma}
The following integrals are equal: 
\begin{equation}
\label{eq:MeasEqual}
\int_{[0,1]^2} \int_0^{\frac \pi 2} \#(\cV_\theta^\xi) \, d \theta d \xi = \sum_{e \in \Z^2} \int_{m \in \Omega} (|\Theta_m^e| + |\Theta_m'^e| + |\Theta_m''^e|) \,dm,
\end{equation}
where $|\Theta|$ denotes the Lebesgue measure of a Borel set $\Theta \subset \R$.
\end{Lemma}

\begin{proof}
Consider $m \in \Omega$, $e \in \Z^2$, and $\theta \in \Theta_m^e$. Then there exists a unique $\xi \in [0,1[^2$ such that $m \in A_\theta^\xi(\cZ)$. This uniquely determines the point $x := m-\frac 1 2 R_\theta e \in X_\theta^\xi$ such that $e \in \cV_\theta^\xi(x)$. 
Likewise for $\Theta_m'^e$, $\Theta_m''^e$. 
Conversely, the data of $\theta$, $\xi$, $x\in X_\theta^\xi$ and $e \in \cV_\theta^\xi(x)$ uniquely determines $m := x+R_\theta e/2$, and also by Lemma \ref{lem:HalfGrid} a unique set among $\Theta_m^e$, $\Theta_m'^e$, $\Theta_m''^e$ containing $\theta$. As a result the left and right hand side of \eqref{eq:MeasEqual} are just two different expressions of the measure of 
\begin{equation*}
\{ (m,e, i, \theta); \, \theta \in  \Theta_m^{(i)e}\} \subset \Omega \times \Z^2 \times \{0,1,2\} \times [0,\pi/2[ ,
\end{equation*}
where $\Theta_m^{(0)e} := \Theta_m^e$, $\Theta_m^{(1)e} := \Theta_m'^e$, and $\Theta_m^{(2)e} := \Theta_m''^e$. Implicitly, we equipped $\Z^2$ and $\{0,1,2\}$ with the counting measure, and $[0,\pi/2[$ and $\Omega$ with the Lebesgue measure (which in the latter case is preserved by the rotations $R_\theta$).
\end{proof}

In order to estimate \eqref{eq:MeasEqual}, we bound in the next lemma the size of the sets $\Theta_m^e$, $\Theta_m'^e$, $\Theta_m''^e$. 

\begin{Lemma}
\label{lem:DoubleQuad}
Let $e\in \Z^2$ be irreducible, with $\|e\|>1$, of parents $f,g$. Let also $m \in \Omega$. Then for any $\theta, \vp \in \Theta_m^e$, one has $\sin |\theta-\vp| \leq 2/ \min\{\<e,f\>, \<e,g\>\}$. Likewise for $\Theta_m'^e$, $\Theta_m''^e$.
\end{Lemma}

\begin{proof}
Without loss of generality, we may assume that $m$ is the origin of $\R^2$. Let $Q$ be the parallelogram of vertices $\{\pm R_\theta e, \pm R_\vp e\}$; note that $\frac 1 2 Q \subset \Omega$. A point $x \in \R^2$ belongs to $Q$ iff
\begin{equation}
\label{eq:InQ}
|\det(x, R_\theta e \pm R_\vp e)| \leq |\det(R_\theta e, R_\vp e)| = \|e\|^2 \sin |\theta- \vp|.
\end{equation}
Indeed $\sin |\theta-\vp| = |\sin(\theta-\vp)|$, since $\theta, \vp \in [0, \pi/2[$ by construction. 
We assume without loss of generality that $\<e,f\> \leq \<e,g\>$. Introducing $h:=e- 2 f = g-f$ we observe that $\<e,h\>  \geq 0$, and compute
\begin{align*}
|\det(R_\theta h, R_\theta e)| &= |\det(h,e)| = |\det(e-2 f,e)| = 2.\\
|\det(R_\theta h, R_\vp e)| & \leq |\det(h,e)| \cos(\vp-\theta) + |\<h,e\>| \sin |\vp-\theta| \leq 2+( \|e\|^2 - 2 \<e,f\>) \sin | \vp-\theta|.
\end{align*}
In the second line, we used the identity $\sin(a+b) = \sin(a) \cos(b) + \cos(a) \sin(b)$, where $a$ denotes the angle between $e$ and $h$, and $b:=\vp-\theta$.
Combining these two estimates with \eqref{eq:InQ}, and assuming for contradiction that $\sin |\theta-\vp| \geq 2 / \<e,f\>$, we obtain that $R_\theta h \in Q$. By symmetry, $-R_\theta h \in Q$, and likewise $\pm R_\vp h \in Q$. 

In the following, we denote $x := -R_\theta e/2$, $y := -R_\vp e/2$, $p := \pm R_\theta h/2$, $q := \pm R_\vp h/2$, where the signs for $p$ and $q$ are chosen so that $p,q \in [x,-x,y]$. 
Denoting by $\alpha, \beta, \gamma$ (resp.\ $\alpha', \beta', \gamma'$) the barycentric coordinates of $p$ (resp.\ $q$) in this triangle, convexity implies
\begin{align}
\label{eq:xqqr}
\uC(p) &\leq \alpha \uC(x) + \beta \uC(-x)+ \gamma \uC (y),\\
\label{eq:yqqr}
\uC(q) &\leq \alpha' \uC(x) + \beta' \uC(-x)+ \gamma' \uC (y).
\end{align}
Let $\xi\in [0,1[^2$ be such that $m \in A_\theta^\xi(\cZ)$. 
Then $x \in X_\theta^\xi$, $e \in \cV_\theta^\xi(x)$, and $m=x+R_\theta e/2$ (recall that we fixed $m=0$). Using the characterization of minimal stencils of Proposition \ref{prop:MinimalStencilsCharacterization}, we obtain
\begin{equation}
\label{eq:PUNeg}
0 > P_{x, \theta}^e(U)  = U(x) - U(x+R_\theta f)-U(x+R_\theta g) + U(x + R_\theta e),
\end{equation}
provided this linear form is supported on $X_\theta^\xi$.
Note that $x+R_\theta e = -x \in X_\theta^\xi$, that $x+R_\theta f = \ve p$, and that $x+R_\theta g = -\ve p$ for some $\ve \in \{-1,1\}$. Since $\pm p \in \frac 1 2 Q \subset \Omega$ this confirms that \eqref{eq:PUNeg} is supported on $X_\theta^\xi := \Omega \cap A_\theta^\xi (\Z^2)$.
Inserting in \eqref{eq:PUNeg} the values of $-x,p,-p$, and proceeding likewise for $y$ and $q$, we obtain
\begin{align}
\label{eq:ppxx}
\uC(x) + \uC(-x) &< U(p)+U(-p)\\
\label{eq:qqyy}
\uC(y) + \uC(-y) &< U(q)+U(-q).
\end{align}
Up to adding an affine map to $\uC$, we may assume that $\uC (x) = \uC(-x) = 0$, and $\uC(p) = \uC(-p)$. From \eqref{eq:ppxx} we obtain $\uC(p) > 0$. Hence also $U(y)>0$ using \eqref{eq:xqqr}, and therefore $\uC(q) \leq \gamma' \uC(y) \leq \uC(y)$ using \eqref{eq:yqqr}. Likewise $\uC(-q) \leq \uC(-y)$, which contradicts \eqref{eq:qqyy} and concludes the proof.
\end{proof}

We finally conclude the proof of Theorem \ref{th:Avg}, by combining \eqref{eq:MeasEqual} with the next lemma.

\begin{Lemma}
For any $m\in \Omega$, with $r := \max\{1,\diam(\Omega)\}$: (likewise for $\Theta_m'^e$, $\Theta_m''^e$)
\begin{equation}
\label{eq:SumOmega}
\sum_{e \in \Z^2} | \Theta_m^e| \leq 2 \pi + 4 \pi^2 (1+\ln r)^2.
\end{equation}
\end{Lemma}

\begin{proof}
Note that $\Theta_m^e \subset [0,\pi/2[$ for any $e\in \Z^2$, and that $\Theta_m^e = \emptyset$ if the $e$ is not irreducible, or if $\|e\| > r$. 
For any two vectors $e,e' \in \Z^2$, we write $e' \lhd e$ iff $e'$ is a parent of $e$ (which implies that $e,e'$ are irreducible, and that $\|e\| >1$).  Isolating the contributions to \eqref{eq:SumOmega} of the four unit vectors, and applying Lemma \ref{lem:DoubleQuad} to other vectors, we thus obtain 
\begin{equation}
\label{eq:SumVariablesExchange}
\sum_{e \in \Z^2} |\Theta_m^e| \leq 4 \times \frac \pi 2 + \sum_{\|e\| \leq r} \sum_{e' \lhd \, e} \arcsin\left(\frac 2 {\<e,e'\>}\right) \leq 2 \pi +  \sum_{e'\in \Z^2} \sum_{\substack{e \, \rhd e'\\ \|e\| \leq r}} \frac \pi {\<e,e'\>},
\end{equation}
where we used the concavity bound $\arcsin(x) \leq \frac \pi 2 x$ for all $x \in [0,1]$ (and slightly abused notations for arguments of $\arcsin$ larger than $1$). 

Consider a fixed irreducible $e' \in \Z^2$, and denote by $f,g$ its parents if $\|e'\| > 1$, or the two orthogonal unit vectors if $\|e'\|=1$, so that $\det(f,e') = 1 = \det(e',g)$. If $e \in \Z^2$ is such that $e' \lhd e$, then $|\det(e,e')|= 1$; assuming $\det(e,e')=1$ (resp.\ $-1$) we obtain that $e = f + k e'$  (resp.\ $e=g+k e'$) for some scalar $k$ which must be (i) an integer since $e'$ is irreducible, (ii)  non-negative since $\<e,e'\> \geq 0$, and (iii) positive since $\|e'\| < \|e\|$. Assuming $\|e\| \leq r$, we obtain in addition $k\|e'\| \leq r$, thus $k\leq r$ and $\|e'\| \leq r$. As a result 
\begin{equation}
\label{eq:SumEp}
\sum_{\substack{e \, \rhd e'\\ \|e\| \leq r}} \frac 1 {\<e,e'\>} \leq \sum_{1 \leq k \leq r} \left(\frac 1 {\<f+k e', e'\>} + \frac 1 {\<g+k e', e'\>}\right) \leq \frac 2 {\|e'\|^2} \sum_{1 \leq k \leq r} \frac 1 k. 
\end{equation}
Inserting \eqref{eq:SumEp} into \eqref{eq:SumVariablesExchange} yields the product of the two following sums, which are easily bounded via comparisons with integrals: isolating the terms for $k=1$, and for all $\|e'\| \leq \sqrt 2$
\begin{equation*}
\sum_{1\leq k \leq r} \frac 1 k \leq 1+\int_1^r \frac {dt} t = 1+\ln r, \qquad 
\sum_{\substack{0 < \|e'\| \leq r\\ \text{irreducible}}} \frac 1 {\|e'\|^2} \leq 4+4 \times \frac 1 2 + \int_{1 \leq \|x\| \leq r} \frac {dx} {\|x\|^2} \leq 6 + 2 \pi \ln r.
\end{equation*}
Noticing that $2 \pi \geq 6$ we obtain \eqref{eq:SumOmega} as announced.
\end{proof}


\section{Numerical experiments}
\label{sec:Numerics}

Our numerical experiments cover the classical formulation \cite{Rochet:1998uj} of the monopolist problem, as well as several variants, including lotteries \cite{Manelli:2006ib,Thanassoulis:2004uy}, or the pricing of risky assets \cite{Carlier:2007gy}. We choose to emphasize this application in view of its appealing economical interpretation, and the often surprising qualitative behavior.
Our algorithm can also be applied in a straightforward manner to the computation of projections onto the cone of convex functions defined on some square domain, with respect to various norms as considered in \cite{Carlier:2001tq,MERIGOT:tr,Oberman:2011wi,Oberman:2011wy} (this amounts to denoising under a convexity prior). It may however not be perfectly adequate for the investigation of geometric conjectures \cite{LachandRobert:2005bi,Wachsmuth:2013ta,MERIGOT:tr}, due to the use of a grid discretization.

The hierarchical cones of discrete convex functions introduced in this paper are combined with a simple yet adaptive and anisotropic stencil refinement strategy, described in \S\ref{sec:strategy}. The monopolist model is introduced in \S \ref{sec:Monopolist}, and illustrated with numerous experiments. We compare in \S \ref{sec:Comparison} our implementation of the constraint of convexity, with alternative methods proposed in the literature, in terms of computation time, memory usage, and solution quality.

\subsection{Stencil refinement strategy}
\label{sec:strategy}

We introduce two algorithms which purpose is to minimize a given lower semi-continuous proper convex functional $\cE : \cF(\OmegaD) \to \R \cup \{+\infty\}$, on the $N$-dimensional cone $\Conv(\OmegaD)$, $N := \#(X)$, without ever listing the $\cO(N^2)$ linear constraints which characterize this cone. 
They both generate an increasing sequence of stencils $\cV_0 \subsetneq \cV_1 \subsetneq \cdots \subsetneq \cV_n$ on $\OmegaD$, and minimizers $(u_i)_{i=0}^n$ of $\cE$ on cones defined by $\cO(\#(\cV_i))$ linear constraints. The subscript $i$ refers to the loop iteration count in Algorithms \ref{algo:SubCones} and \ref{algo:SuperCones}, and the loop ends when the stencils are detected to stabilize: $\cV_n = \cV_{n+1}$. The final map $u_n$ is guaranteed to be the global minimum of $\cE$ on $\Conv(\OmegaD)$. 

Our first algorithm is based on an increasing sequence $\Conv(\cV_0) \subset \cdots \subset \Conv(\cV_n)$ of sub-cones of $\Conv(\OmegaD)$. 
If constraints of type $P_x^e$, $x \in X$, $e \in \cV_i(x)$ are \emph{active} for the minimizer $u_i$ of $\cE$ on $\Conv(\cV_i)$ (i.e.\ the corresponding Lagrange multipliers are positive), then refined stencils $\cV_{i+1}$ are adaptively generated from $\cV_i$; otherwise $u_i$ is the global minimizer of $\cE$ on $\Conv(X)$, and the method ends.
Note that the optimization of $\cE$ on $\Conv(\cV_{i+1})$ can be hot-started from the previous minimizer $u_i \in \Conv(\cV_i) \subset \Conv(\cV_{i+1})$.

\begin{algorithm}
\caption{Sub-cones approach to stencil refinement}
\label{algo:SubCones}
Start with the minimal stencils: $\cV \leftarrow \cV_{\min}$. (See Proposition \ref{prop:MinimalStencils})\\
\textbf{Until} the stencils $\cV$ stabilize\\
\phantom{ab} \textbf{Find} a minimizer $u$ of the energy $\cE$ on $\Conv(\cV)$, \\
\phantom{ab Fi} and extract the Lagrange multipliers $\lambda$ associated to the constraints $P_x^e$, $x\in X$, $e \in \hat \cV(x)$.\\
\phantom{ab} \textbf{Set} $\cV(x) \leftarrow \cV(x) \cup \{e \in \hat \cV(x);\, \lambda(P_x^e) > 0\}$, for all $x\in \OmegaD$.
\end{algorithm}

\begin{Definition}
For any family $\cV$ of stencils on $X$, we denote by $\Conv'(\cV) \subset \cF(X)$ the cone 
defined by the non-negativity of: for all $x \in X$, and all $e\in \cV(x)$, the linear forms $S_x^e$ and $T_x^e$ if $\|e\|>1$, provided they are supported on $\OmegaD$. Note that $\Conv(\cV) \subset \Conv(X) \subset \Conv'(\cV)$.

\end{Definition}

Algorithm \ref{algo:SuperCones} is based on a \emph{decreasing} sequence $\Conv'(\cV_0) \supset \cdots \supset \Conv'(\cV_n)$ of super-cones of $\Conv(\OmegaD)$. 
The minimizer $u_i$ of $\cE$ on $\Conv'(\cV_i)$ may not belong to $\Conv(\OmegaD)$, even less to $\Conv(\cV_i)$, and in particular some of the values $P_x^e(u_i)$, $x\in \OmegaD$, $e \in \hat \cV_i(x)$, may be negative. In that case,
 refined stencils $\cV_{i+1}$ are adaptively generated from $\cV_i$; otherwise, $u_i$ is the global minimizer of $\cE$ on $\Conv(\OmegaD)$, and the method ends. 

\begin{figure}
\includegraphics[width=3cm]{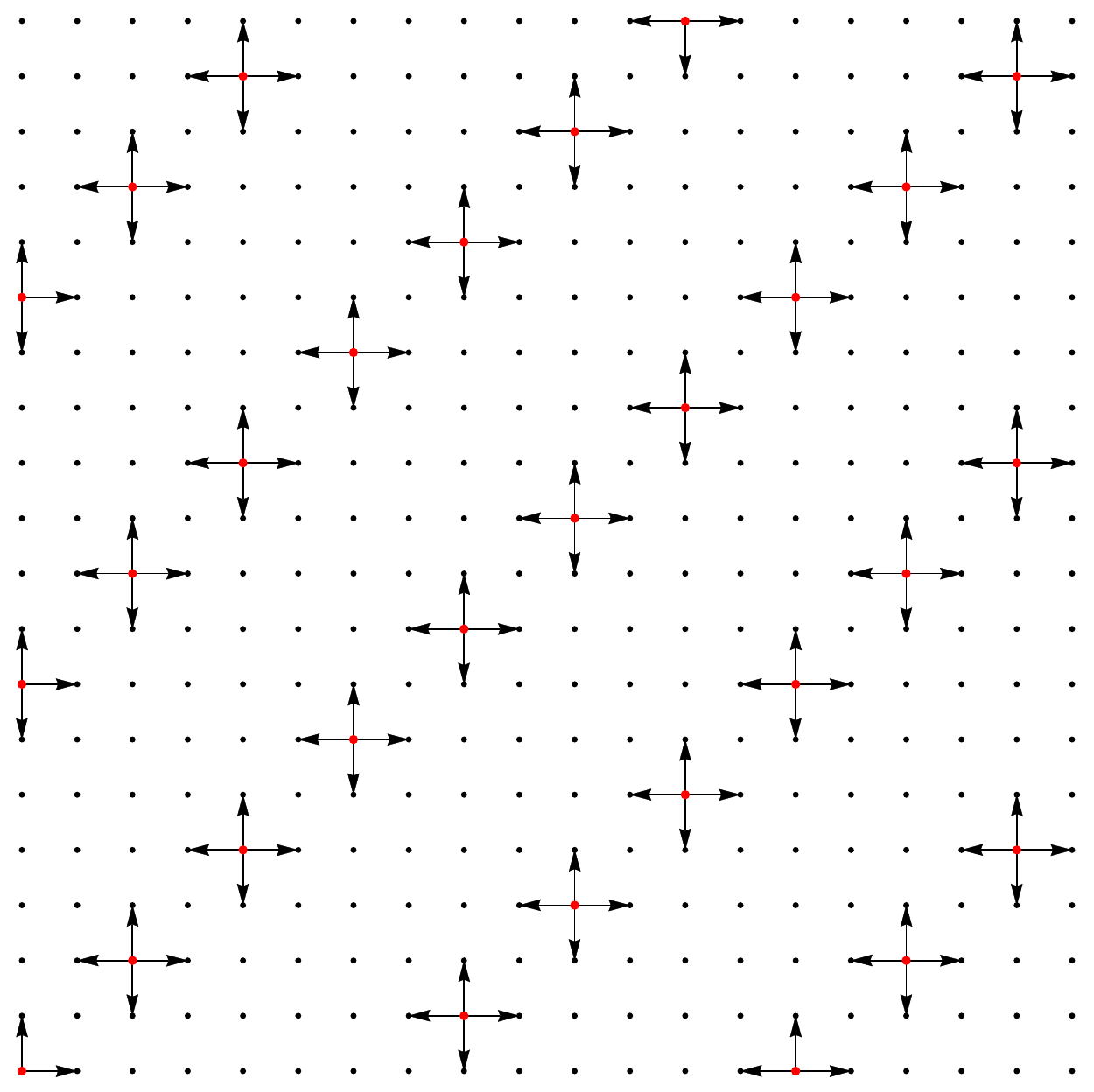}
\includegraphics[width=3cm]{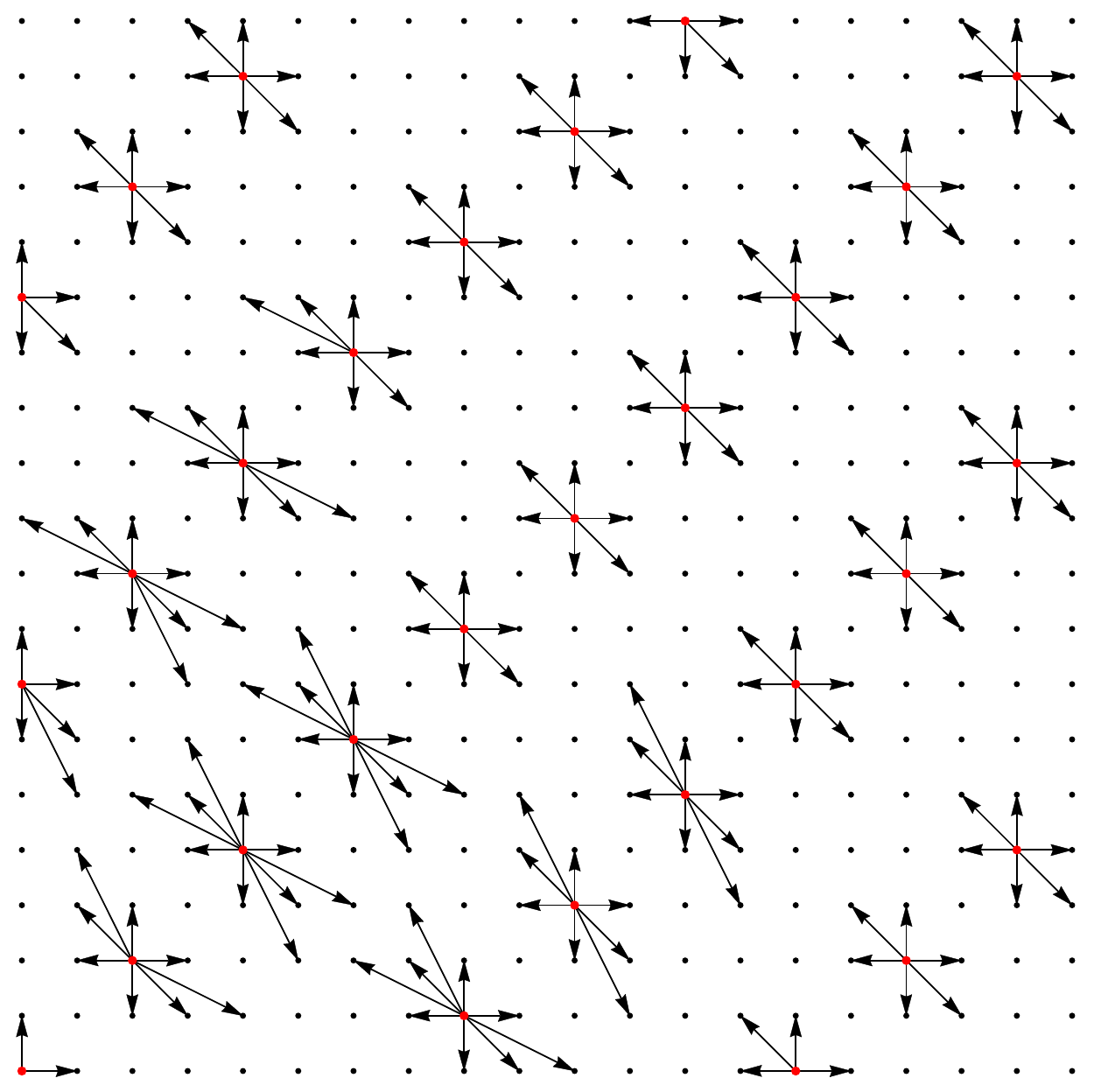}
\includegraphics[width=3cm]{\pathPic/Convex/PrincipalAgent/20/Stencils/Default_2.pdf}
\includegraphics[width=3cm]{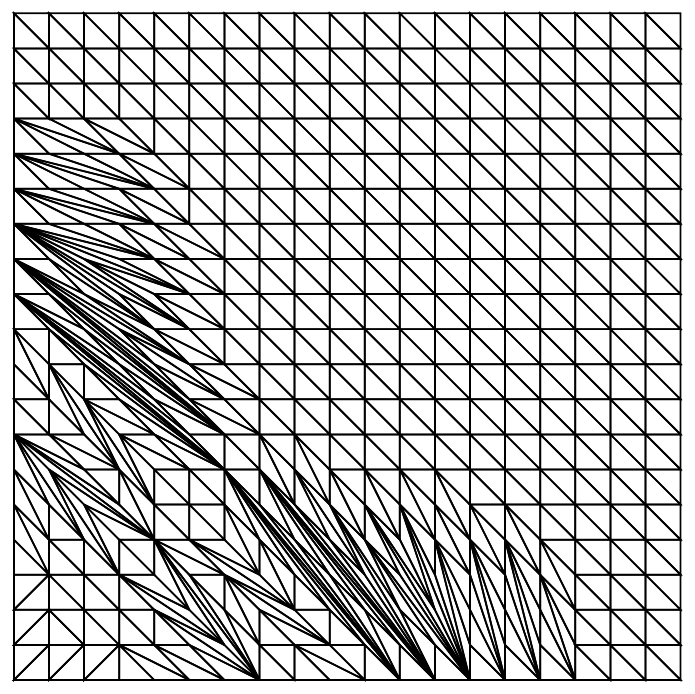}
\includegraphics[width=3.3cm]{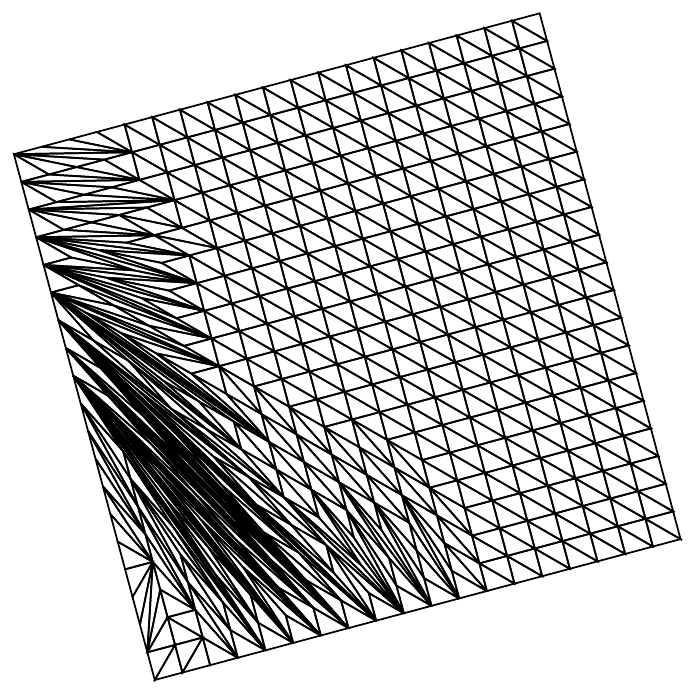}\\
\includegraphics[width=3.3cm]{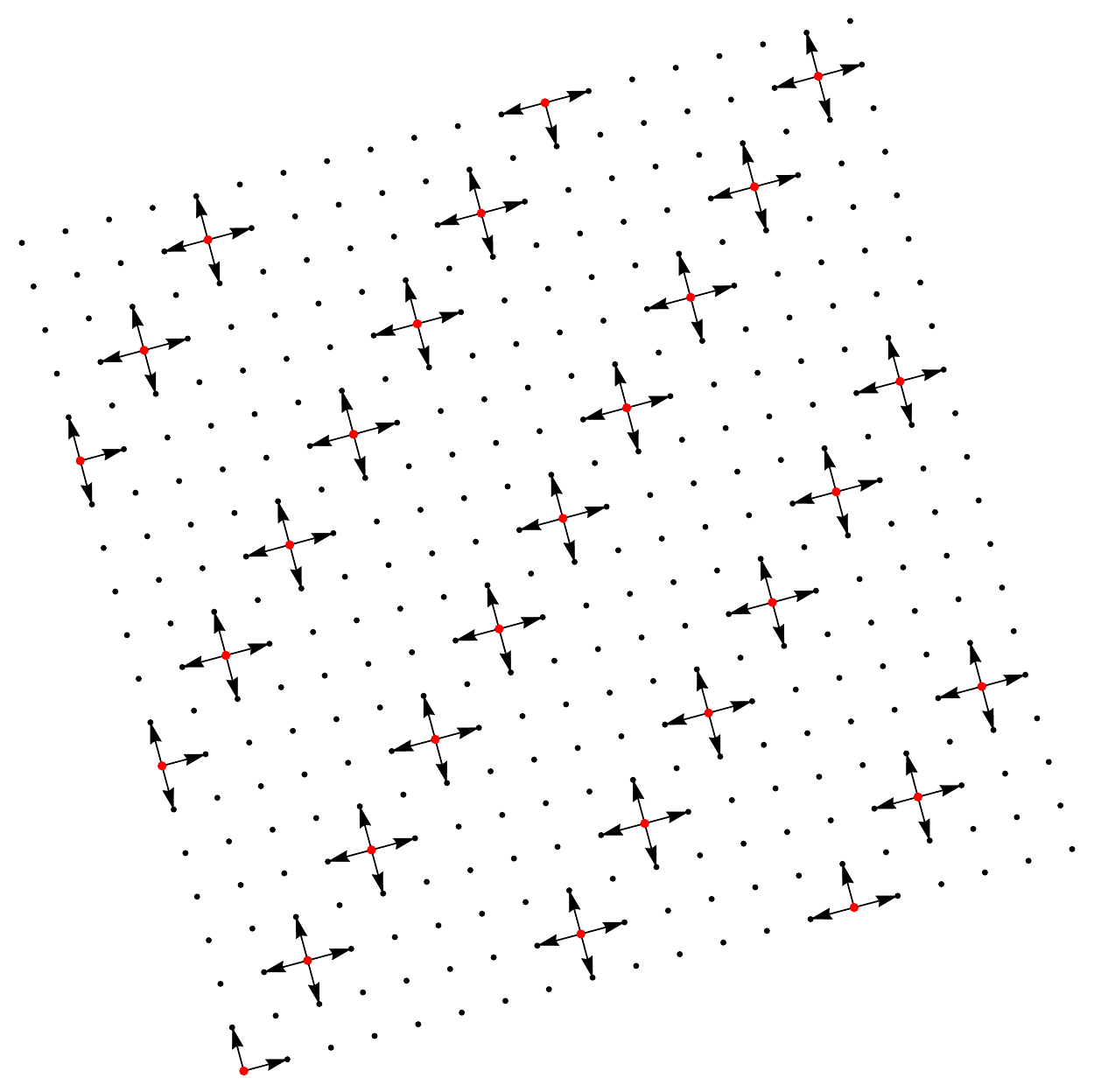}
\hspace{-0.6cm}
\includegraphics[width=3.3cm]{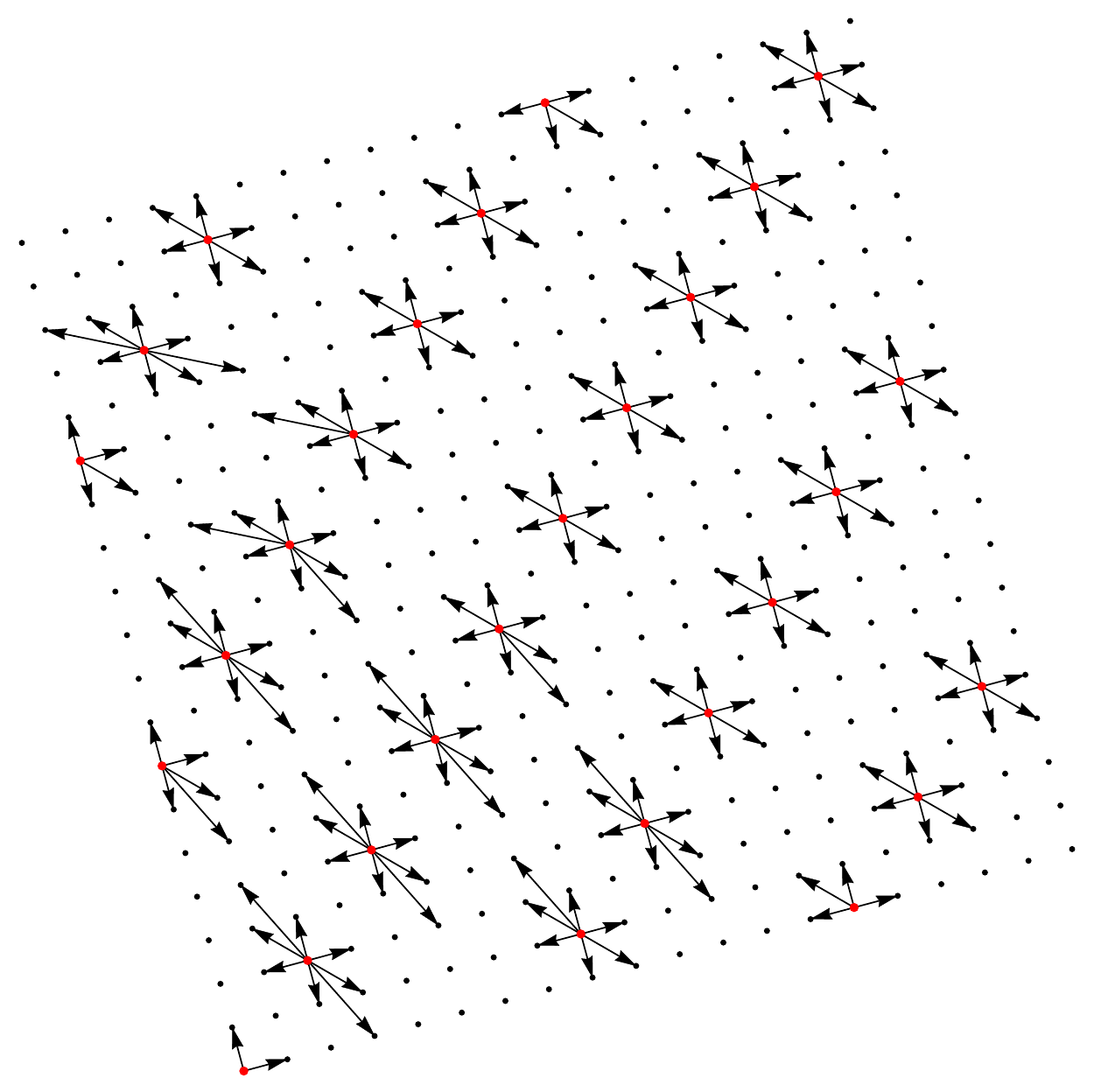}
\hspace{-0.6cm}
\includegraphics[width=3.3cm]{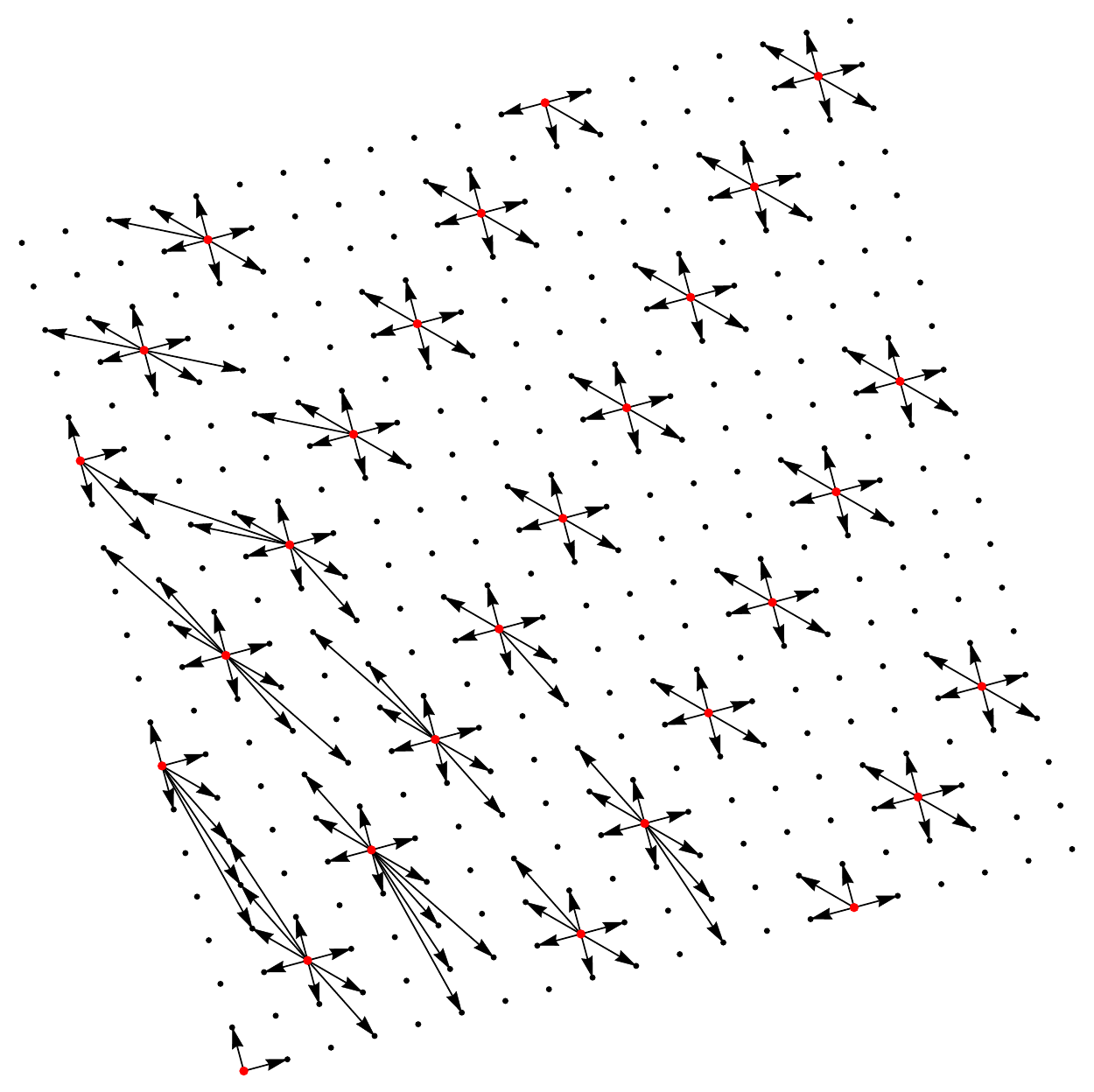}
\hspace{-0.6cm}
\includegraphics[width=3.3cm]{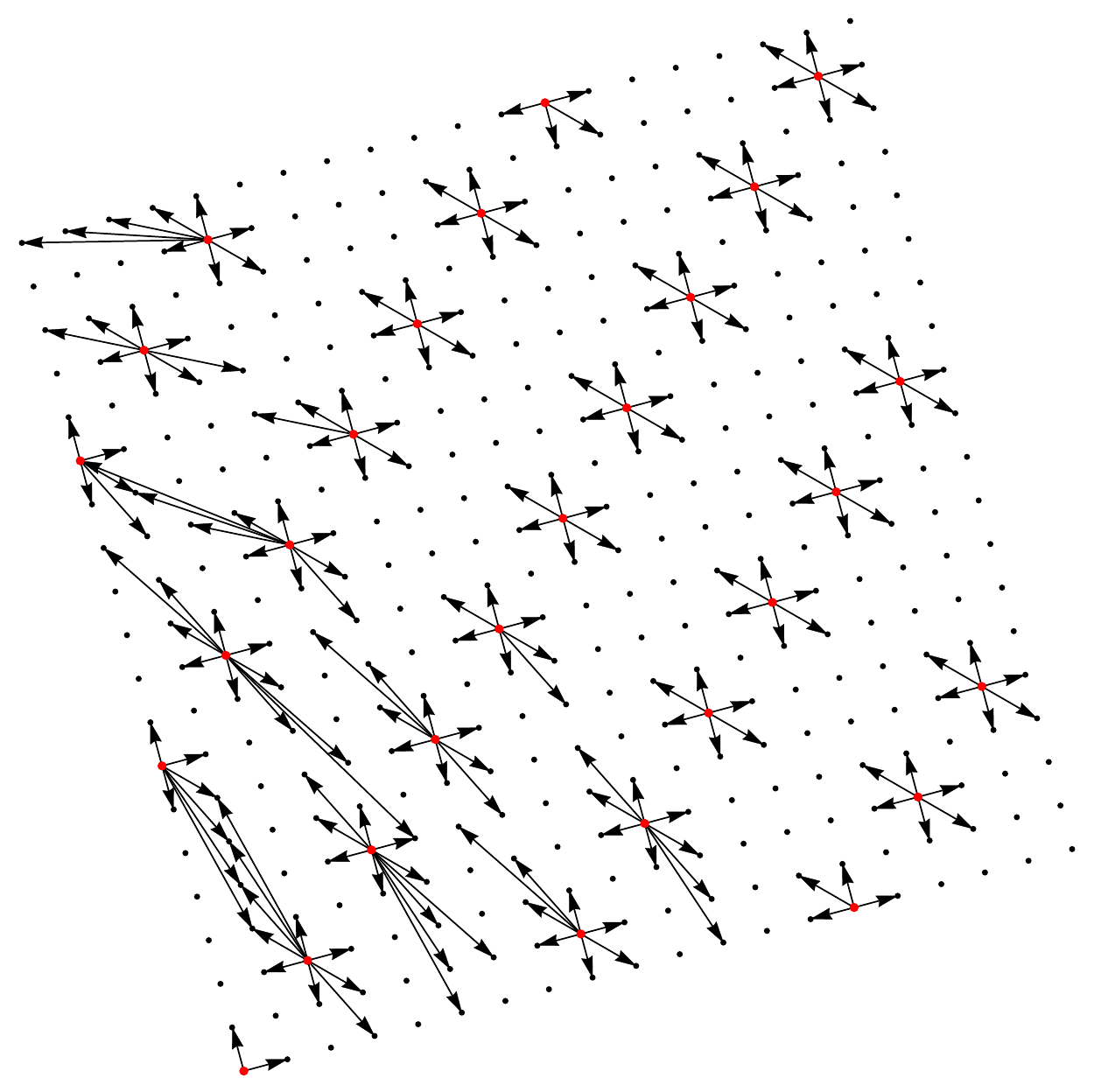}
\hspace{-0.6cm}
\includegraphics[width=3.3cm]{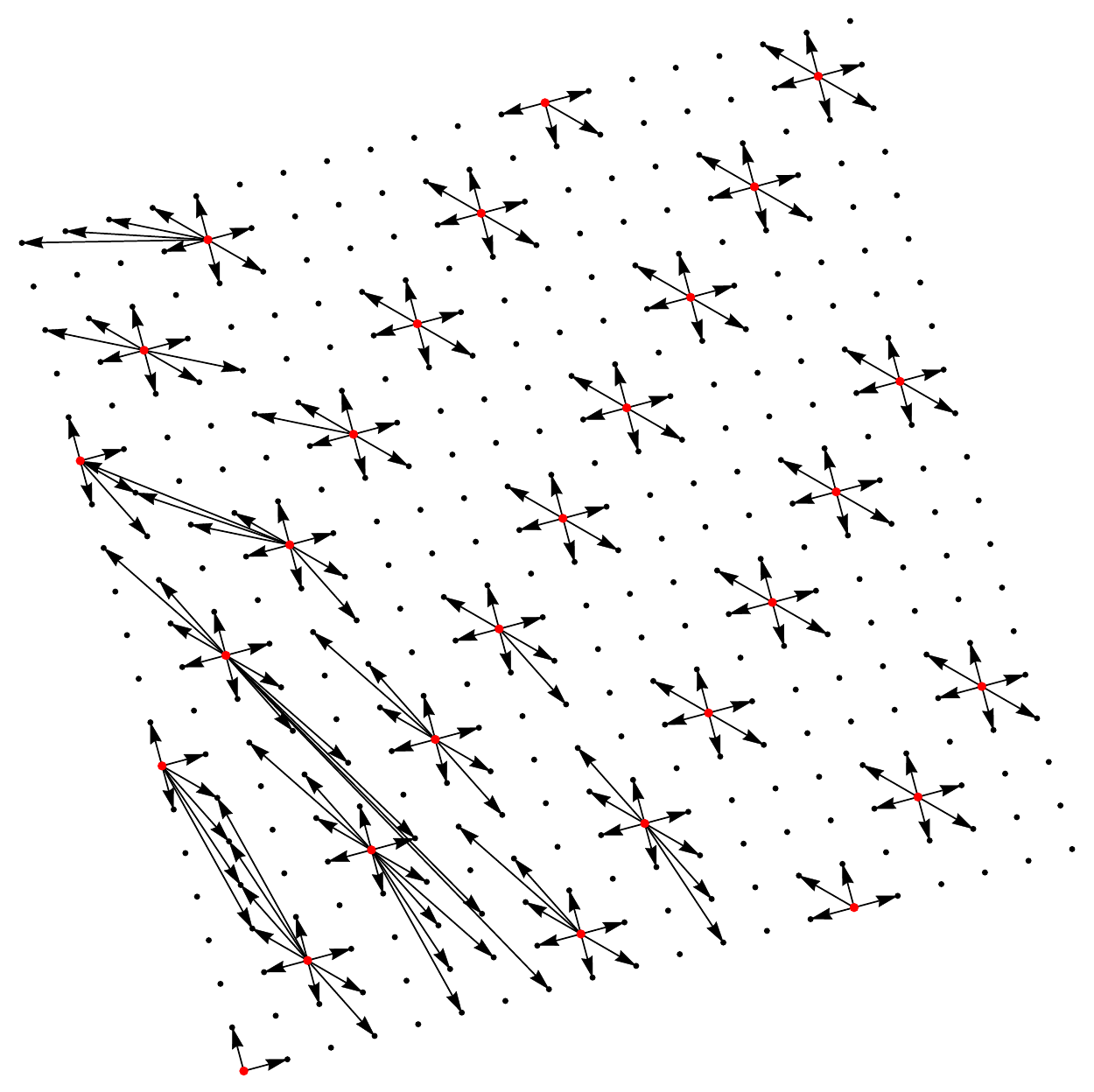}
\caption{
Top: Algorithm 2, for the classical monopolist problem \eqref{eq:PrincipalAgent} on $[1,2]^2$ with a $20 \times 20$ grid, converges in $2$ stencils refinement steps (using the extended candidates $\hat \cV_\rho$, $\rho:=1.5$). Bottom: $4$ refinement steps are needed with a different density of customers, uniform on the square $[1,2]^2$ rotated by $\pi/12$. Top right: $u$-Delaunay triangulations associated with the respective discrete solutions, for illustration of Theorem \ref{th:Decomp}.
}
\label{fig:StencilRefinement}
\end{figure}

\begin{algorithm}
\caption{Super-cones approach to stencil refinement}
\label{algo:SuperCones}
Start with the minimal stencils $\cV$.\\
\textbf{Until} the stencils $\cV$ stabilize\\
\phantom{ab} \textbf{Find} a minimizer $u$ of the energy $\cE$ on $\Conv'(\cV)$.\\
\phantom{ab} \textbf{Set} $\cV(z) \leftarrow \cV(z) \cup \{e \in \hat \cV(z); \, P_z^e(u) < 0\}$, for all $z\in \OmegaD$.
\end{algorithm}

Algorithms \ref{algo:SubCones} and \ref{algo:SuperCones} are provided ``as is'', without any complexity guarantee. 
Our numerical experiments are based on algorithm \ref{algo:SuperCones}, because the numerical test ``$P_x^e(u) < 0$'' turned out to be more robust than ``$\lambda(P_x^e)>0$''. 
In order to limit the number $n$ of stencil refinement steps, we use a slightly extended set $\hat \cV_\rho(x)$ of candidates for refinement, with $\rho := 1.5$, see Definition \ref{def:ExtendedCandidates} below (note that $\hat \cV_1 = \hat \cV$).
With this modification, $n$ remained below 10 in all our experiments.
The constructed stencils were generally sparse, highly anisotropic, and almost minimal for the discrete problem solution $u\in \Conv(X)$ eventually found, see Figure \ref{fig:StencilRefinement}. 
Observation of Figure \ref{fig:Comparison2} suggests that $n$ grows logarithmically with the problem dimension $N$, and that the final stencils cardinality $\#(\cV_n)$ depends quasi-linearly on $N$, as could be expected in view of Remark \ref{rem:OptStrategy} and Theorem \ref{th:Avg}. 
However, we could not mathematically establish such complexity estimates.

\begin{Definition}[Extended candidates]
\label{def:ExtendedCandidates}
Let $\cV$ be a family of stencils, let $\rho \geq 1$, and let $x \in X$. A vector $e \in \cV_{\max}(x) \sm \cV(x)$, of parents $f,g$, belongs to the extended candidates $\hat \cV_\rho(x)$ iff there exists trigonometrically consecutive $f',g' \in \cV(x)$ such that $f' \preceq f \prec g \preceq g'$ and $\|f\| \|g\| \leq \rho \|f'\| \|g'\|$.
\end{Definition}


\subsection{The monopolist problem} 
\label{sec:Monopolist}

A monopolist has the ability to produce goods, which have two characteristics and may thus be represented by a point $q \in \R^2$. The manufacturing cost $\Cost(q) : \R^2 \to \R\cup \{+\infty\}$ is known and fixed a-priori. Infinite costs account for products which are ``meaningless'', or impossible to build.
The selling price $\pi : \R^2 \to \R\cup \{+\infty\}$ is fixed unilaterally by the monopolist except for the ``null'' product $(0,0)$, which must be available for free ($\pi(0,0) \leq 0$).
The characteristics of the consumers are also represented by a point $z\in \R^2$, and the utility of product $q$ to consumer $z$ is modeled by the scalar product between their characteristics
\begin{equation}
\label{eq:Utility}
\cU(q,z) := \<q,z\>.
\end{equation}
More general utility pairings $\cU$ are considered in \cite{Figalli:2011tz}, yet the numerical implementation of the resulting optimization problems remains out of reach, see \cite{MERIGOT:tr} for a discussion.
All consumers $z$ are rational, ``screen'' the proposed price catalog $\pi$, and choose the product of maximal \emph{net} utility $\<q,z\> - \pi(q)$ (i.e.\ raw utility minus price). 
Introducing the Legendre-Fenchel dual $U$ of the prices $\pi$: for all $z\in \R^2$
\begin{equation*}
U(z): = \pi^*(z) := \sup_{q\in \R^2} \, \<q,z\> - \pi(z),
\end{equation*}
we observe that the optimal product%
\footnote{%
Strictly speaking, the optimal product $Q(z)$ is an element of the subgradient $\partial_z U$, which (Lebesgue-)almost surely is a singleton $\{\nabla U(z)\}$. Hence we may write \eqref{eq:TotalProfit} in terms of $\nabla U(z)$, provided the density $\mu$ of customers is absolutely continuous with respect to the Lebesgue measure.
} 
for consumer $z$ is $\nabla U(z)$, which is sold at the price
\begin{equation}
\label{eq:PriceU}
\pi(\nabla U(z)) = \<\nabla U(z), z\> - U(z).
\end{equation}
The net utility function $U$ is convex and non-negative by construction, and uniquely determines the products bought and their prices. 
Conversely, any convex non-negative $U$ defines prices $\pi = U^*$ satisfying the admissibility condition $\pi(0) \leq 0$, and such that $\pi^* = U^{**} = U$.
The distribution of the characteristics of the potential customers is known to the monopolist, under the form of a bounded measure $\mu$ on $\R^2$. He aims to maximize his total profit: the integrated difference (sales margin) between the selling price \eqref{eq:PriceU}, and the production cost 
\begin{equation}
\label{eq:TotalProfit}
\sup\left\{ \int_{\R^2} \Big(\<\nabla U(z),z\> - U(z) - \Cost(\nabla U(z))\Big) d\mu(z); \, \, U \in \Conv(\R^2), \, U \geq 0\right\}
\end{equation}
If production costs are convex, then this amounts to maximizing a concave functional of $U$ under convex constraints; see \cite{Carlier:2001gv} for precise existence results.
If $U$ maximizes \eqref{eq:TotalProfit}, then an optimal catalog of prices is given by $U^*$. Quantities of particular economic interest are the monopolist margin, and the distribution of product sales: 
\begin{equation}
\label{eq:ProductMarginDistribution}
{\rm Margin} = U^* - \Cost, \qquad 
{\rm SalesDistribution} = (\nabla U)_{\#} \mu,
\end{equation}
where $\#$ denotes the push forward operator on measures.
The regions defined by $\{U=0\}$ and $\{\det(\Hessian u)=0\}$ are also important, as they correspond to different categories of customers, see below.
We present numerical results for three instances of the monopolist problem, associated to different product costs. These three models are clearly simplistic idealizations of real economy. Their interest lies in their, striking, qualitative properties, which are stable and are expected to transfer to more complex models. 

For implementation purposes, we observe that the maximum profit \eqref{eq:TotalProfit} is unchanged if one considers $U$ only defined on a convex set $K \supseteq \supp (\mu)$, and imposes the additional constraint $\Cost(\nabla U(x)) < \infty$ for all $x \in K$.
The chosen discrete domain is a square grid $X$, such that $\supp(\mu) \subset \Hull(X)$. This density is represented by non-negative weights $(\mu_x)_{x\in X}$, set to zero outside $\supp(\mu)$. 
The integral appearing in \eqref{eq:TotalProfit} is discretized using finite differences, see  \cite{Carlier:2001tq} for convergence results. 
The resulting convex program is solved by combining Mosek software's interior point (for linear problems) or conic (for quadratic%
\footnote{%
Following the indications of Mosek's user manual, quadratic functionals are implemented under the form of linear functionals involving auxiliary variables subject to conic constraints.
} problems)
optimizer, with the stencil refinement strategy of Algorithm \ref{algo:SuperCones}, \S \ref{sec:strategy}.

\paragraph{Classical model.}
The produced goods are cars (for concreteness), which characteristics $q=(q_1,q_2)$ are non-negative and account for the engine horsepower $q_1$ and the upholstery quality $q_2$. Production cost is quadratic: $\Cost(q) := \frac 1 2 \|q\|^2$ for all $q \in \R_+^2$, and $\Cost(q)=+\infty$ otherwise (cars with negative characteristics are unfeasible). Consumer characteristics $x=(x_1,x_2)$ are their appetite $x_1$ for car performance, and $x_2$ for comfort, consistently with \eqref{eq:Utility}. 
The qualitative properties of this model are the following  \cite{Rochet:1998uj}: denoting by $U$ a solution of \eqref{eq:TotalProfit}, and ignoring regularity issues in this heuristic discussion
\begin{itemize}
\item (Desirability of exclusion) The optimal monopolist strategy often involves neglecting a positive proportion of potential customers - which ``buy'' the null product $0$ at price $0$. In other words, the solution of \eqref{eq:TotalProfit} satisfies $U=0$ on an open subset of $\supp(\mu)$, hence also $\nabla U=0$.
The economical interpretation is that introducing (low end) products attractive to this population would reduce overall profit, because other customers currently buying expensive, high margin products, would change their minds and buy these instead.
\item (Bunching)
``Wealthy'' customers $z$ generally buy products which are specifically designed for them, in the sense that $\nabla U$ is a local diffeomorphism close to $z$. ``Poor'' potential customers are excluded from the trade: $U=0$ close to $z$, see the previous point. There also exists an intermediate category of customers characterized by $\det(\Hessian U)=0$ close to $z$, so that the same product $q=\nabla U(z)$ is bought by a one dimensional ``bunch'' of customers $(\nabla U)^{-1}\{q\}$. 
The image of this category of customers, by $\nabla U$, is one dimensional product line. 
From an economic point of view, the optimal strategy limits the variety of intermediate range products in order, again, to avoid competing with high margin sales.
\end{itemize}

Considering, as in \cite{Rochet:1998uj,MERIGOT:tr}, a uniform density of customers on the square $[1,2]^2$, we illustrate%
\footnote{%
With this customer density, \cite{Rochet:1998uj} expected the bunching region to be triangular, and the image $\nabla U$ to be the union of the segment $[(0,0),(1,1)]$ and of the square $[1,2]^2$. After discussion with the author, and in view of the numerical experiments, we believe that these predictions are erroneous.
} on Figure  \ref{fig:Monopolist} the estimated solution $U$ (left), $\det(\Hessian U)$ (center left), the sales distribution (center right) and the monopolist margin (right), see \eqref{eq:ProductMarginDistribution}. 
The phenomena of exclusion $U=0$ and of bunching $\det(\Hessian U)=0$ are visible (center left subfigure) as a white triangle and as the darkest level set of $\det(\Hessian U)$ respectively. 
The image by $\nabla U$ of customers subject to bunching appears (center right subfigure) as a one dimensional red structure in the product sales distribution.

\begin{figure}
\includegraphics[width=3.8cm]{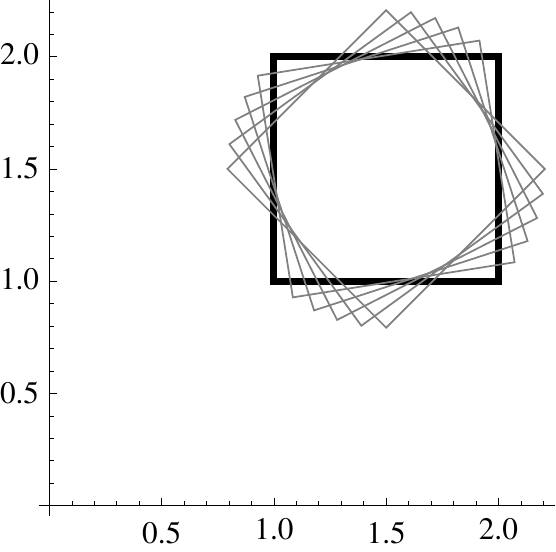}
\includegraphics[width=3.8cm]{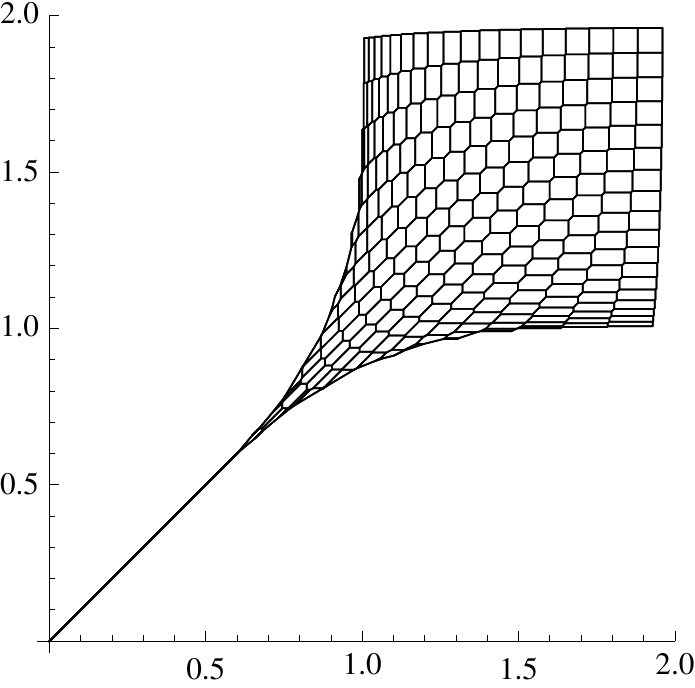}
\includegraphics[width=3.8cm]{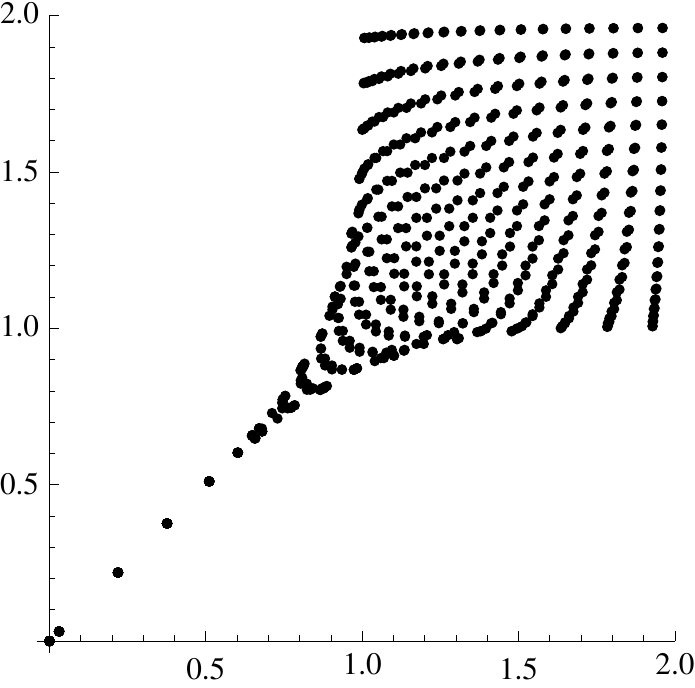}
\includegraphics[width=3.8cm]{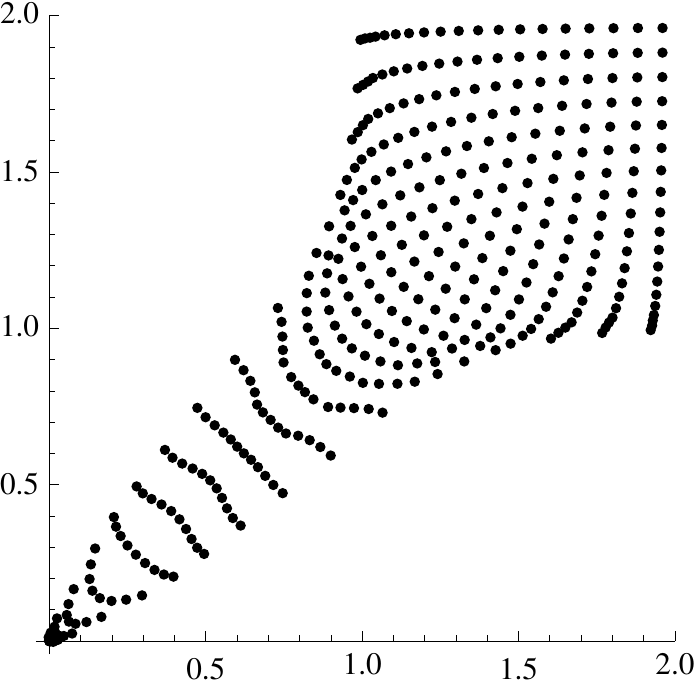}
\caption{
Left: domain $[1,2]^2$ (thick black), and rotated domains used in numerical experiments for the classical principal agent problem, see Figure \ref{fig:PARotated}. 
We computed a minimizer $u \in \Conv(X)$ of a discretization of the classical monopolist problem on $[1,2]^2$, see \eqref{eq:PrincipalAgent}, on a $20 \times 20$ grid $X$, and an $u$-Delaunay triangulation $\cT$, see Figure \ref{fig:Flipping}. 
Center left: subgradients cells $\partial_x U$, $x \in X$, with $U := \interp_\cT u$. Center right: the gradients $\nabla \interp_T u$, $T \in \cT$ (vertices of the previous cells). Right: the less precise numerical method OF$_3$, see \ref{sec:Comparison}, thickens the product line and hides the bunching phenomenon.
}
\label{fig:DTGC}
\end{figure}

\begin{figure}
\phantom{blabla} $\phantom{bla bla } \theta = \pi/8$ \phantom{bla bla bla } $\theta=13\pi/80$ \phantom{bla bla bla } $\theta=\pi/5$ \phantom{bla bla bla b} $\theta=\pi/4$\\
\includegraphics[width=3.6cm]{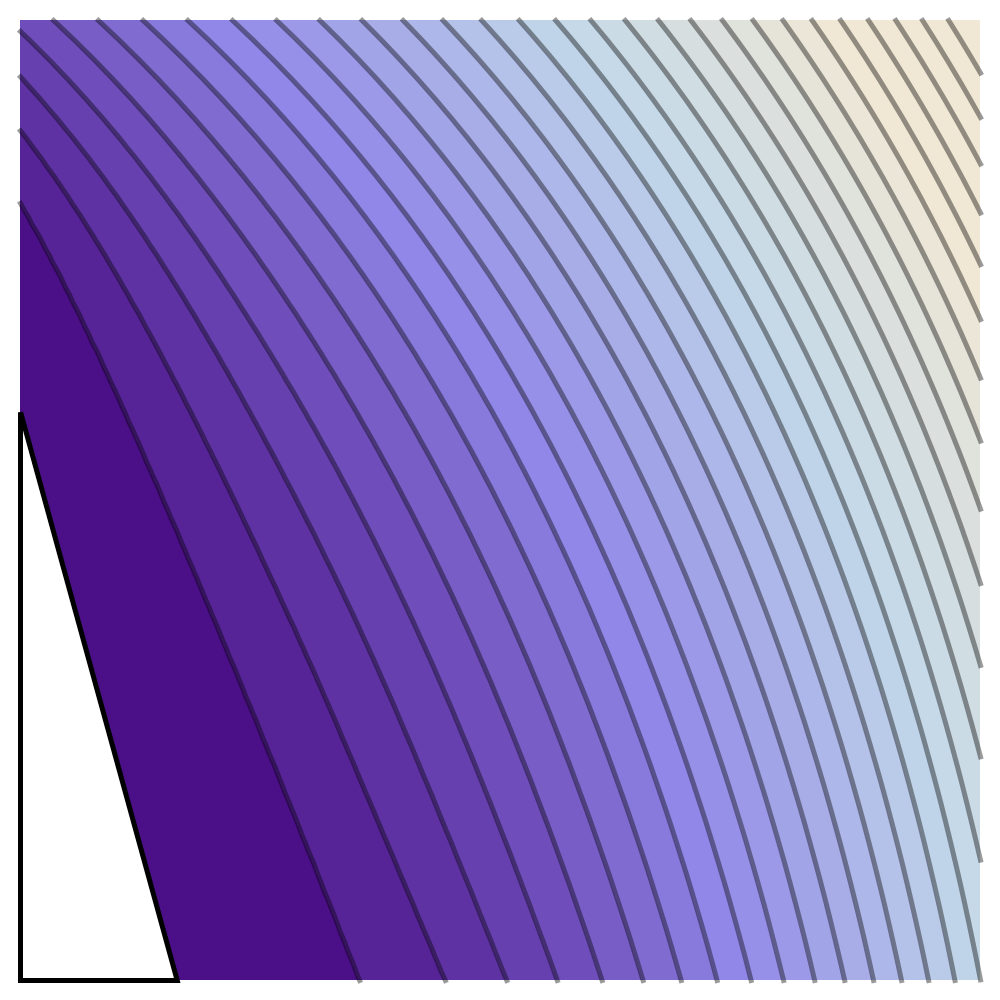}
\includegraphics[width=3.6cm]{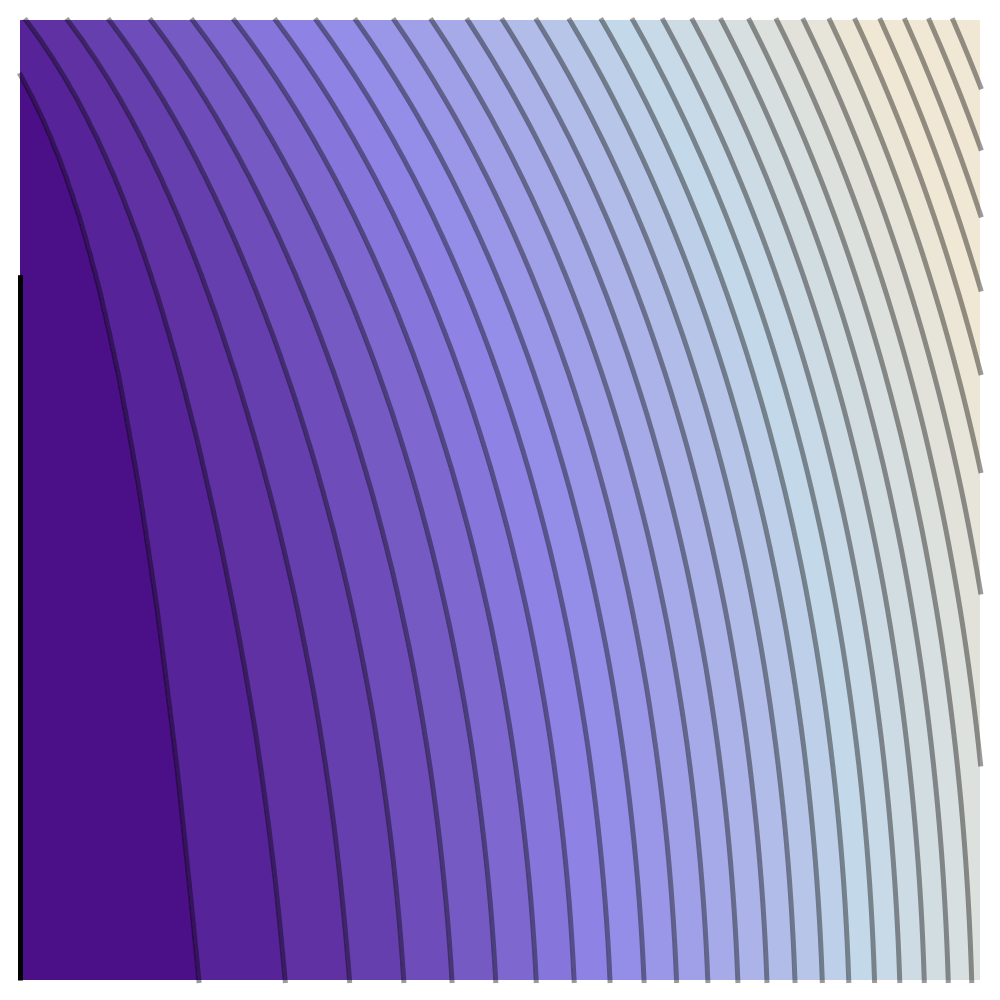}
\includegraphics[width=3.6cm]{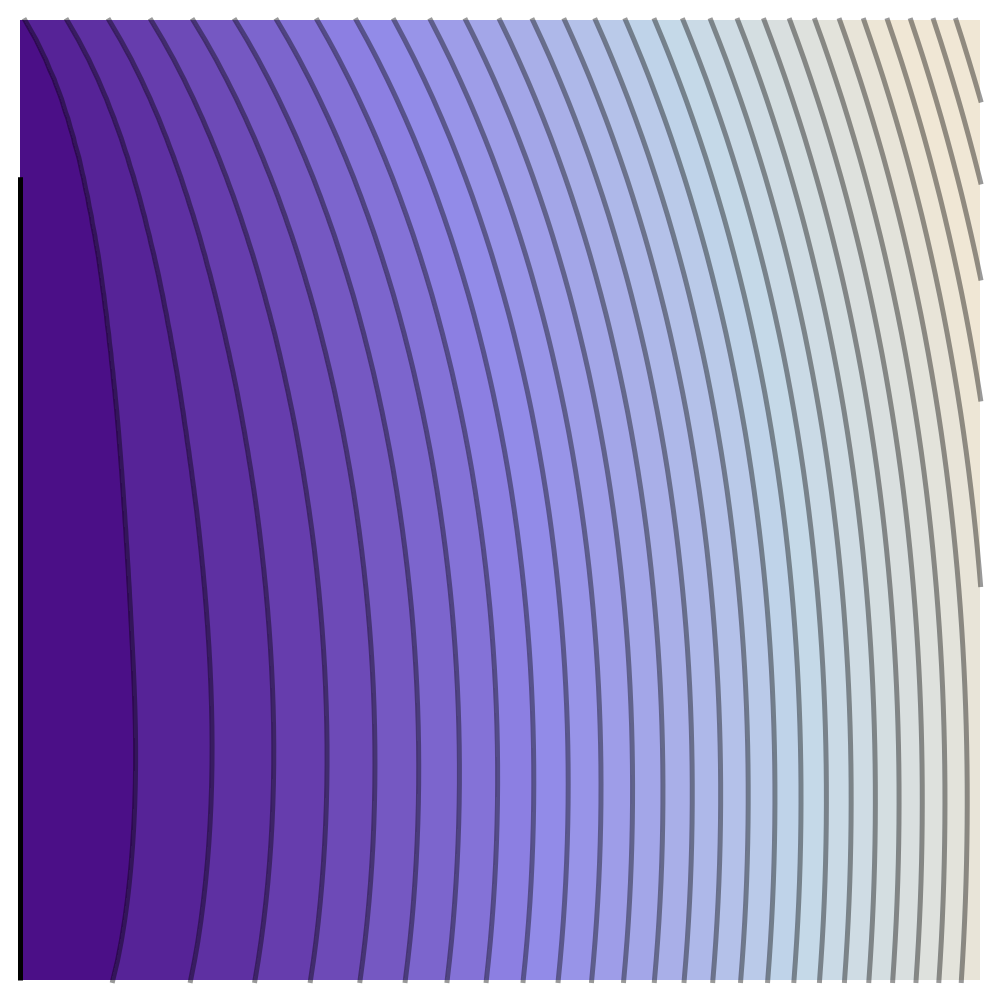}
\includegraphics[width=3.6cm]{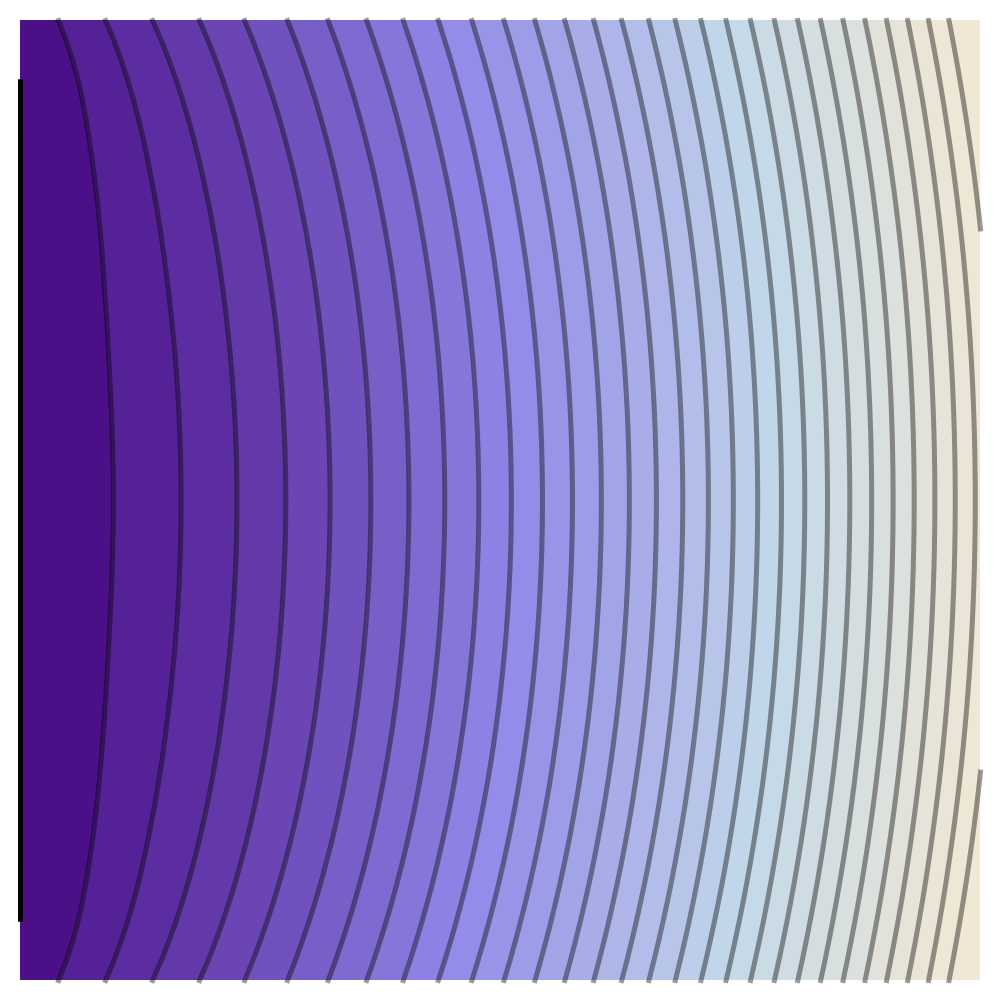}
\includegraphics[width=0.6cm]{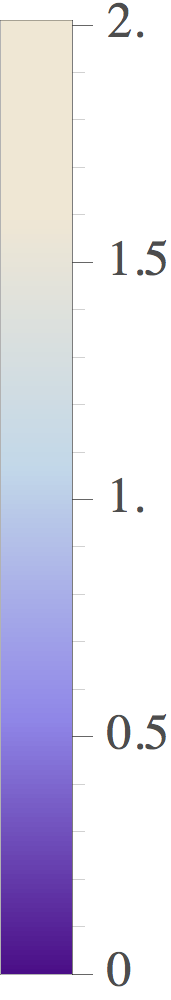}
\\

\includegraphics[width=3.6cm]{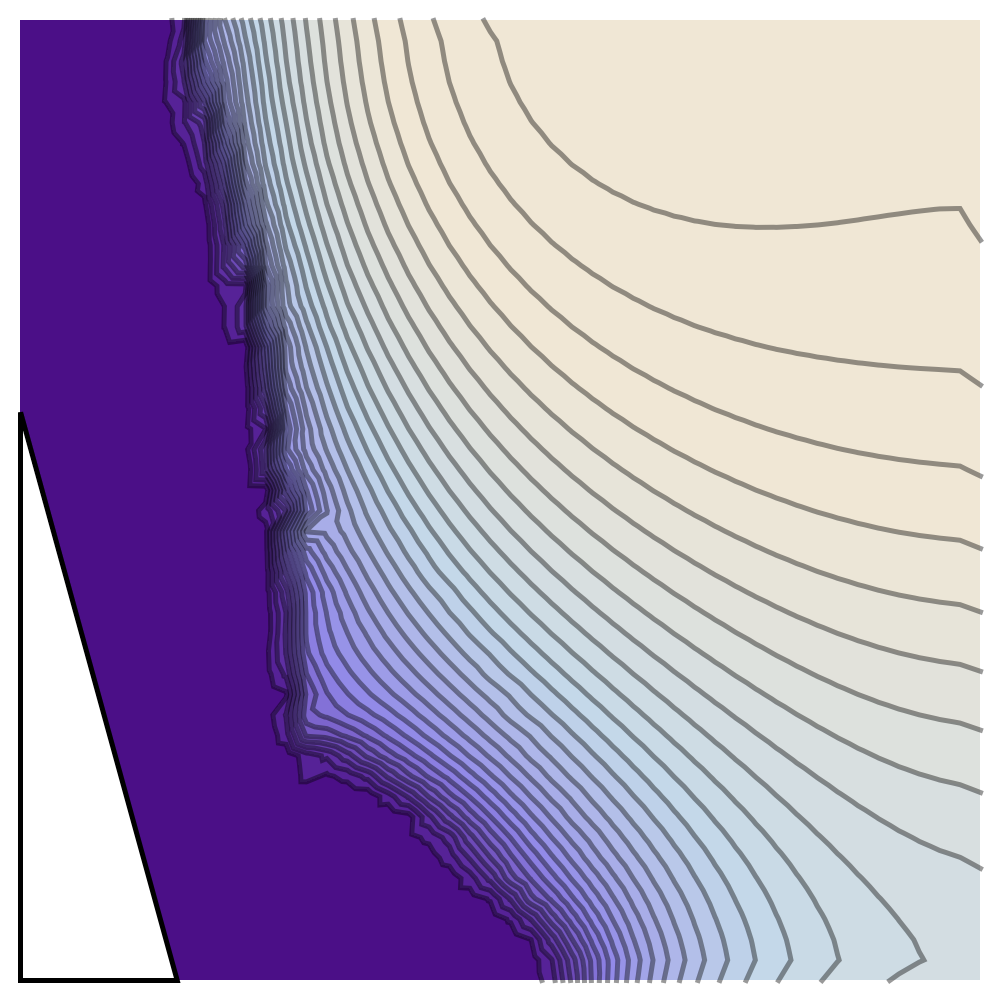}
\includegraphics[width=3.6cm]{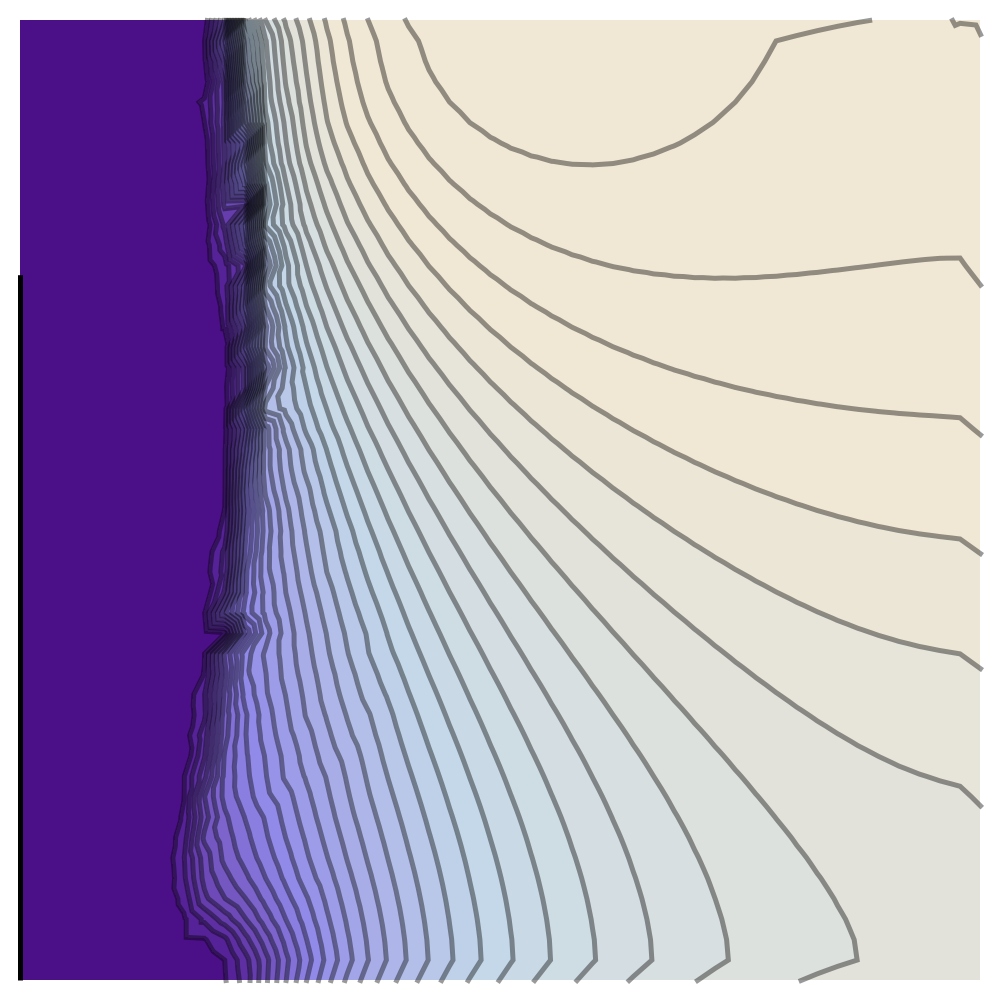}
\includegraphics[width=3.6cm]{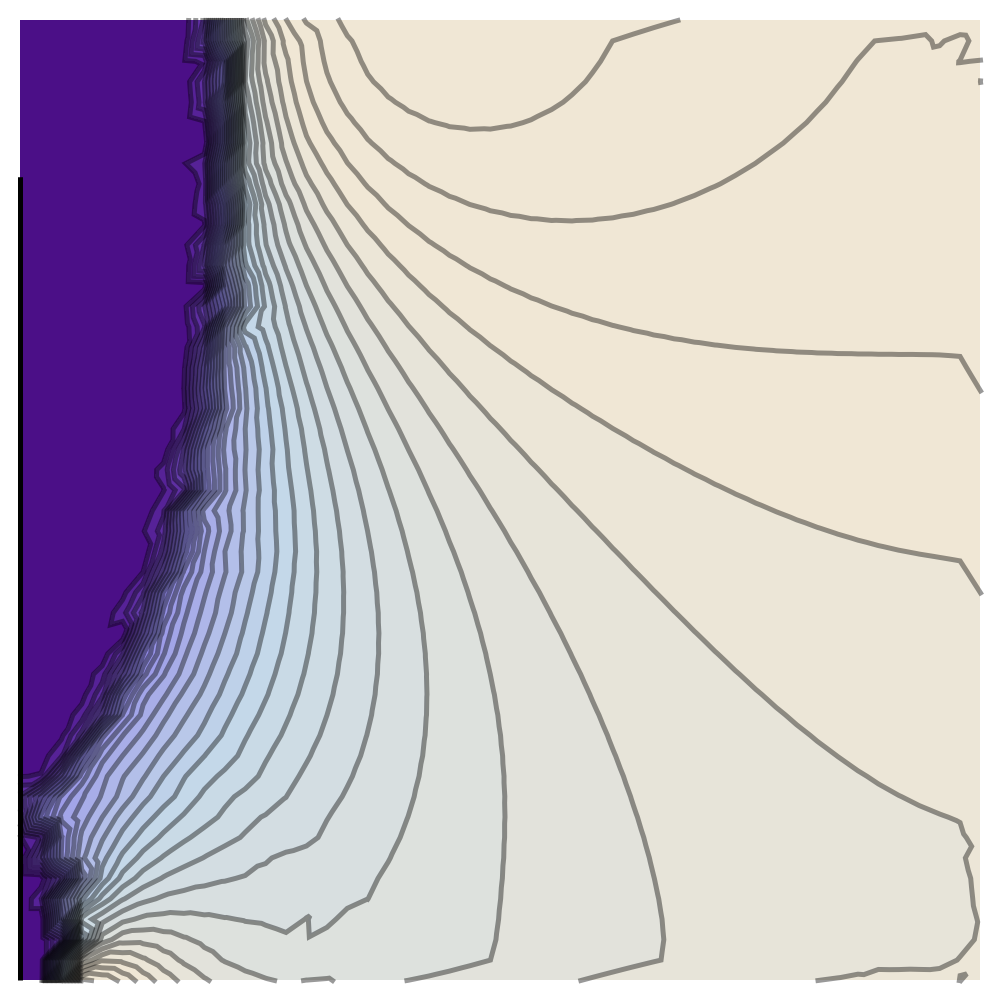}
\includegraphics[width=3.6cm]{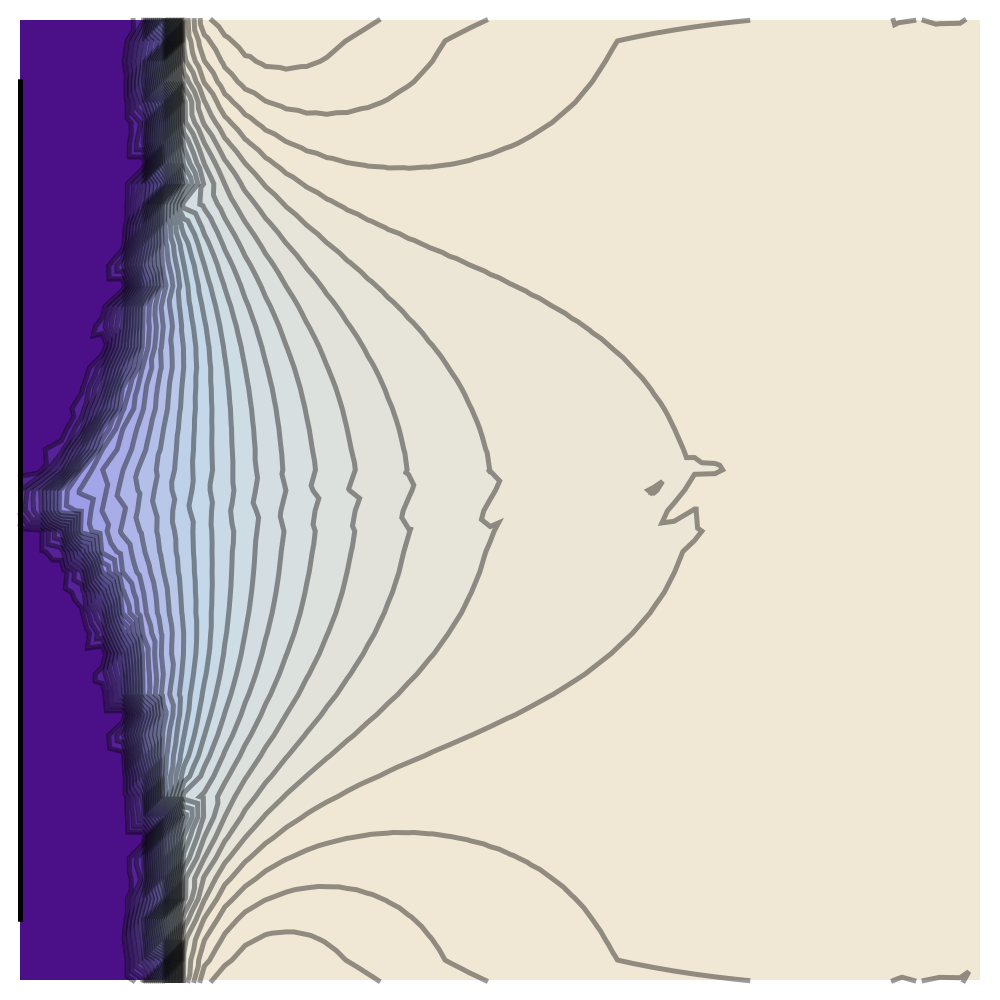}
\includegraphics[width=0.6cm]{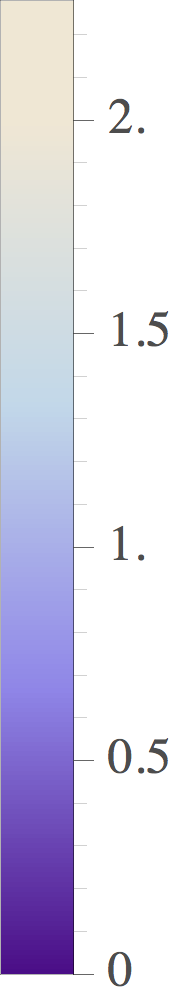}
\\

\includegraphics[width=3.6cm]{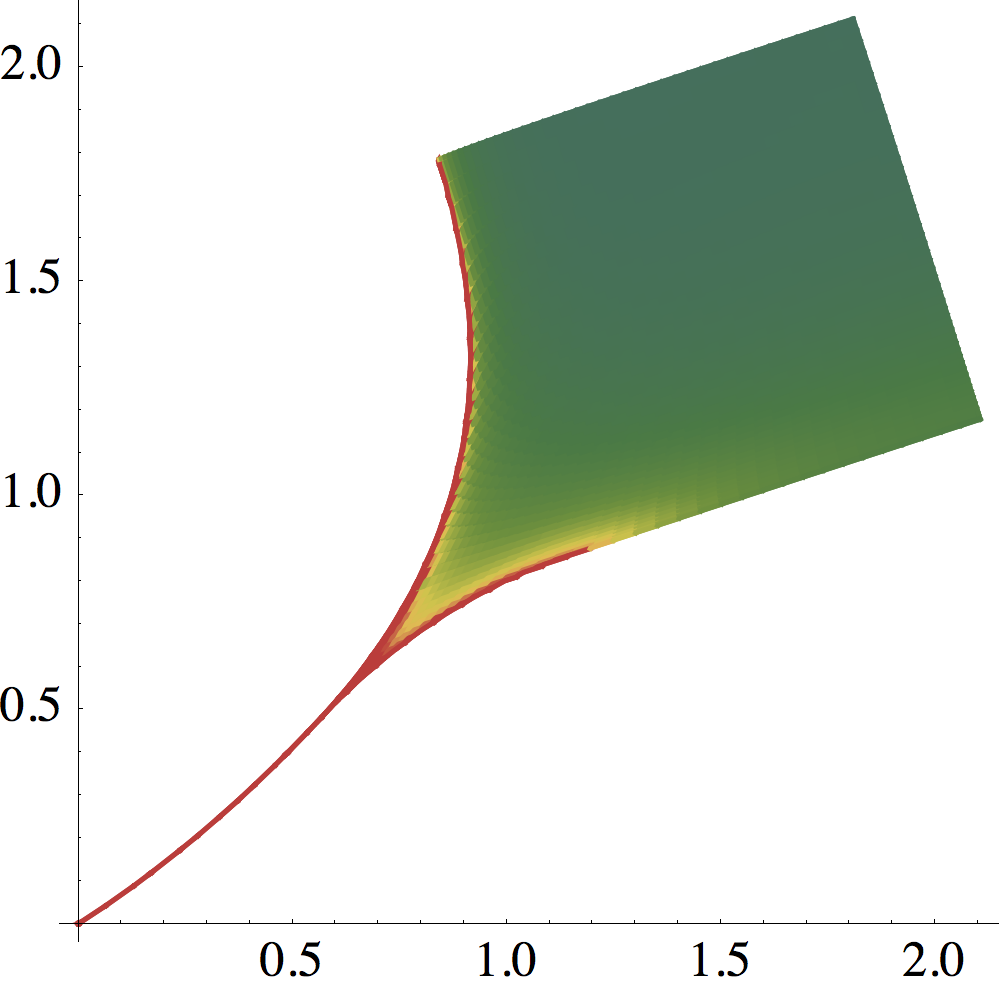}
\includegraphics[width=3.6cm]{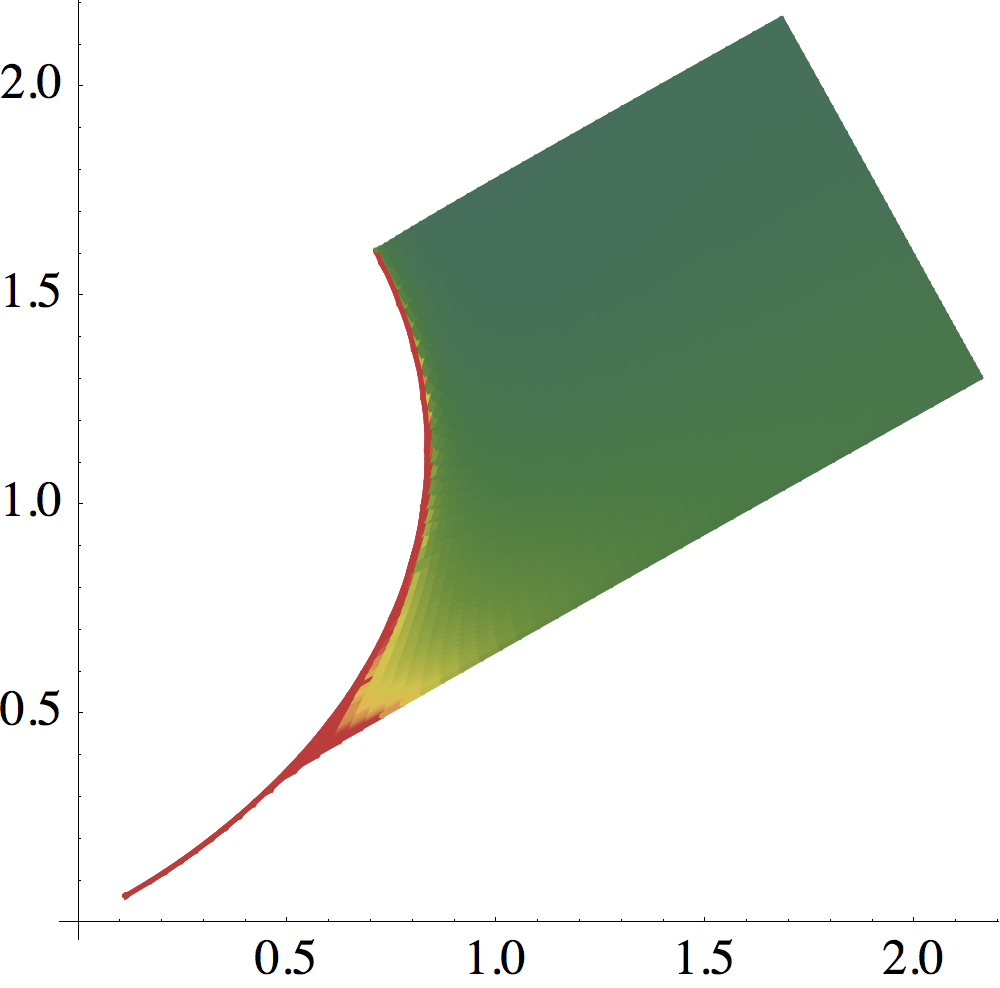}
\includegraphics[width=3.6cm]{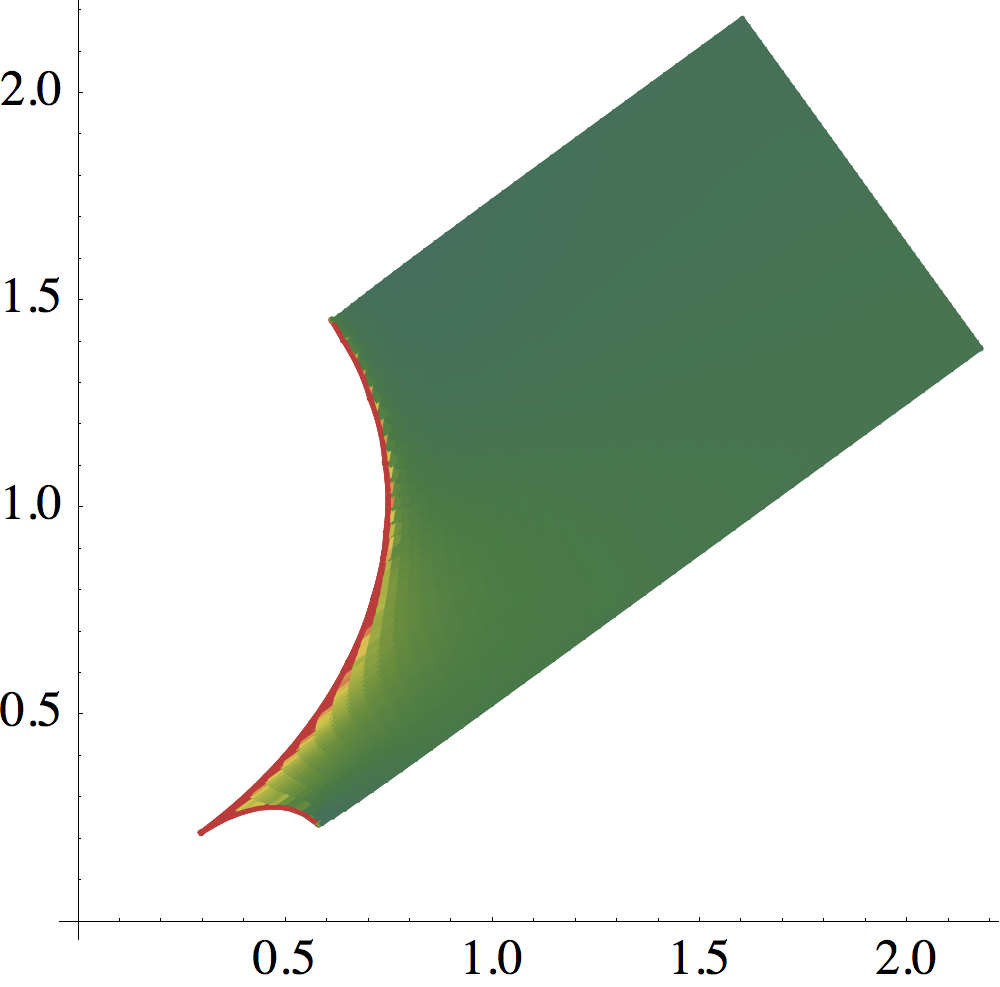}
\includegraphics[width=3.6cm]{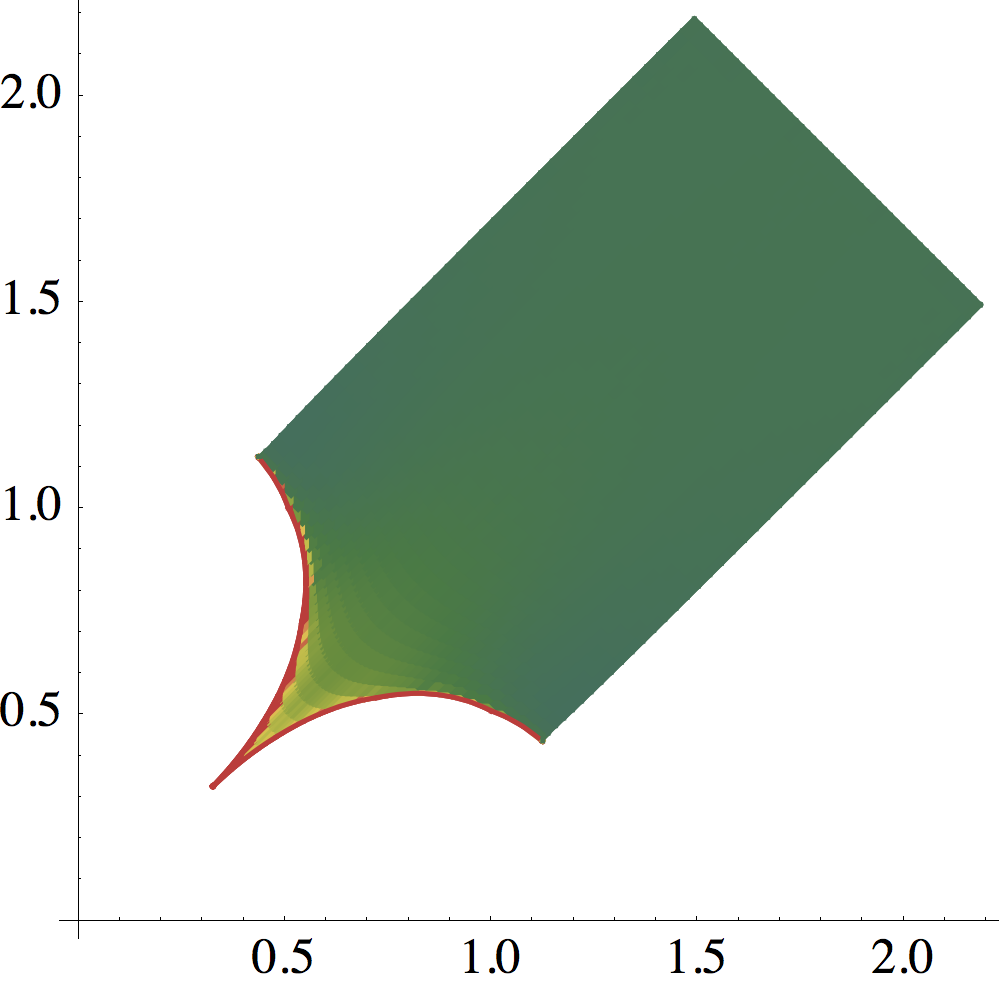}
\includegraphics[width=0.6cm]{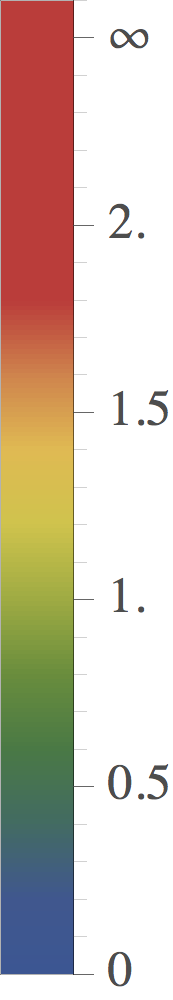}
\\

\includegraphics[width=3.6cm]{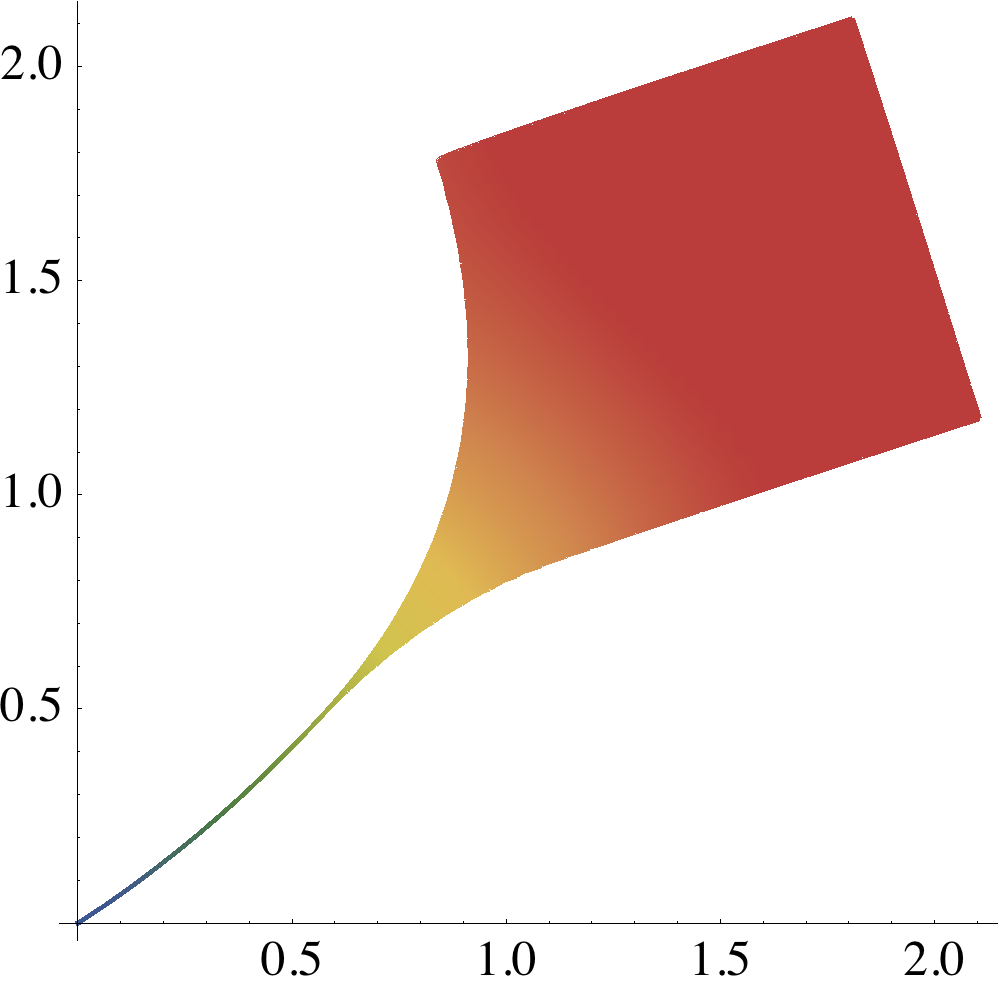}
\includegraphics[width=3.6cm]{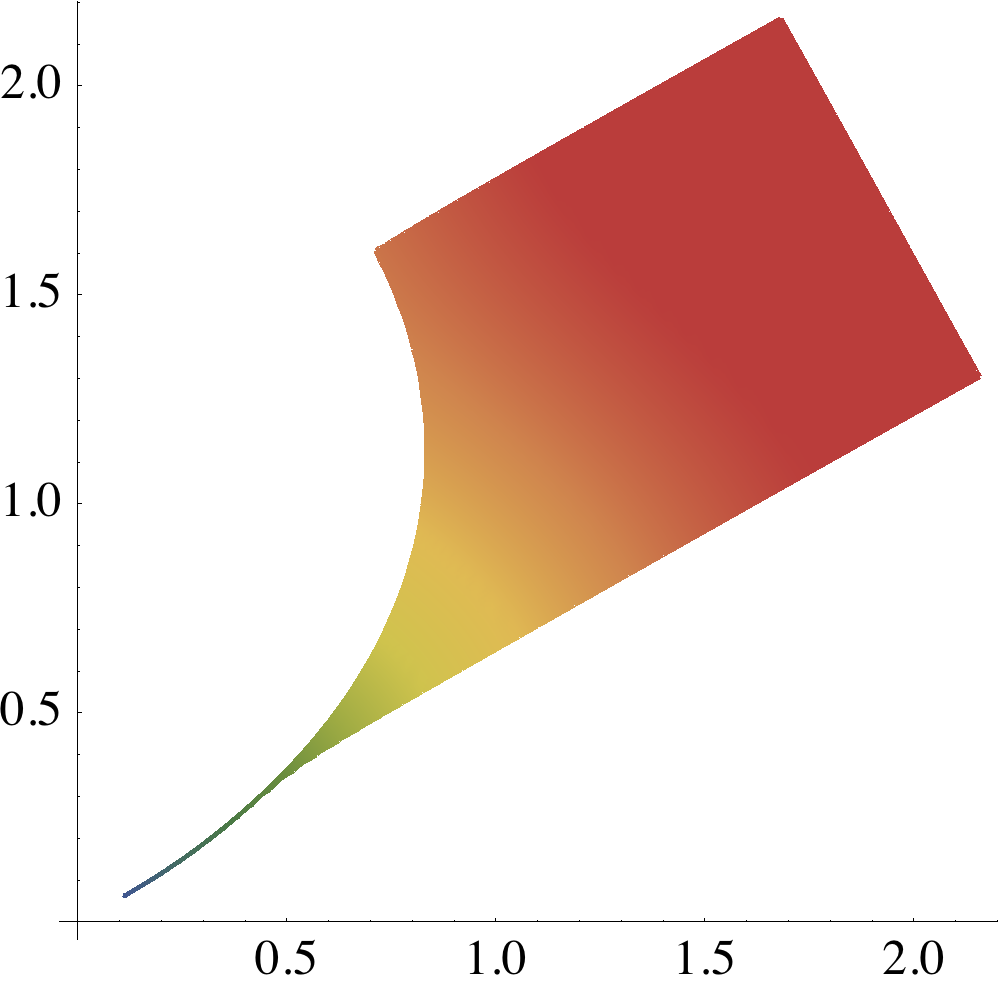}
\includegraphics[width=3.6cm]{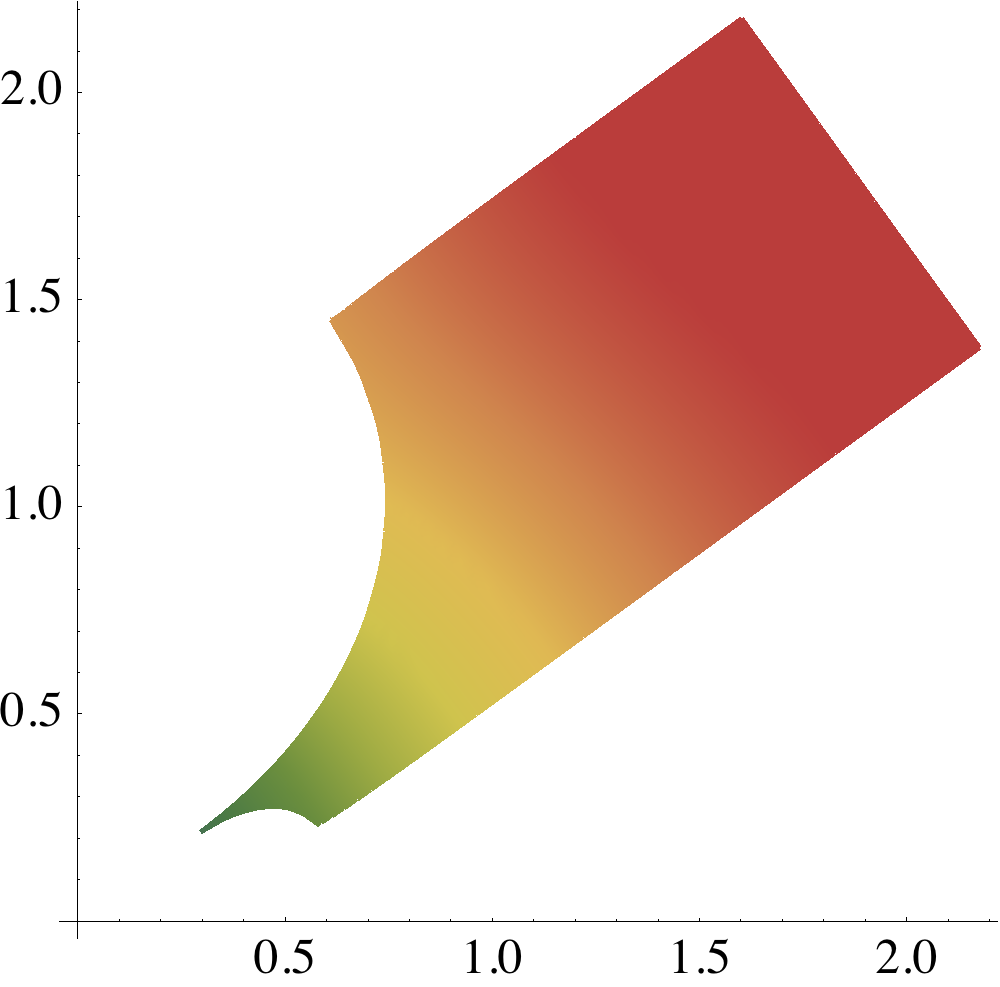}
\includegraphics[width=3.6cm]{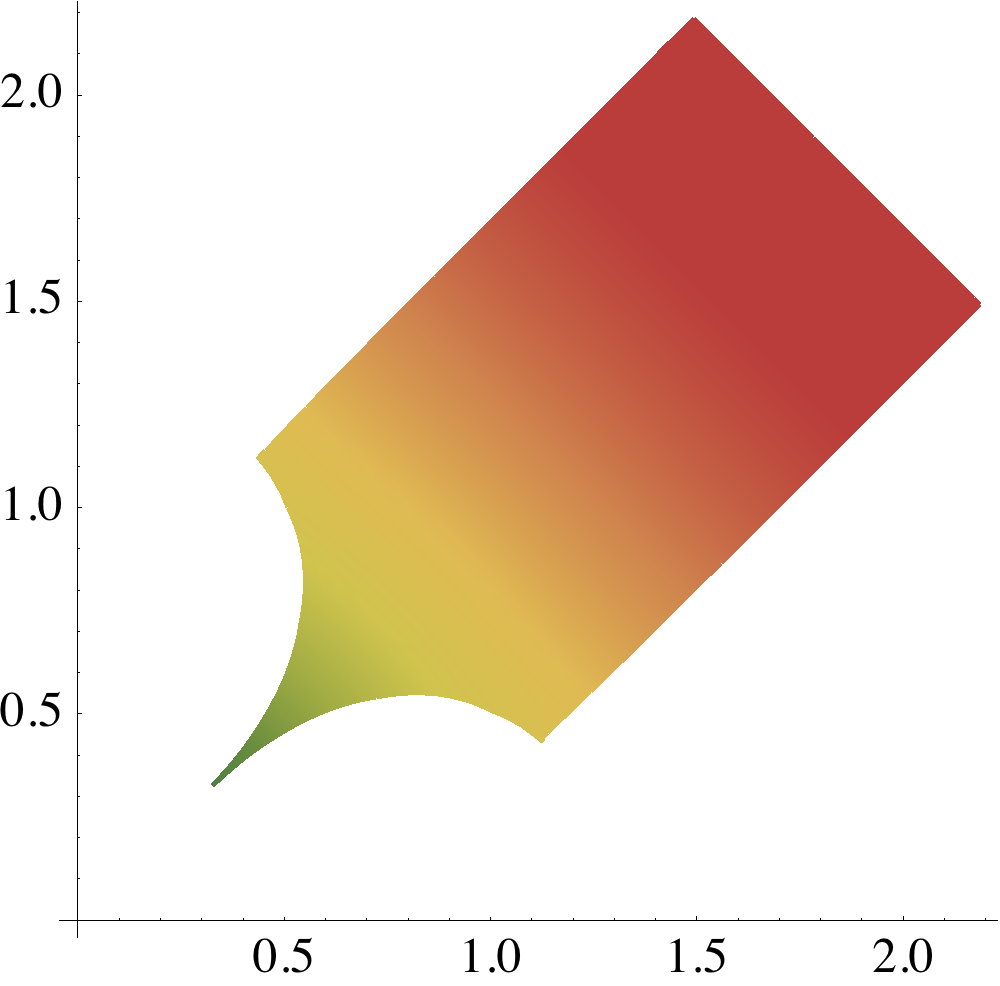}
\includegraphics[width=0.6cm]{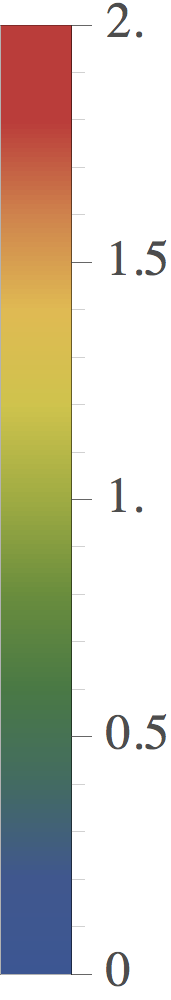}

\caption{
Results of the classical principal agent model, with a uniform density of customers on the domain obtained by rotating the square $[1,2]^2$ around its center $(3/2,3/2)$ by the indicated angle $\theta$. 
From top to bottom: level sets of $U$ (vanishing set indicated in white), level sets of $\det(\Hessian U)$ (allows to discriminate between the three categories of customers), density of products bought (which explodes on some curves), margin of the monopolist.
}
\label{fig:PARotated}
\end{figure}

We also consider variants where the density of customers is uniform on the square $[1,2]^2$ \emph{rotated} by an angle $\theta\in [0, \pi/4]$ around its center, see Figures \ref{fig:DTGC} (left) and \ref{fig:PARotated}.
Our experiments suggest that exclusion occurs iff $\theta \in [0, \theta_0]$, with $\theta_0 \approx 0.47$  rad. Bunching is always present, yet two regimes can be distinguished: the one dimensional product line, associated to the bunching phenomenon, is included in the boundary of the two dimensional one iff $\theta\in [\theta_1, \pi/4]$, with $\theta_0 < \theta_1 \approx 0.55$ rad. Proving mathematically this qualitative behavior is an open problem. 

\begin{figure}
\includegraphics[width=3.9cm]{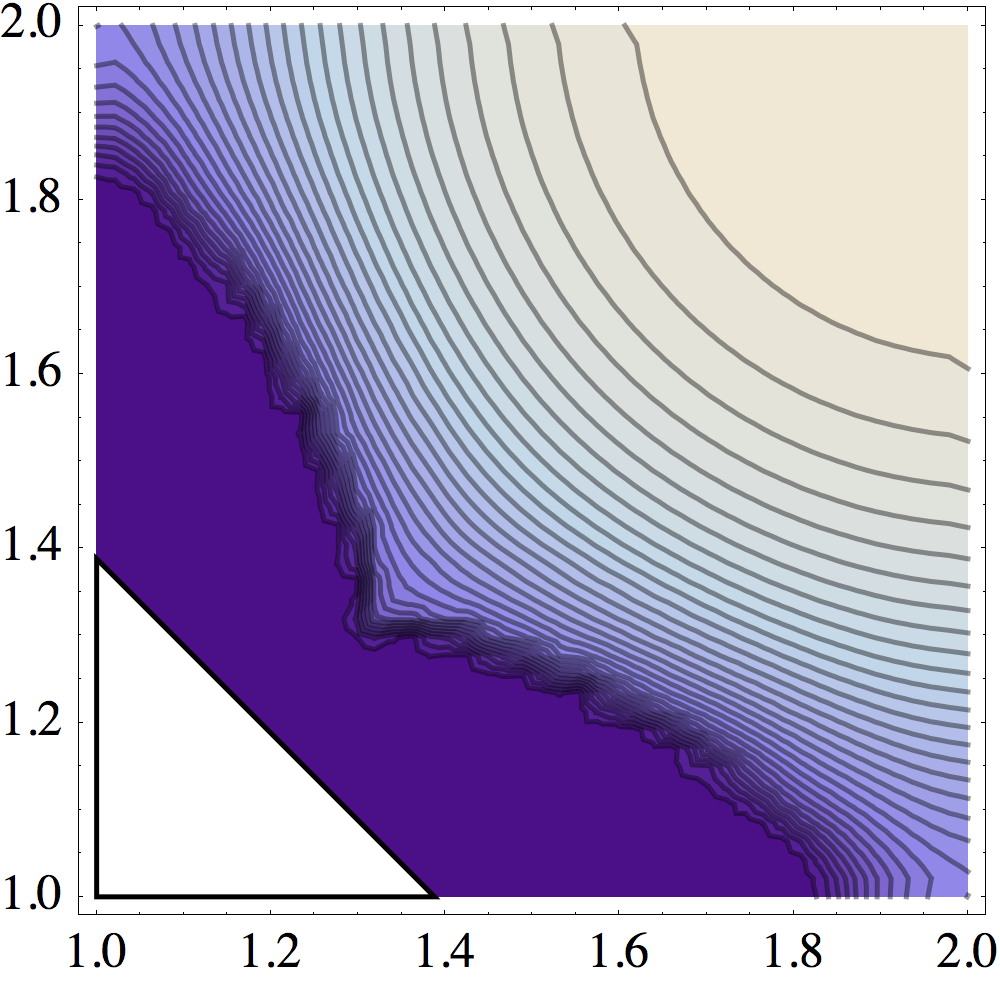}
\includegraphics[width=3.9cm]{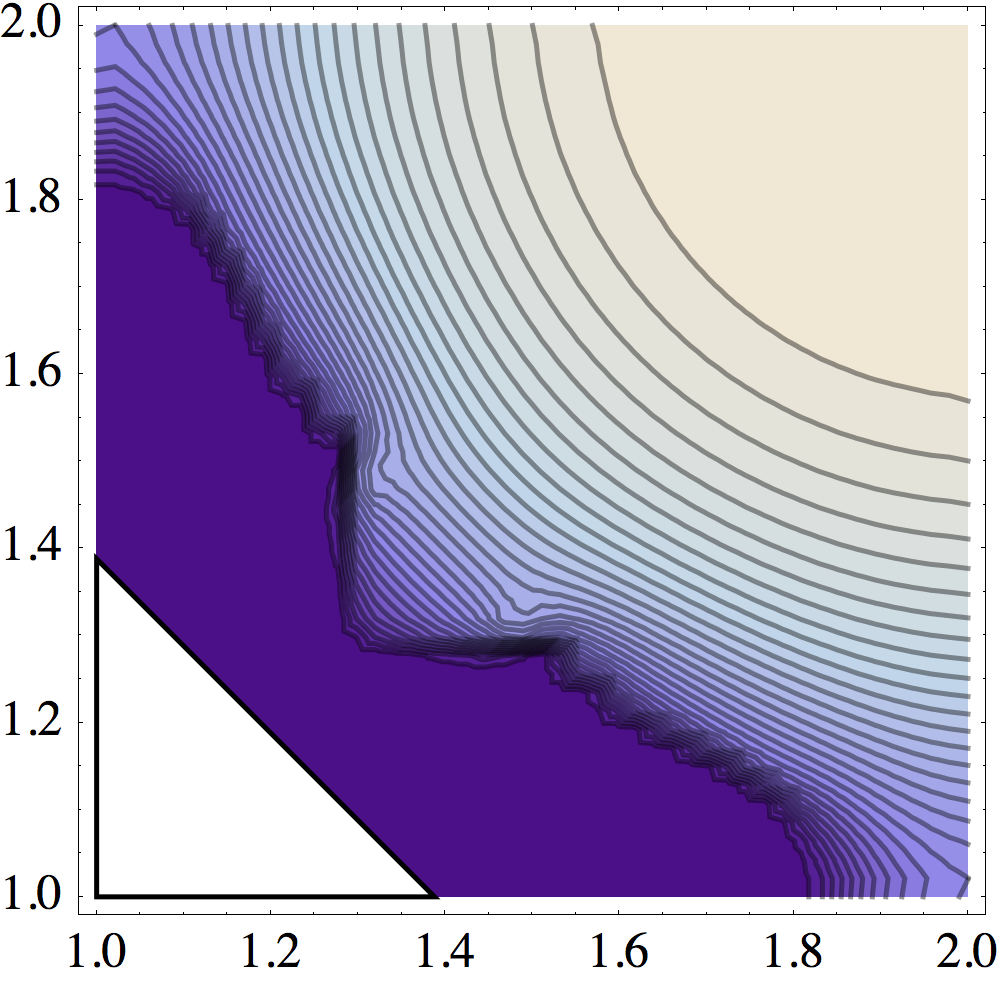}
\includegraphics[width=3.9cm]{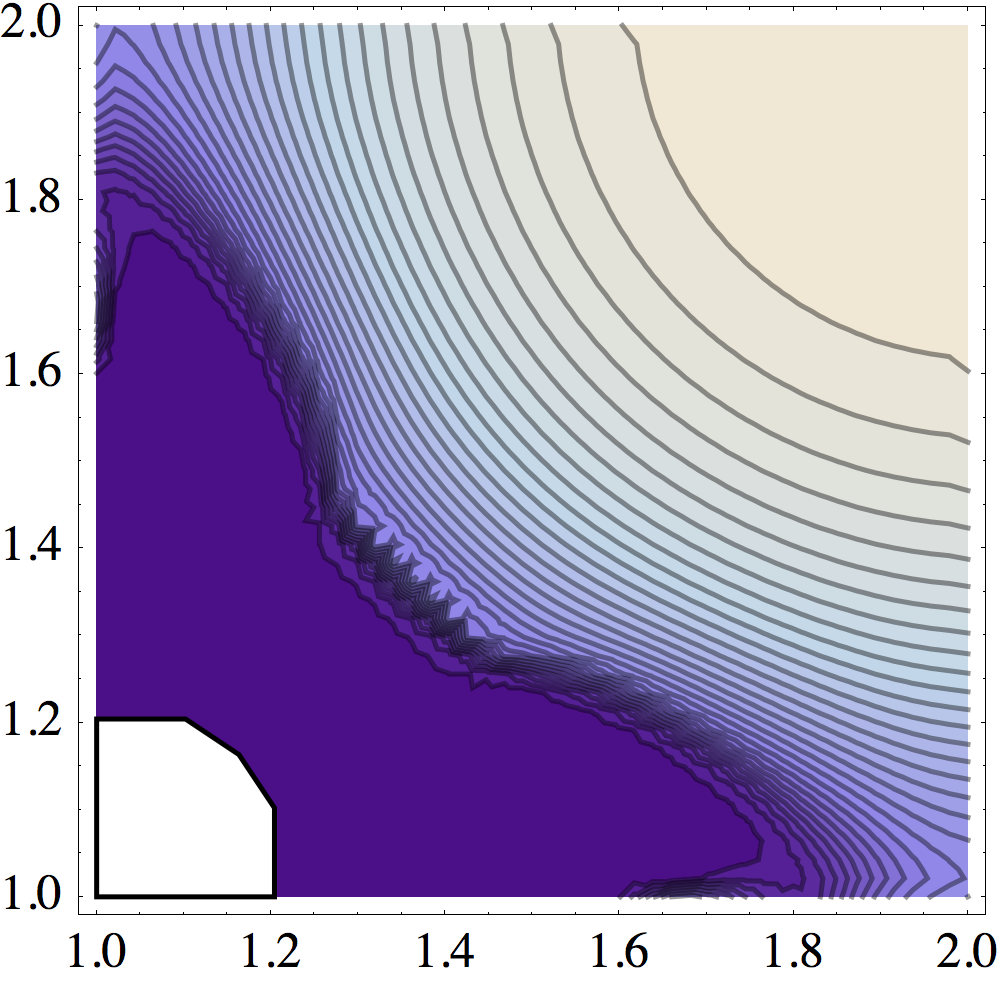}
\includegraphics[width=3.9cm]{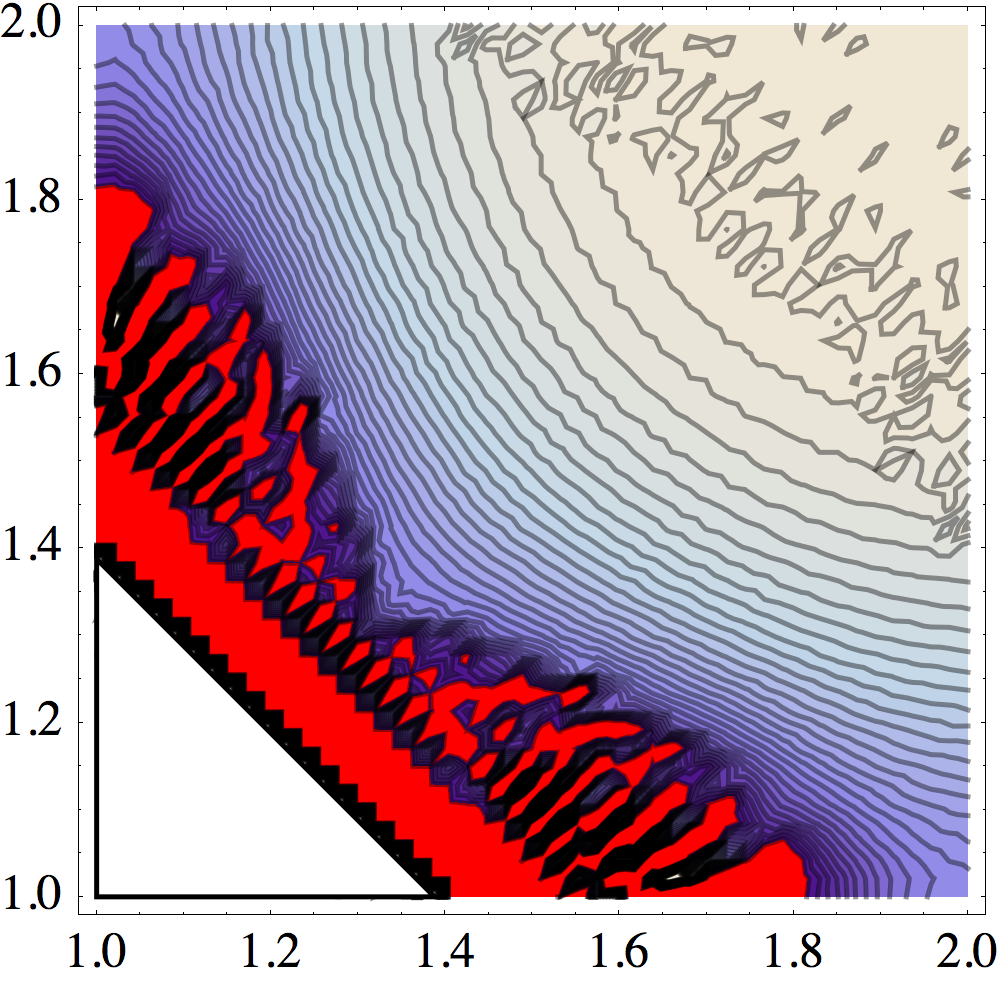}
\caption{
Level set $U<10^{-4}$, in white, approximating the region $U=0$ of excluded customers. Other level sets: $\{k\eta \leq \det(\Hessian U) \leq (k+1) \eta\}$, with $\eta=0.07$. The dark blue one, for $k=0$, approximates the region $\det(\Hessian U) = 0$ of customers subject to bunching. 
Hessian determinant extracted via subgradient measures, see Remark \ref{rem:Subgradients}, and the numerical methods $\Conv$ (left), OF$_2$ (center left) and AM (center right). Right: extraction by taking the determinant of a finite differences discrete Hessian (and here the method $\Conv$); this naïve procedure is unstable and produces negative values (in red).
}
\label{fig:BadHessian}
\end{figure}

\begin{Remark}[Subgradient measure]
\label{rem:Subgradients}
Studying the ``bunching'' phenomenon requires to estimate the hessian determinant $\det(\Hessian U) \geq 0$ of the solution $U$ of \eqref{eq:TotalProfit}, and to visualize the degenerate region $\det(\Hessian U)=0$. The hessian determinant also appears in the density of product sales \eqref{eq:ProductMarginDistribution}. These features need to be extracted from a minimizer $u \in \Conv(X)$ of a finite differences discretization of \eqref{eq:TotalProfit}, which is a delicate problem since (i) the hessian determinant is a ``high order'' quantity, and (ii) equality to zero is a numerically unstable test.
Naïvely computing a discrete hessian $H_u$ via second order finite differences, we obtain an oscillating, non-positive and overall imprecise approximation $\det(H_u)$, see Figure \ref{fig:BadHessian} (right).
 

The following approach gave better results, see Figure \ref{fig:BadHessian} (left):
compute the largest convex $U: \Hull(X) \to \R$ such that $U \leq u$ on $X$ (if $u \in \Conv(X)$, then $U=\interp_\cT u$ for any $u$-Delaunay triangulation).
Then for all $x \in X \sm \partial \Hull(X)$
\begin{equation*}
h^2 \det(\Hessian U (x)) \approx | \{\nabla U (x+e); \, \|e\|_\infty \leq h/2\} | \approx  | \{\nabla \hat U (x+e);\,  \|e\|_\infty \leq h/2\} | = |\partial_x \hat U|,
\end{equation*}
where $\partial$ denotes the sub-gradient, set-valued operator on convex functions, and $|\cdot|$ the two dimensional Lebesgue measure. The sub-gradient sets $\partial_x \hat U$ are illustrated on Figure \ref{fig:DTGC}.
\end{Remark}

\paragraph{Product bundles and lottery tickets.}
Two types of products $P_1$, $P_2$ are considered, which the consumer of characteristics $x=(x_1,x_2)\in \R^2$ respectively values $x_1$ and $x_2$. The two products are indivisible, and consumers are not interested in buying more than one of each. The monopolist sells them in bundles $q = (q_1,q_2) \in \{0,1\}^2$ which characteristics are the presence ($q_i=1$) of product $P_i$, or its absence ($q_i=0$), for $i\in \{1,2\}$. In order to maximize profit, the monopolist also considers probabilistic bundles, or lottery tickets, $q\in [0,1]^2$ for which the product $P_i$ has the probability $q_i$ of being present. This is consistent with \eqref{eq:Utility}, provided consumers are risk neutral. Production costs are neglected, so that $\Cost(q) = 0$ if $q\in [0,1]^2$ and $\Cost(q) = +\infty$ otherwise. Three different customer densities were considered, see below.
The qualitative property of interest is the presence, or not, of probabilistic bundles in the monopolist's optimal strategy.

\begin{figure}
\begin{tabular}{|c||cc||cc|}
\hline
\includegraphics[width=3cm]{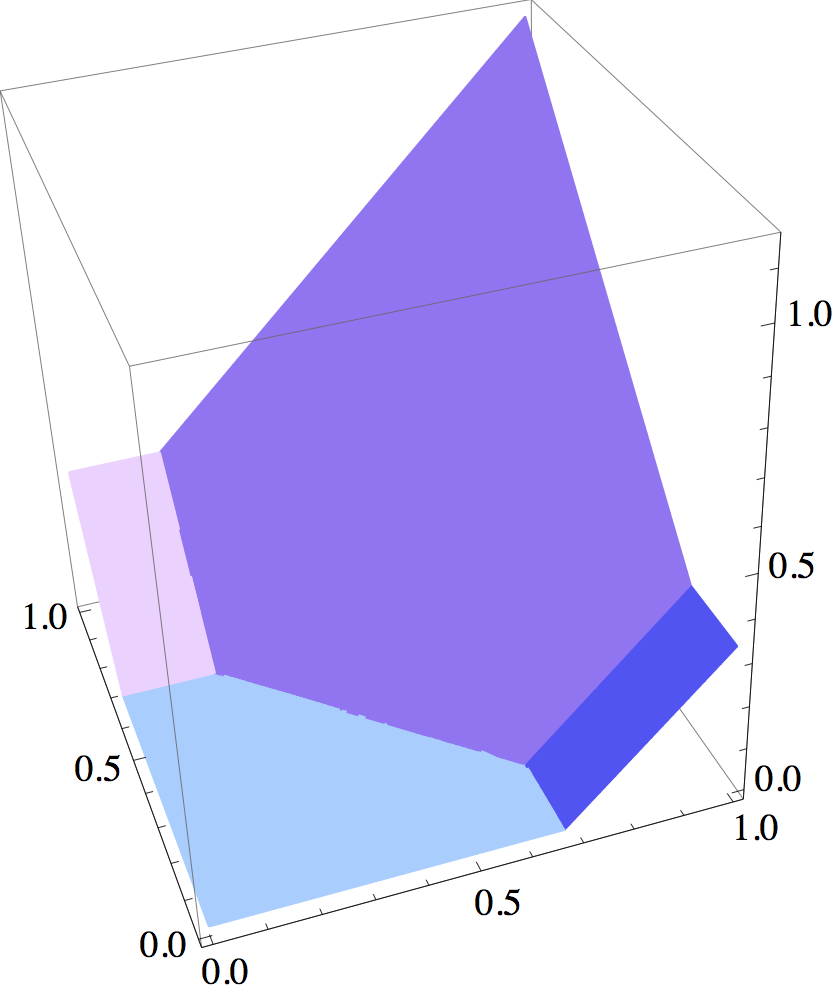} &
\includegraphics[width=2cm]{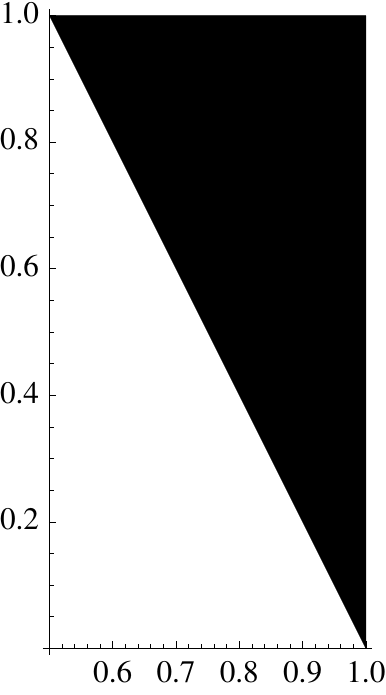} &
\includegraphics[width=2.5cm]{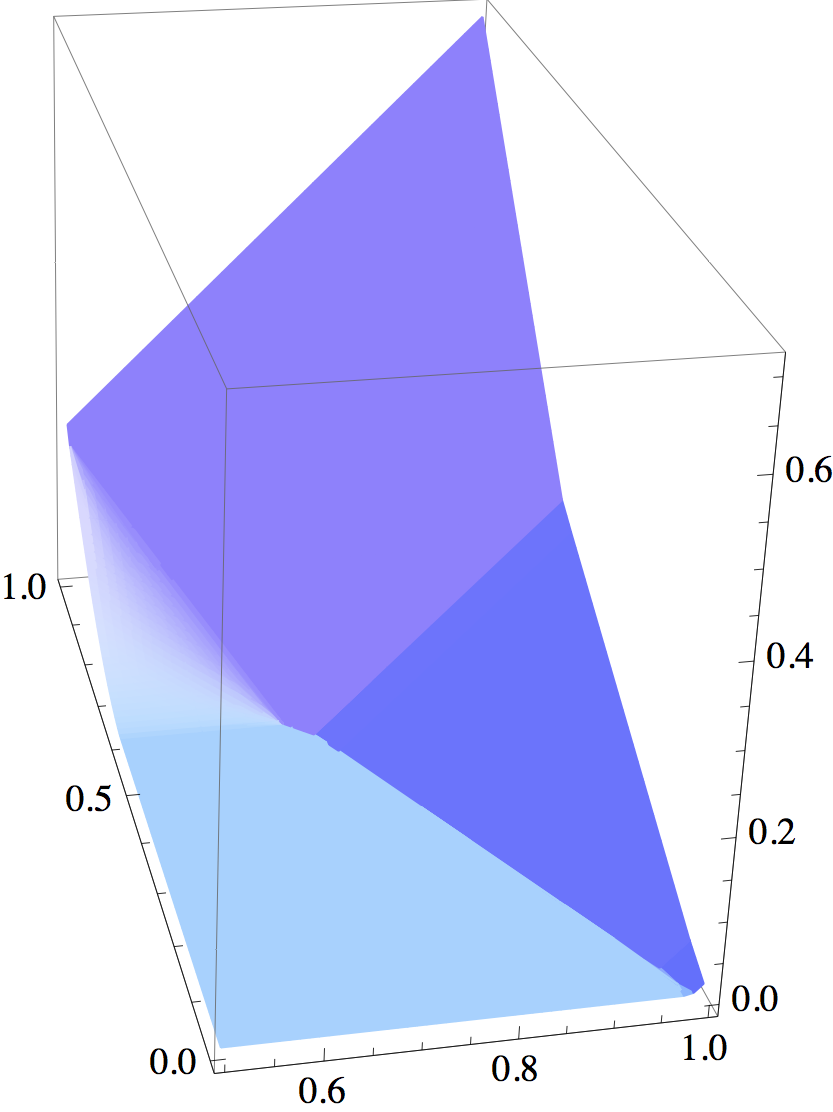} &
\includegraphics[width=3cm]{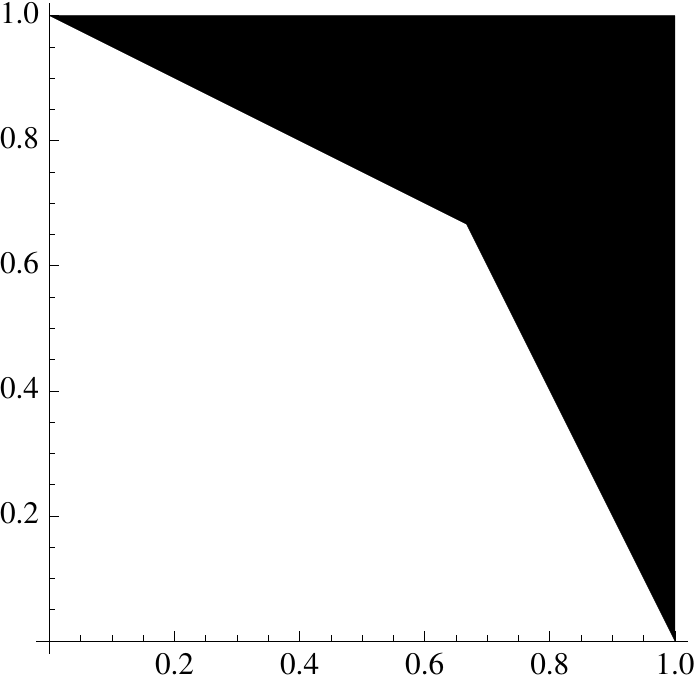} &
\includegraphics[width=3.5cm]{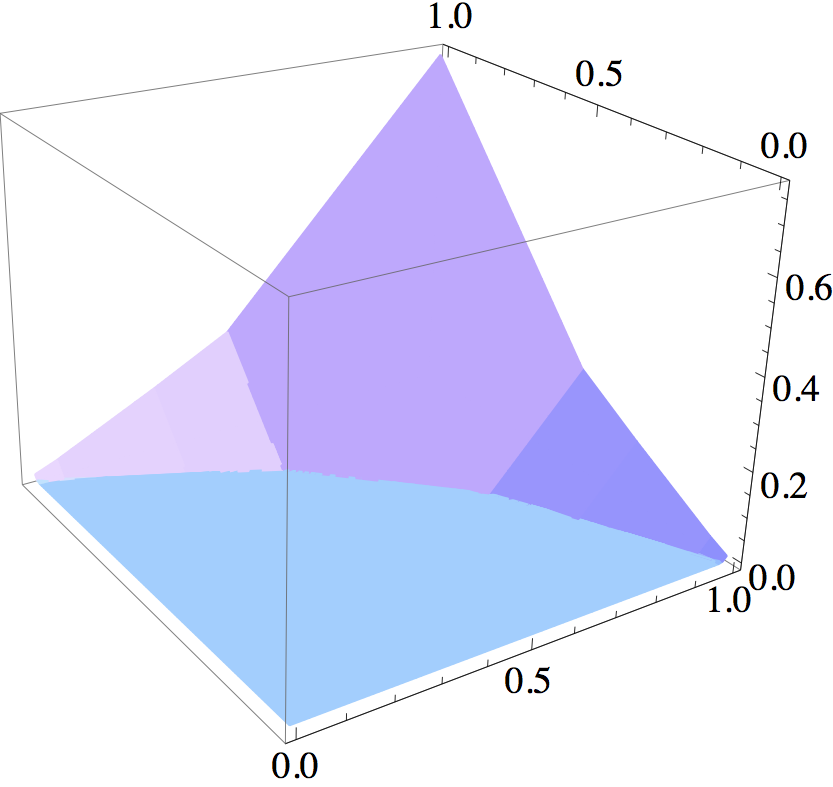}\\
\hline
\end{tabular}
\caption{Three dimensional plot of the optimal $U$ for the product-bundles variant of the monopolist problem, with respect to various customer distributions. Left : uniform distribution on $[0,1]^2$. Center, and right: distribution uniform on the illustrated black polygon.}
\label{fig:Bundles}
\end{figure}

\begin{itemize}
\item Uniform customer density on $[0,1]^2$. We recover the known exact minimizer \cite{Manelli:2006ib}:
\begin{equation*}
U(x,y) := \max\{0,\, x-a,\, y-a,\, x+y-b\}, \text{ with } a:=2/3, \text{ and } b:=(4-\sqrt 3)/2,
\end{equation*}
up to numerical accuracy, see Figure \ref{fig:Bundles} (left).
This optimal strategy does \emph{not} involve lottery tickets: $\nabla U(x,y) \in \{0,1\}^2$, wherever this gradient is defined. 
The uselessness of lottery tickets is known for similar 1D problems and was thought to extend to higher dimension, 
until the following two counter-examples were independently found \cite{Manelli:2006ib,Thanassoulis:2004uy}.
\item 
Uniform customer density on the triangle $T := \{(x,y) \in [0,1]^2; \, x+y/2\geq 1\}$. The monopolist strategy associated to 
\begin{equation}
\label{eq:BestBundleT}
U(x,y) := \max \{0, \, x+y/2-1, \, x+y-b\}, \text{ with } b=1+1/(2\sqrt 3),
\end{equation}
which involves the lottery ticket $(1,1/2)$, yields better profits that any strategy restricted to deterministic bundles \cite{Manelli:2006ib}. The triangle $T$, and the numerical best $U$, are illustrated on Figure \ref{fig:Bundles} (center).
These experiments suggest that \eqref{eq:BestBundleT} is a%
\footnote{%
Optimal solutions of \eqref{eq:TotalProfit} are not uniquely determined outside the customer density support $T=\supp(\mu)$.
}
globally optimal solution. 
\item 
Uniform customer density on the kite shaped domain $\{(x,y) \in [0,1]^2;\,  x+y/2 \geq 1 \text{ or } x/2+y \geq 1\}$, see Figure \ref{fig:Bundles} (right). The optimal monopolist strategy is proved in \cite{Thanassoulis:2004uy} to involve probabilistic bundles, but it is not identified. Our numerical experiments suggest that it has the form
\begin{equation*}
U(x,y) := \max \{0, \, x+y/2-1, \, x/2+y-1, \, x+y-b\},
\end{equation*}
which involves the lottery tickets $(1,1/2)$ and $(1/2,1)$. Under this assumption, the optimal value $b=1+1/(3\sqrt 2)$ is easily computed.
\end{itemize}
Numerous qualitative questions remain open. Is there a distribution of customers for which the optimal strategy involves a continuum of distinct lottery tickets $\{(1,\alpha);\, \alpha_0\leq \alpha \leq \alpha_1\}$ ?

%

\paragraph{Pricing of risky assets.}

\begin{figure}
\centering
\includegraphics[width=3.8cm]{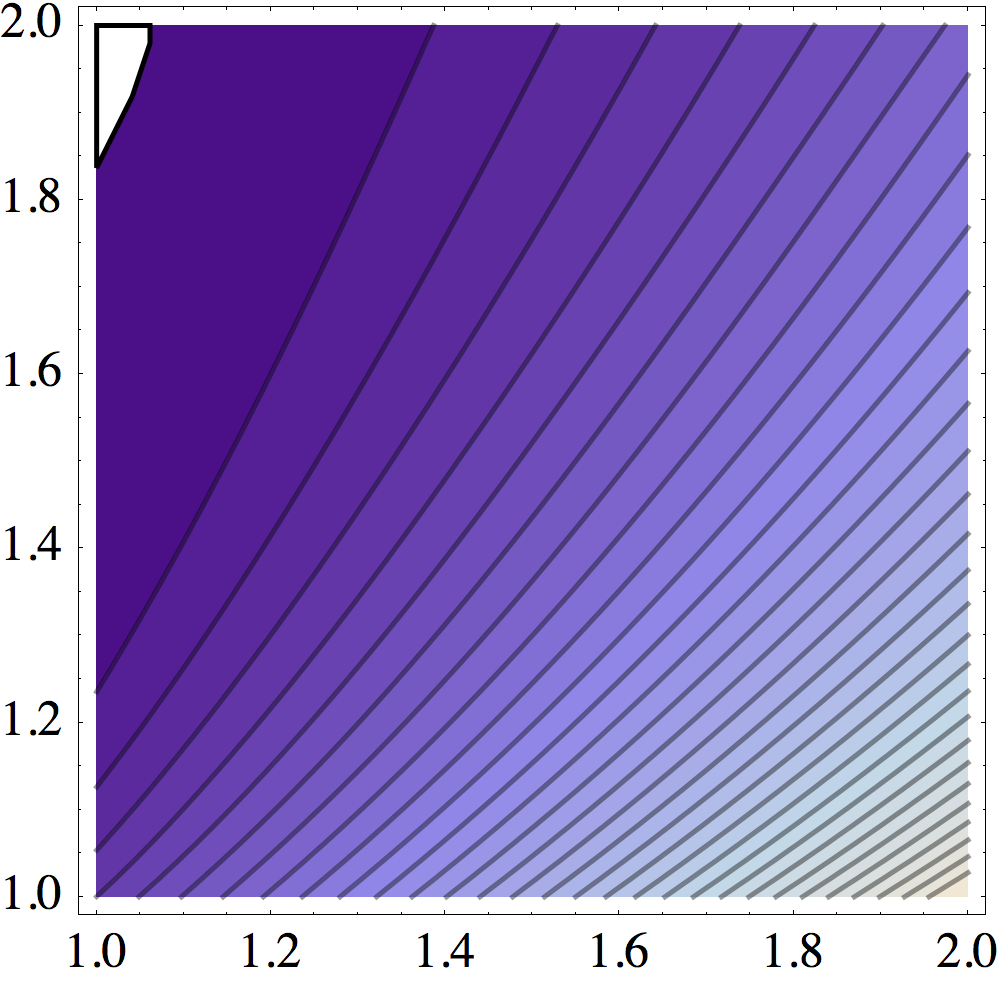}
\includegraphics[width=3.8cm]{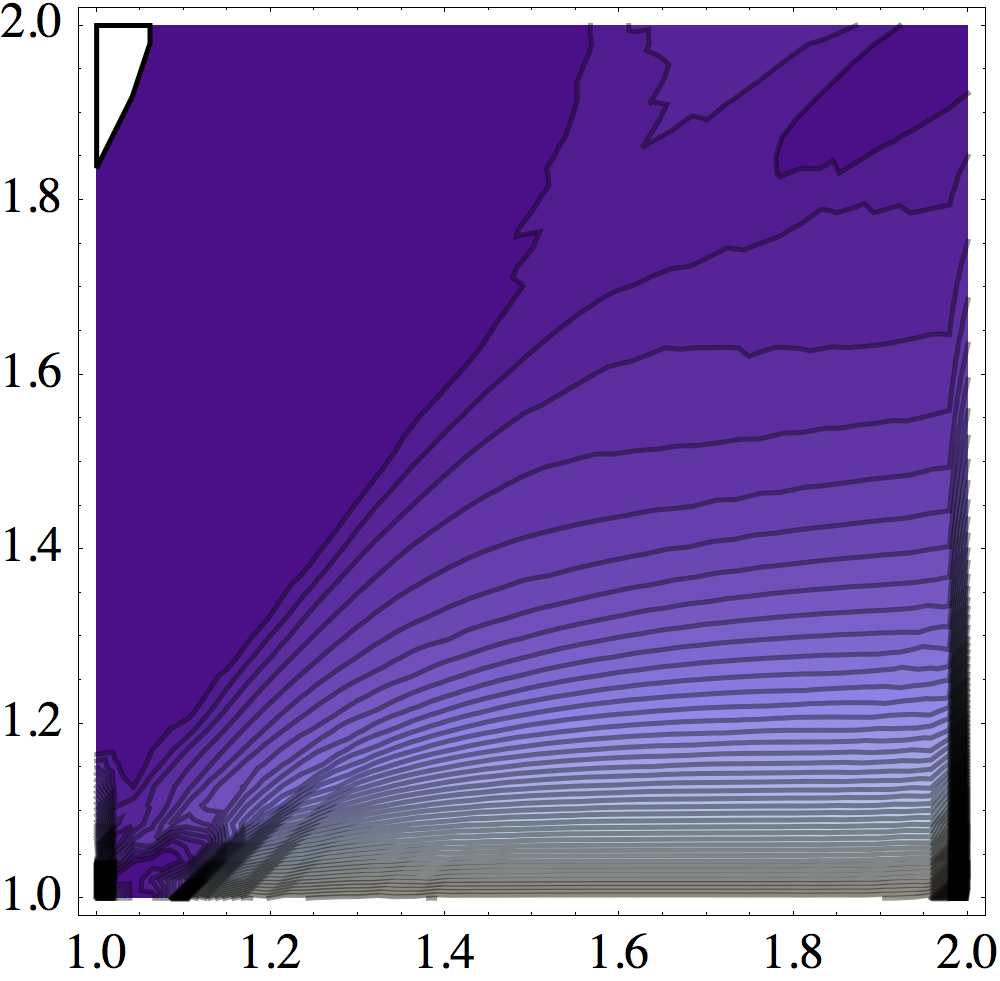}
\includegraphics[width=2.4cm]{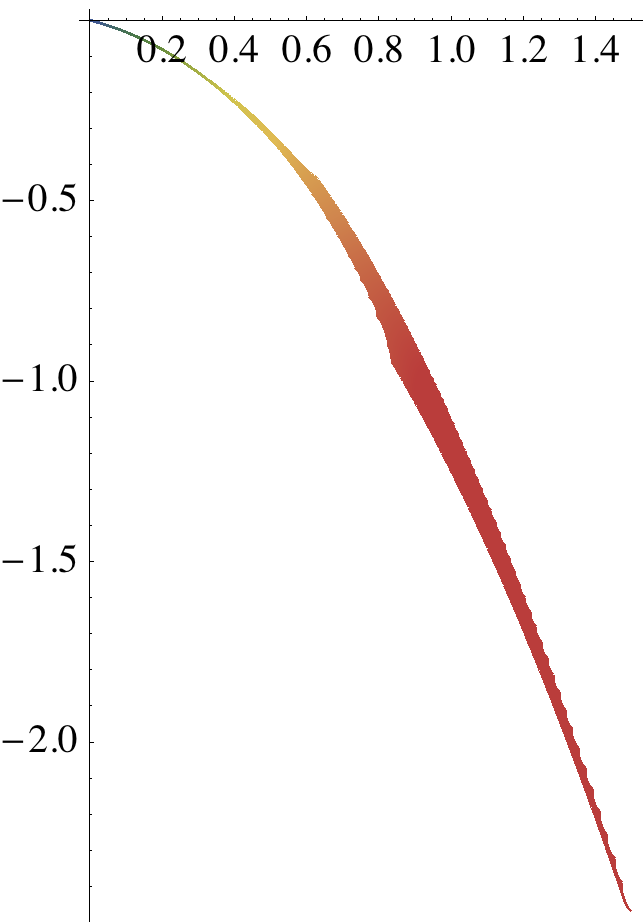}
\includegraphics[width=4cm]{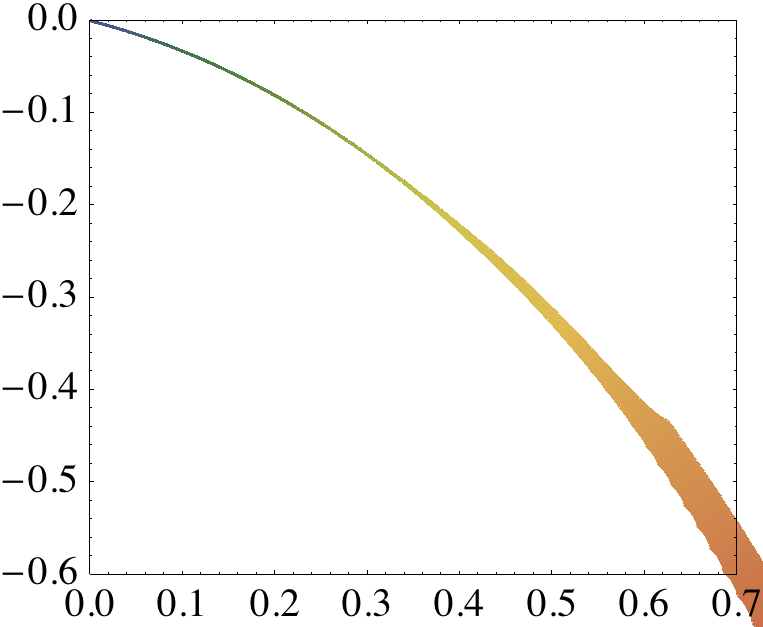}
\caption{Optimal pricing of risky assets, with the parameters $\alpha=1$, $\xi=1$, and a uniform customer density on $[1,2]^2$.
Left: level sets of the estimated solution $U$ of \eqref{eq:Assets}, with exclusion region $U=0$ in white. Center left: level sets of $\det(\Hessian U)$, the darkest one is the estimated bunching region $\det(\Hessian U)=0$. Center right, and right (detail): optimal product line, colored with the monopolist margin.}
\label{fig:Assets}
\end{figure}

A more complex economic model is considered in \cite{Carlier:2007gy}, where financial products, characterized by their expectancy of gain and their variability, are sold to agents characterized by their risk aversion and their initial risk exposure. 
We do not give the details of this model here, but simply point out that it fits in the general framework of \eqref{eq:TotalProfit} with the cost function 
$\Cost(a,b) := - \alpha (\xi a + \sqrt{-(a^2+b)})$,
 if $a^2+b \leq 0$, and $+\infty$ otherwise, where $\xi \in \R$ and $\alpha \geq 0$ are parameters, see Example 3.2 in \cite{Carlier:2007gy}.
Observing that, for $a^2+b < 0$
\begin{equation*}
\det(\Hessian \Cost(a,b)) = \frac 1 {4 (a^2+b)^2},
\end{equation*}
we easily obtain that this cost is convex\footnote{This property was not noticed in the original work \cite{Carlier:2007gy}.}. The lack of smoothness of the square-root appearing in the cost function is a potential issue for numerical implementation, hence the problem \eqref{eq:TotalProfit} is reformulated using an additional variable $V$ subject to a (optimizer friendly) conic constraint
\begin{equation}
\label{eq:Assets}
\max \left\{ \int \left(\<\nabla U,z\> - U + \alpha\xi\, \partial_x U + \alpha V\right) d \mu ; \, U \in \Conv_0(\R^2), \, V^2 + (\partial_x U)^2+ \partial_y U \leq 0\right\}.
\end{equation}
A numerical solution, presented Figure \ref{fig:Assets}, displays the same qualitative properties (Desirability of exclusion, Bunching) as the classical monopolist problem with quadratic cost.

\subsection{Comparison with alternative methods}
\label{sec:Comparison}
We compare our implementation of the constraint of convexity with alternative methods that have been proposed in the literature. The compared algorithms are the following: 

%
%

\begin{itemize}
\item (Adaptive constraints)
The optimization strategy ($\Conv$) described in Algorithm \ref{algo:SuperCones} section \S \ref{sec:strategy}, based on the hierarchy of cones $\Conv(\cV)$, and used in our numerical experiments \S \ref{sec:Monopolist}. The adaptation ($\DConv$) of this strategy to the hierarchy $\DConv(\cV)$ of cones of ``directionally convex'' functions, see Appendix \ref{sec:Directional}.
\item (Local constraints)
The approach of Aguilera and Morin (AM, \cite{Aguilera:2008uq}) based on semi-definite programming. A method of Oberman and Friedlander (OF$_2$, OF$_3$, \cite{Oberman:2011wi}), where OF$_k$ refers to minimization over the cone $\DConv(\cV_k)$ associated to the fixed stencil $\cV_k(x) := \{e \in \cV_{\max}(x); \, \|e\| \leq k\}$. A modification of OF$_3$ by Oberman (Ob$_3$, \cite{Oberman:2011wy}), with additional constraints ensuring that the output is truly convex.
\item (Global constraints)
Direct minimization over the full cone $\Conv(X)$, as proposed by Carlier, Lachand-Robert and Maury (CLRM, \cite{Carlier:2001tq}). Minimization over $\GradConv(X)$, see \eqref{def:GradConvX}, following%
\footnote{%
We use the description of $\GradConv(X)$ by $\cO(N^2)$ linear constraints given in \cite{Ekeland:2010tl}, but (for simplicity) not their energy discretization, nor their method for globally extending elements $u \in \GradConv(X)$.
}
Ekeland and Moreno (EM, \cite{Ekeland:2010tl}). 
\end{itemize}

The numerical test chosen is the classical model of the monopolist problem, with quadratic cost, on the domain $[1,2]^2$, see Figure \ref{fig:Monopolist} and \S \ref{sec:Monopolist}. This numerical test case is classical and also considered in \cite{Ekeland:2010tl,MERIGOT:tr,Oberman:2011wi}. It is discretized on a $n \times n$ grid, for different values of $n$ ranging from $10$ to $100$ ($10$ to $50$ for global constraint methods due to memory limitations).

\begin{figure}
\centering
\includegraphics[width=3.9cm]{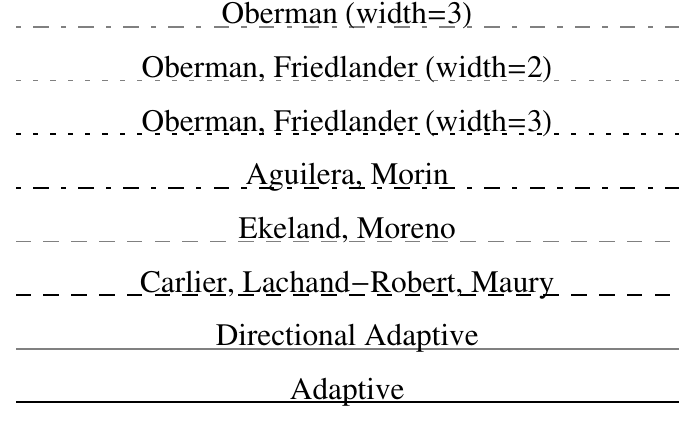}
\includegraphics[width=3.9cm]{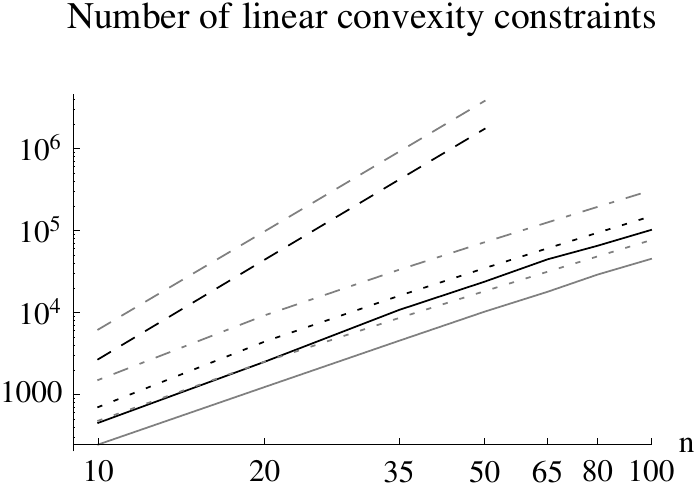}
\includegraphics[width=3.9cm]{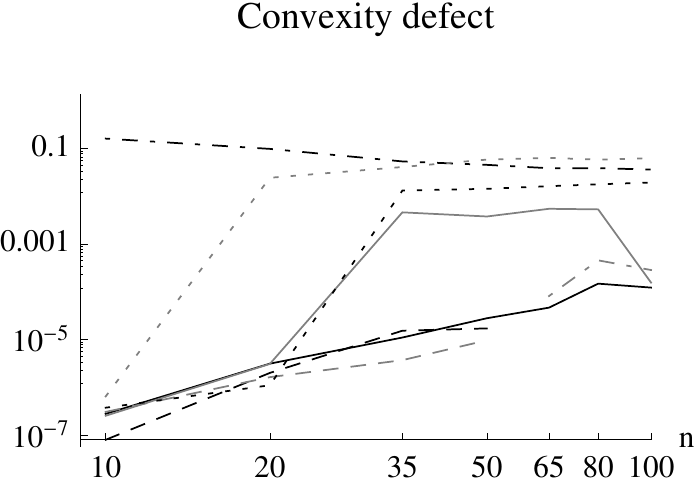}
\includegraphics[width=3.9cm]{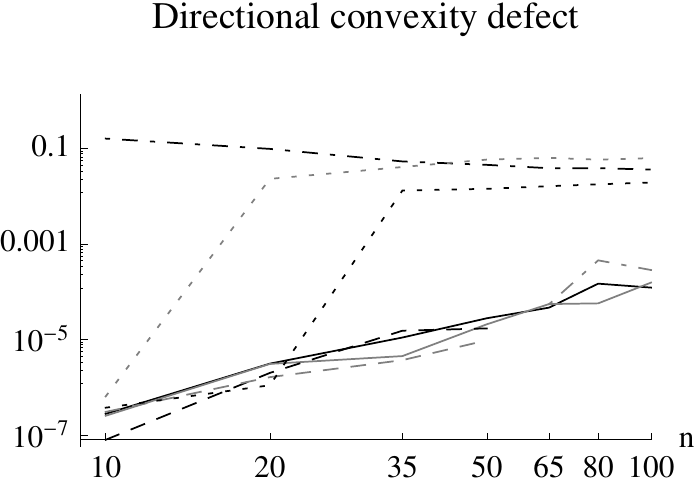}
\vspace{0.2cm}
\\
%
\includegraphics[width=4.1cm]{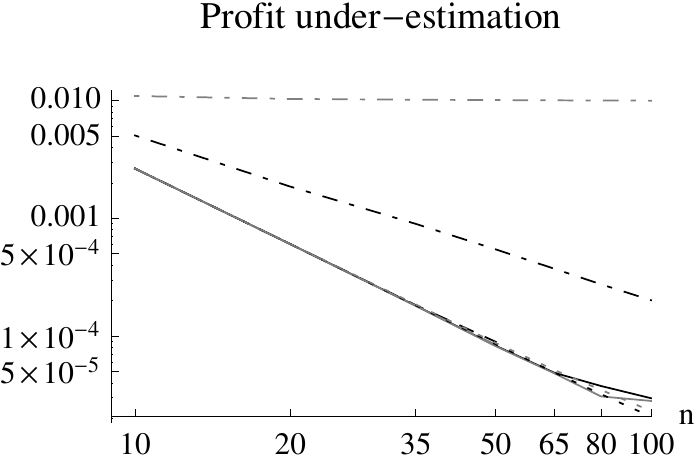}
\includegraphics[width=4.1cm]{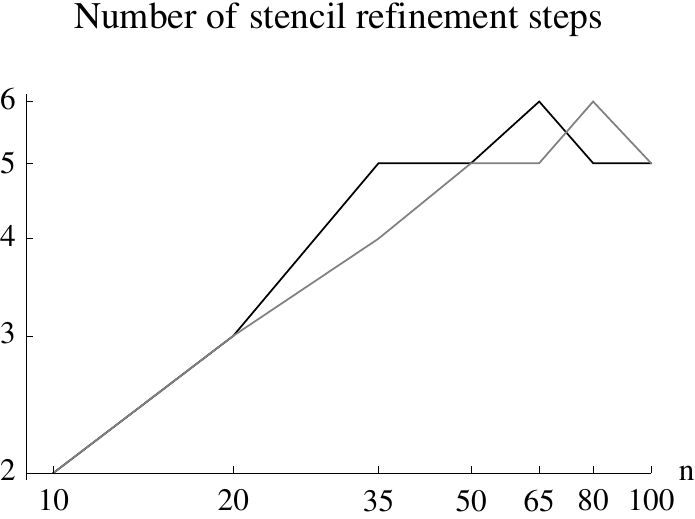}
\includegraphics[width=4.1cm]{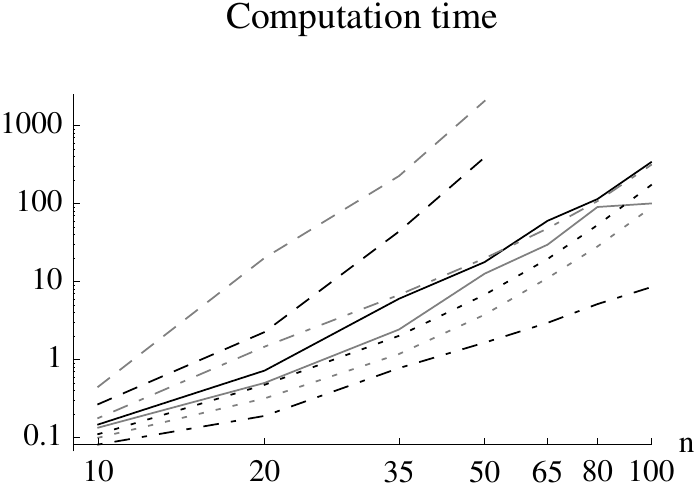}
\vspace{0.2cm}
\\

\begin{tabular}{|c|cccccccc|}
\hline
$n=50$ $\sm$ Method & 
$\Conv$ & 
$\DConv$ & 
AM & 
OF$_2$ & 
OF$_3$ & 
Ob$_3$ &
CLRM & 
EM 
\\
\hline
Constraints $\times 10^{-3}$  &  
24 & 
10 & 
NA & 
18 & 
35 & 
72 & 
1738 & 
3803  
%
\\
Defect $\times 10^{3}$ & 
0.03 & 
3.8 & 
46 & 
59 & 
14 & 
0 & 
0.02& 
0.01  
\\
Profit under estimation $\times 10^{3}$& 
0.11 & 
0.11 & 
0.57 & 
0.11 & 
0.11 & 
10 & 
0.12 & 
0.11  
\\
Computation time & 
18s & 
13s & 
1.7s & 
3.8s & 
6.8s & 
20s & 
391s & 
2070s  
\\ \hline
\end{tabular}

\caption{Comparison of different numerical methods for the classical Monopolist problem.}
\label{fig:Comparison1}
\label{fig:Comparison2}
\end{figure}

The number of linear constraints of the optimization problems assembled by the methods is shown on Figure \ref{fig:Comparison1} (center left). For adaptive strategies, this number corresponds to the final iteration. The semi-definite approach AM is obviously excluded from this comparison.
Two groups are clearly separated: Adaptive and Local methods on one side, with quasi-linear growth, and Global methods on the other side, with quadratic growth. Let us emphasize that, despite the similar cardinalities, many constraints of the adaptive methods are not local, see Figure \ref{fig:StencilRefinement}. The method $\DConv$ generally uses the least number of constraints, followed by OF$_2$ and then $\Conv$.

\begin{Definition}
\label{def:Defect}
The convexity defect of a discrete map $u : X \to \R$, is the smallest $\ve\geq 0$ such that $u+\ve q \in \Conv(X)$. 
The directional convexity defect of $u$ is the smallest $\ve \geq 0$ such that $u+\ve q \in \DConv(X)$, see Appendix \ref{sec:Directional}.
\end{Definition}

Figure \ref{fig:Comparison1} displays the convexity defect of the discrete solutions produced by the different algorithms, at several resolutions.
This quantity stabilizes at a positive value for the methods OF$_2$ and OF$_3$, which betrays their non-convergence as $n\to \infty$. We expect the convexity defect of the method AM to tend to zero, as the resolution increases, since this method benefits from a convergence guarantee \cite{Aguilera:2008uq}; for practical resolutions, it remains rather high however. 
Other methods, except $\DConv$, have a convexity defect several orders of magnitude smaller, and which only reflects the numerical precision of the optimizer (some of the prescribed linear constraints are slightly violated by the optimizer's output). 
Finally, the method $\DConv$ has a special status since it often exhibits a large convexity defect, but its \emph{directional} convexity defect vanishes (up to numerical precision).

We attempt on figure \ref{fig:BadHessian} to extract, with the different numerical methods, the regions of economical interest: potential customers excluded 
from the trade $\{U=0\}$, and of customers subject to bunching $\{\det (\Hessian U)=0\}$. 
While the features extracted from the method $\Conv$ are (hopefully) convincing, the coordinate bias of the method OF$_2$ is apparent, whereas the method AM does not recover the predicted triangular shape of the set of excluded customers \cite{Rochet:1998uj}. 
The other methods $\DConv$, CLRM, EM, (not shown)  perform similarly to $\Conv$; the method OF$_3$ (not shown) works slightly better than OF$_2$, but still suffers from coordinate bias. The method Ob$_3$ (not shown) seems severely inaccurate%
\footnote{%
The methods (Ob$_k$)$_{k\geq 1}$ are closely related to our approach since they produce outputs with zero convexity defect (up to numerical precision), and the number of linear constraints only grows linearly with the domain cardinality: $\cO(k^2 N)$, with $N:= \#(X)$. We suspect that better results could be obtained with these methods by selecting adaptively and locally the integer $k$.
}%
: indeed the hessian matrix condition number with Ob$_k$, $k \geq 1$, cannot drop below $\approx 1/k^2$, see \cite{Oberman:2011wy}, which is incompatible with the bunching phenomenon, see the solution gradients on Figure \ref{fig:DTGC}.


For each method we compute exactly the monopolist profit \eqref{eq:TotalProfit}, associated with the largest global map $U \in \Conv(\Omega)$ satisfying $U \leq u$ on $X$, where $u\in \cF(X)$ is the method's discrete output. 
It is compared on Figure \ref{fig:Comparison1} with the best possible profit (which is not known, but was extrapolated from the numerical results). 
Convergence rate is numerically estimated to $n^{-1.1}$ for all methods%
\footnote{%
The (presumed) non-convergence of the methods OF$_2 $ and OF$_3$ is not visible in this graph.
}
except (i) the semi-definite approach AM for which we find $n^{-0.75}$, and (ii) the method Ob$_3$, for which energy does not seem to decrease. 

In terms of computation time%
\footnote{%
Experiments conducted on a 2.7 GHz Core i7 (quad-core) laptop, equipped with 16 GB of RAM.
}%
, three groups of methods can be distinguished. Global methods suffer from a huge memory cost in addition to their long run times. 
Methods using a (quasi)-linear number of constraints have comparable run times, thanks to the limited number of stencil refinement steps of the adaptive ones (their computation time might be further reduced by the use of appropriate hot starts for the consecutive subproblems). 
Finally, the semi-definite programming based method AM is surprisingly fast%
\footnote{%
The method AM, implemented with Mosek's conic optimizer, takes only 2.5s 
to solve the product bundles variant of the monopolist problem on a $64\times 64$ grid, with a uniform density of consumers on $[0,1]^2$.  This contrasts with the figure, 751s, reported in \cite{Aguilera:2008uq} in the same setting but with a different optimizer.
}, although this is at the expense of accuracy, see above.
For $n=100$, the method CLRM would use $27\times 10^6$ linear constraints, which with our equipment simply do not fit in memory. The proposed method $\Conv$ selects in $5$ refinement steps a subset containing $\approx 0.4 \%$ of these constraints ($100\times 10^3$), and which is by construction guaranteed to include all the active ones; it completes in $6$ minutes on a standard laptop. 

In summary, adaptive methods combine the accuracy and convergence guarantees of methods based on global constraints, with the speed and low memory usage of those based on local constraints.


%


%
%




%



\section*{Conclusion and perspectives}

We in this paper introduced a new hierarchy of discrete spaces, used to adaptively solve optimization problems posed on the cone of convex functions. The comparison with existing hierarchies of spaces, such as wavelets or finite element spaces on adaptively refined triangulations, is striking by its similarities as much as by its differences. 
The cones $\Conv(\cV)$ (resp.\ adaptive wavelet or finite element spaces) are defined through linear inequalities (resp.\ bases), which become increasingly global (resp.\ local) as the adaptation loop proceeds. Future directions of research include improving the algorithmic guarantees, developing more applications of the method such as optimal transport, and generalizing the constructed cones of discrete convex functions to unstructured or three dimensional point sets.

\paragraph{Acknowledgement.}
The author thanks Pr Ekeland and Pr Rochet for introducing him to the monopolist problem, and the Mosek team for their free release policy for public research.

\bibliographystyle{abbrv}
\bibliography{SelectedPapers} 
\appendix

\section{Directional convexity}
\label{sec:Directional}

We introduce and discuss a weak notion of discrete convexity, which involves slightly fewer linear constraints than \eqref{def:ConvX} and seems sufficient to obtain convincing numerical results, see \S \ref{sec:Comparison}.
\begin{Definition}
\label{def:DConv}
We denote by $\DConv(X)$ the collection of elements in $\cF(X)$ on which all the linear forms $S_x^e$, $x \in X$, $e\in \Z^2$ irreducible, supported on $X$, take non-negative values.
\end{Definition}

Some elements of $\DConv(X)$ cannot be extended into global convex maps on $\Hull(X)$. Their existence follows from the second point of Theorem \ref{th:AllConstraints} (minimality of the collection of constraints $S_x^e$, $T_x^e$), but for completeness we give (without proof) a concrete example.
\begin{Proposition}
Let $u \in \cF(\Z^2)$ be defined by $u(1,1) = 1$, $u(-1,0)=u(0,-1) = -1$, and $u(x) := 2\|x\|^2$ for other $x \in \Z^2$. Then  for all $x \in \Z^2$, and all irreducible $e \in \Z^2$, one has $S_x^e (u) \geq 1$, and $T_x^e(u) \geq 2$ if $\|e\|>1$, with the exception $T_0^{(1,1)}(u) = -1$.
\end{Proposition}

%

Elements of $\DConv(\Z^2)$ are nevertheless ``almost'' convex, in the sense that their restriction to a coarsened grid is convex. 
\begin{Proposition}
If $u \in \DConv(X)$, then $u_{|X'} \in \Conv(X')$, with $X' : =X \cap 2 \Z$.
\end{Proposition}
\begin{proof}
Let $x \in X$, and let $e \in \Z^2$, $\|e\|>1$, be irreducible and of parents $f,g$. Assuming that $x+2 e, x-2 f, x-2 g \in X$, and observing that $x \pm e \in X$ by convexity, we obtain
\begin{equation*}
u(x+2 e)+u(x-2 f) +u(x- 2 g) - 3 u(x) = 2 S_x^e(u) + S_{x+e}^e(u) + S_{x-e}^{f-g}(u) \geq 0.
\end{equation*}
Likewise $u(x+2 e)- 2 u(x) + u(x-2 e) = S_{x-e}^e(u)+2 S_x^e(u)+S_{x+e}^e(u) \geq 0$. 
\end{proof}

The cone $\DConv(X)$ of directionally convex functions admits, just like $\Conv(X)$, a hierarchy of sub-cones $\DConv(\cV)$ associated to stencils. 
\begin{Definition}
Let $\cV$ be a family of stencils on $X$, and let $u \in \cF(X)$. 
The cone $\DConv(\cV)$ is defined by the non-negativity of the following linear forms: for all $x \in X$
\begin{itemize} 
\item For all $e \in \cV(x)$, the linear form $S_x^e$, if supported on $X$.
\item For all $e \in \hat \cV(x)$, the linear form 
$
H_x^e :=  P_x^e + P_x^{-e}
$, if supported on $X$.
\end{itemize}
\end{Definition}

\begin{Proposition}
\begin{itemize}
\item 
For any stencils $\cV, \cV'$, one has $\DConv(\cV) \subset \DConv(X)$,  $\DConv(\cV) \cap \DConv(\cV') = \DConv(\cV \cap \cV')$, and $\DConv(\cV) \cup \DConv(\cV') \subset \DConv(\cV \cup \cV')$.
\item 
For any stencils $\cV$, one has 
\begin{equation}
\label{eq:DConvH}
\DConv(\cV) = \{u \in \DConv(X); \, H_x^e(u) \geq 0 \text{ for all } x \in X, \, e \in \cV_{\max}(x) \sm \cV(x)\}.
\end{equation}
\end{itemize}
\end{Proposition}

\begin{proof}
As observed in \S \ref{sec:CombiningIntersecting}, the second point of this proposition implies the first one. We denote by $\mP(\cV)$ the identity \eqref{eq:DConvH}, and prove it by decreasing induction over $\#(\cV)$. Since $\mP(\cV_{\max})$ clearly holds, we consider stencils $\cV \subsetneq \cV_{\max}$.

Let $x \in X$, $e \in \cV_{\max}(x) \sm \cV(x)$, be such that $\|e\|$ is minimal. Similarly to Proposition \ref{prop:ConvLComplement} we find that $e$ belongs to the set $\hat \cV(x)$ of candidates for refinement at $x$, and define stencils $\cV'$ by $\cV(x) := \cV(x) \cup \{e\}$, and $\cV'(y) := \cV(y)$ for $y \neq x$. The cones $\DConv(\cV)$ and $\DConv(\cV')$ are defined by a common collection of constraints, with the addition respectively of $H_x^e$ for $\DConv(\cV)$, and $S_x^e$, $H_x^{e+f}$, $H_x^{e+g}$ for $\DConv(\cV')$. Expressing the latter linear forms as combinations of those defining $\DConv(\cV)$
\begin{equation*}
S_x^e = H_x^e + S_x^f + S_x^g,\quad H_x^{e+f} = H_x^e + S_{x+e}^f + S_{x-e}^f, \quad H_x^{e+g} = H_x^e + S_{x+e}^g + S_{x-e}^g,
\end{equation*}
and observing that $\mP(\cV')$ holds by induction, we conclude the proof of $\mP(\cV)$.
\end{proof}

\end{document}